\documentclass[3p,times]{elsarticle}

\newcommand{\rred}[1]{\textcolor{red}{#1}}

\newcommand{\ema}{\textcolor{orange}}

\usepackage{array}
\usepackage{subcaption}
\usepackage[dvipsnames]{xcolor}
\usepackage{tabularx}
\usepackage{graphicx,wrapfig,lipsum}
\usepackage{ifplatform}
\usepackage{preview}
\usepackage{environ}
\usepackage{graphics}
\usepackage{epstopdf}
\usepackage{amsmath}
\usepackage{amssymb}
\usepackage{amsfonts}	
\usepackage{moreverb}
\usepackage{dsfont}
\usepackage{grffile}
\usepackage{bm}
\usepackage{multirow}
\usepackage{frcursive}
\usepackage{mathrsfs}
\usepackage{nicefrac}
\usepackage[toc]{appendix}

\usepackage[textsize=footnotesize backgroundcolor=yellow!70, bordercolor=orange]{todonotes}

\newcommand\BibTeX{{\rmfamily B\kern-.05em \textsc{i\kern-.025em b}\kern-.08em
T\kern-.1667em\lower.7ex\hbox{E}\kern-.125emX}}

\usepackage{lmodern}  
\usepackage{xcolor}   
\usepackage{listings}
\usepackage{hyperref} 
\hypersetup{
    colorlinks=true,                          
    linkcolor=blue, 
    citecolor=red, 
    urlcolor=blue  } 
\usepackage[square,numbers]{natbib}
\RequirePackage[hyperpageref]{backref}
\backreffrench
\renewcommand*{\backref}[1]{}  
\renewcommand*{\backrefalt}[4]{
  \ifcase #1 %
  \relax
  \or
  (Cited page~#2.)%
  \else
  (Cited page~#2.)%
  \fi}

\newcommand{\Frac}{\displaystyle\frac}

\newcommand{\U}{\mathbf{U}}
\newcommand{\W}{\mathbf{W}}

\newcommand{\x}{\mathbf{x}}
\newcommand{\n}{\mathbf{n}}
\renewcommand{\u}{\mathbf{u}}

\def\bcdot{{*}}
\newcommand{\Z}{\mathbb{Z}}
\newcommand{\ha}{\frac{1}{2}}

\newcommand{\bi}{\mathbf{i}}
\newcommand{\bj}{\mathbf{j}}
\newcommand{\bhalf}{\mathbf{\ha}}
\newcommand{\bef}{\mathbf{e}_1}
\newcommand{\bes}{\mathbf{e}_2}
\newcommand{\bzero}{\mathbf{0}}
\newcommand{\bone}{\mathbf{1}}

\newtheorem{rem}{Remark}

\def\be{\begin{equation}}
\def\ee{\end{equation}}
\newcommand{\ei}[0]{\end{itemize}}
\newcommand{\beann}[0]{\begin{eqnarray*}}
\newcommand{\eeann}[0]{\end{eqnarray*}}
\def\bea{\begin{eqnarray}}
\def\eea{\end{eqnarray}}
\def\ba{\begin{array}{l}\displaystyle}
\def\ea{\end{array}}

\newcommand{\apriori}{\textit{a priori} }
\newcommand{\aposteriori}{\textit{a posteriori} }

\newfont{\numerikEleven}{ecrm1000}
\newfont{\numerikTen}{cmss10}
\newfont{\numerikNine}{cmss9}
\newfont{\numerikEight}{cmss8}
\newfont{\numerikSeven}{cmss7}
\newfont{\numerikSix}{cmss6}
\usepackage[normalem]{ulem}

\journal{Journal of Computational Physics}

\begin{document} 
\begin{frontmatter}
\title{An almost fail-safe a-posteriori limited high-order CAT scheme} 

\author[catania]{Emanuele Macca$^{*}$}
\ead{emanuele.macca@unict.it}
\cortext[cor1]{Corresponding author}

\author[bordeaux]{Rapha{\"e}l Loub{\`e}re}
\ead{raphael.loubere@math.u-bordeaux.fr}

\author[malaga]{Carlos Par{\'e}s}
\ead{pares@uma.es}

\author[catania]{Giovanni Russo}
\ead{giovanni.russo1@unict.it}


\address[bordeaux]{Universit{\'e} de Bordeaux, CNRS, Bordeaux INP, IMB, UMR 5251, F-33400 Talence, France}
\address[malaga]{Departamento de An{\'a}lisis Matem{\'a}tico, Estad{\'i}stica e Investigaci{\'o}n Operativa, y Matem{\'a}tica aplicada, Universidad de M{\'a}laga, Bulevar Louis Pasteur, 31, 29010, Málaga, Spain}
\address[catania]{Dipartimento di Matematica ed Informatica Universit{\'a} di Catania, Viale Andrea Doria 6, 95125, Catania, Italy}

\begin{abstract}
In this paper we blend the high order Compact Approximate Taylor (CAT) numerical schemes with an \textit{a posteriori} Multi-dimensional Optimal Order Detection (MOOD) paradigm to solve hyperbolic systems of conservation laws in 2D. The resulting scheme presents high accuracy on smooth solutions, essentially non-oscillatory behavior on irregular ones, and, almost fail-safe property concerning positivity issues. 
The numerical results on a set of sanity test cases and demanding ones are presented assessing the appropriate behavior of the CAT-MOOD scheme.
\end{abstract}

\begin{keyword}
 High-order scheme \sep
 CAT \sep
 MOOD \sep
 HLL/HLLC \sep Rusanov \sep
 Hyperbolic system of conservation laws \sep
 Hydrodynamics.
\end{keyword}
\end{frontmatter}


\section{Introduction} \label{sec:introduction}
Peter Lax and Burton Wendroff have presented their seminal finite difference numerical method more than 80 years ago in \cite{LaxWendroff}.
This scheme was designed to solve generic hyperbolic systems of conservation laws.
At the core of the Lax-Wendroff (LW) scheme lays the so-called LW procedure which relies on Taylor expanding  the solution in time up to second-order of accuracy, then replacing the time derivative by the space derivative according to the governing equations, and, finally, approximating the space derivatives by finite differences.
This procedure revealed extremely fruitful. This is still used to design modern numerical schemes, for instance the Approximate Taylor methods introduced by Zorio et al. \cite{ZBM2017}, or, the Compact Approximate Taylor (CAT) family \cite{Carrillo-Pares,CPZMR2020,Macca-Pares,TesiPhD}, and others. 

When dealing with discontinuous solutions which may occur for any hyperbolic system of Partial Differential Equations (PDEs), the key point of most of numerical methods is their ability to dissipate appropriately. In other words extra dissipation must be added. The questions of where, when, and how much dissipation is to be added are of paramount importance to assure that the numerical method can handle smooth flows and discontinuous solutions equally well.
Limiting second-order schemes has been achieved with slope or flux limiters relying on maximum principle preservation, or alternative related procedures. 
Beyond second-order accuracy, the limiting is not anymore a well-agreed subject of research. The most known technique is presumably the ENO/WENO procedure \cite{Shu1} for finite volume (FV) or finite differences (FD) schemes. Most of the limiting techniques rely on blending the first-order scheme/flux/reconstruction with a high-order one, using some \textit{a priori} sensor to determine where this blending would be appropriate.
The limiting entirely depends on the 
quality of the \textit{a priori} sensor which must determine "where to act?" and the amount of blending, i.e.\ "how much dissipation?" to ensure that the numerical solution is physically and numerically acceptable.
Based on this philosophy, high-order CAT schemes have been supplemented with an automatic limiting procedure in \cite{Macca-Pares}, called ACAT.
However this procedure suffers from the same drawback: it relies on \textit{a priori} sensors which are difficult to design. \\
Contrarily, in this work we operate a shift and will couple the high-order CAT schemes with the \textit{a posteriori} MOOD limiting procedure \cite{CDL1,CDL2,CDL3}. 
The philosophical statement in MOOD is that it is easier to observe the inappropriate consequences of using a high-order scheme rather than predicting them with \textit{a priori} sensors. Consequently within a MOOD loop the unlimited high-order explicit scheme is used for the current time-step to produce a candidate numerical solution which is further tested against detection criteria stating if some cells are invalid. While the valid ones are kept, the invalid ones are recomputed with a low-accurate but more robust scheme, possibly a second-order limited scheme, or, more drastically a robust first order one.
The goal of this paper is to provide a proof of concept for the design and validation of a CAT-MOOD scheme for systems of conservation laws.

The rest of this paper is organized as follows. 
The second section introduces the governing equations.
Next, the Lax-Wendroff procedure and the CAT schemes are presented.
We re-derive the second-, fourth- and sixth-order accurate versions of CAT schemes
and the generic CAT2P scheme, along with their adaptive limiter.
The fourth section presents the blending of CAT schemes with the \textit{a posteriori} MOOD procedure to replace the adaptive limiter.
Numerical results are gathered in the fifth section to assess the good behavior of the CAT-MOOD sixth order scheme on smooth or discontinuous solutions.
Conclusions and perspectives are finally drawn.

%
\section{Governing equations} \label{sec:equations}

\subsection{1D linear and non-linear scalar conservation laws}
To simplify the description of the numerical methods we also consider the non-linear scalar conservation law on the $Oxt$-Cartesian frame
\begin{equation}
    \label{sec:CAT_gov_equ}
    u_t + \partial_x f(u) = 0,
\end{equation} 
where $u=u(x,t):\mathbb{R}\times\mathbb{R}^+\rightarrow\mathbb{R}$ denotes the scalar variable, and $f(u)=f(u(x,t))$ the non-linear flux depending on $u$.
$u(x,0)=u_0(x)$ denotes the IC, while BC are prescribed depending on the test case; for instance periodic ones, Dirichlet or Neumann ones. 

~(\ref{sec:CAT_gov_equ}) represents the generic model of non-linear scalar equation, the simplest one being probably Burgers' equation for which the flux is given by: $f(u)=u^2/2$.
~(\ref{sec:CAT_gov_equ}) also represents the generic model of linear scalar advection equation if $f(u)=a u$ with $a\in \mathbb{R}$ being the advection velocity. \\ 
In the following we denote the partial derivative in time and space with under-script letters as 
$u_t \equiv \partial_t u$ and $u_x \equiv \partial_x u$.

\subsection{2D gas-dynamics system of conservation laws}
In this paper we focus on hyperbolic systems of conservation laws (Partial Differential Equations, PDEs) in 1D and 2D of the form
\bea \label{eq:systemPDEs}
    \partial_t \U + \nabla \cdot \mathbb{F}(\U) = 0,
\eea
where $t\in \mathbb{R}^+$ represents the time variable, $\x=(x,y)\in \mathbb{R}^2$ the space variable in $2$ dimensions.
$\U=\U(\x,t)$ is the vector of conserved variables while 
$\mathbb{F}(\U)=(\mathbf{F}(\U(\x,t)),\mathbf{G}(\U(\x,t)) )^t$ is the flux vector.
$\nabla\cdot$ is the divergence operator which allows us to rewrite
(\ref{eq:systemPDEs}) as
\bea \label{eq:systemPDEs_bis}
    \partial_t \U + \partial_x \mathbf{F}(\U) + \partial_y \mathbf{G}(\U) = \mathbf{0}.
\eea
In this work we mainly focus on the gas-dynamics system of PDEs (Euler equations) where $\U=(\rho, \rho u, \rho v, \rho e)^\top$ with $\rho$ the density, $\u=(u,v)$ the velocity vector, $e=\varepsilon + \frac12 \| \u \|^2$ the total energy per unit mass, and $\varepsilon$ the internal one. 
The flux tensor is given by 
\bea
\mathbb{F}(\U)=\left( \begin{array}{cc} \mathbf{F}(\U) & \mathbf{G}(\U) \end{array} \right), 
\quad \text{with} \quad
\mathbf{F}(\U)= \left( \begin{array}{c}
\rho u \\
\rho u^2 +p  \\
\rho uv  \\
(\rho e +p )u  \\
\end{array} \right), 
\quad 
\mathbf{G}(\U)= \left( \begin{array}{c}
\rho v \\
\rho u v \\
 \rho v^2  +p  \\
(\rho e +p )v \\
\end{array} \right).
\eea
The system is closed by an Equation Of State (EOS) which specifies the value of the pressure $p$ as a function of two thermodynamics variables, for instance of the form $p=p(\rho,\varepsilon)$. For a polytropic gas we have $p=(\gamma-1)\rho \varepsilon$ with $\gamma$ the polytropic constant characterizing the type of gas considered. The sound-speed is given by $a^2=\gamma p /\rho$.
The equations in (\ref{eq:systemPDEs}) express the conservation of mass, momentum and total energy.
An entropy inequality is supplemented to the system of PDEs to ensure that the weak solutions are the entropic ones.
This system is hyperbolic with eigenvalues $\lambda^-=\u\cdot\n - a$, $\lambda^0=a$ (multiplicity 2), $\lambda^+=\u\cdot\n + a$ where $\n$ is a generic unit vector indicating a direction in 2D. 
It is well known that the physical states all belong to 
\bea \label{eq:physical_states}
    \mathcal{A} = \left\{ \U\in \mathbb{R}^4, \; \text{such that} \; \rho>0, \; p>0 \right\}.
\eea
The primitive variables are the components of vector $\W=( \rho, u,v, p )^\top$ and are computed from the conservative ones as 
\bea
    u = (\rho u)/\rho, \qquad
    v = (\rho v)/\rho, \qquad
    p = (\gamma - 1)\left( (\rho E) - \frac12 \left((\rho u)^2+(\rho v)^2 \right)/\rho \right).
\eea
System (\ref{eq:systemPDEs}) is further equipped with Initial Conditions (IC) and Boundary Conditions (BC). 
This system of PDEs is the target one, but a simpler one is also considered in the next section for the sake of simplicity.

\subsection{Mesh}
In this article we consider logical rectangular meshes in 1D and 2D. 
The time domain $\mathcal{T}=[0,T]$ with final time $T>0$ is split into time intervals $[t^n,t^{n+1}]$, $n\in \mathbb{N}$, and time-steps $\Delta t=t^{n+1}-t^n$ subject to a CFL (Courant-Friedrichs-Lewy) condition. 
The computational domain denoted $\Omega$ is a segment/rectangle in 1D/2D. 
In 1D,  $\Omega$ is paved with $N_x$ cells. The generic cell is denoted $\omega_i$ and indexed by a unique label $1\leq i \leq N_x$. Classically we identify the cell end-points by half indexes so that $\omega_i=[x_{i-1/2},x_{i+1/2}]$ and the cell center is given by $x_i=\frac12(x_{i+1/2}+x_{i-1/2})$.
The size of a cell is given by $\Delta x_i = x_{i+1/2}-x_{i-1/2}$, that we simply denote as $\Delta x$, since we assume, for simplicity, that the grid is uniform.
In 2D, $\Omega$ is paved with $N_x \times N_y$ cells. The generic cell is denoted $\omega_{i,j}$ and indexed by a double label $1\leq i \leq N_x$ and $1\leq j \leq N_y$. 
The four vertices of any cell $\omega_{i,j}$ are $\bm{x}_{i\pm 1/2,j\pm 1/2}=( x_{i\pm 1/2}, y_{j\pm 1/2})$. $x_{i+1/2}$ represents a vertical mesh line, while $y_{j+1/2}$ a horizontal one. 
Accordingly, $\omega_{i,j}=[x_{i-1/2},x_{i+1/2}]\times [y_{j-1/2},y_{j+1/2}]$, and, the cell center is given by $\bm{x}_{i,j}=(x_i,y_j)= \left( \frac{x_{i+1/2}+x_{i-1/2}}{2},\frac{y_{j+1/2}+y_{j-1/2}}{2}\right)$.
The size of a cell is given by $\Delta x_i\times \Delta y_j$ with $\Delta y_j = y_{j+1/2}-y_{j-1/2}$, that we simply denote as $\Delta x$ and $\Delta y$, since we shall adopt a uniform mesh throughout the paper.
We call an interface or face, the intersection between two neighbor cells, that is a point in 1D and a segment in 2D.
The neighbor cells of a generic one are those with a non-empty intersection. A generic cell has two/eight neighbors in 1D/2D on logical rectangular grids. In 2D we make the difference between the four face-neighbors and the four corner-neighbors.
A "stencil" in 1D is a collection of $K>0$ cells surrounding and including the current one, onto which derivatives are approximated.
\begin{figure}[!h]
\centering
\includegraphics[width=0.4\textwidth]{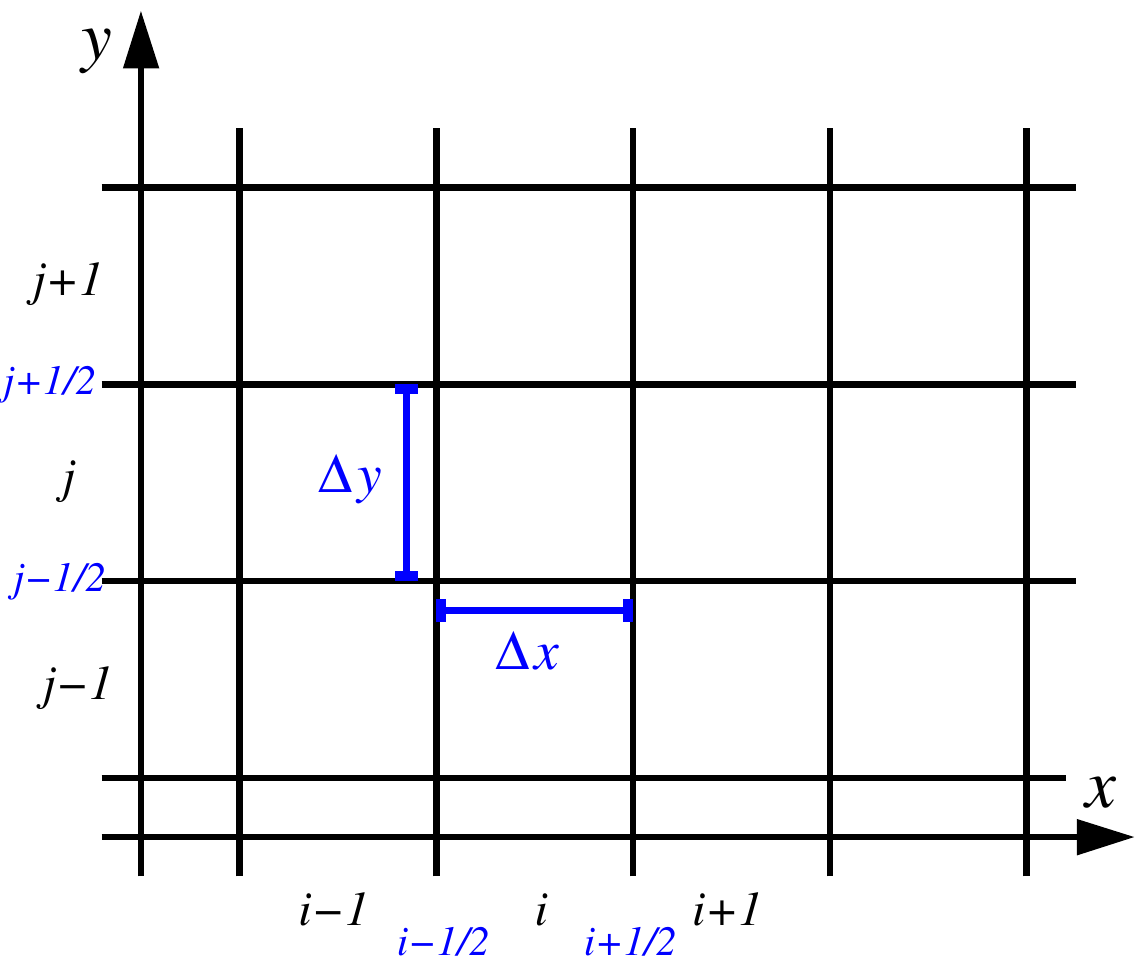}
\includegraphics[width=0.59\textwidth]{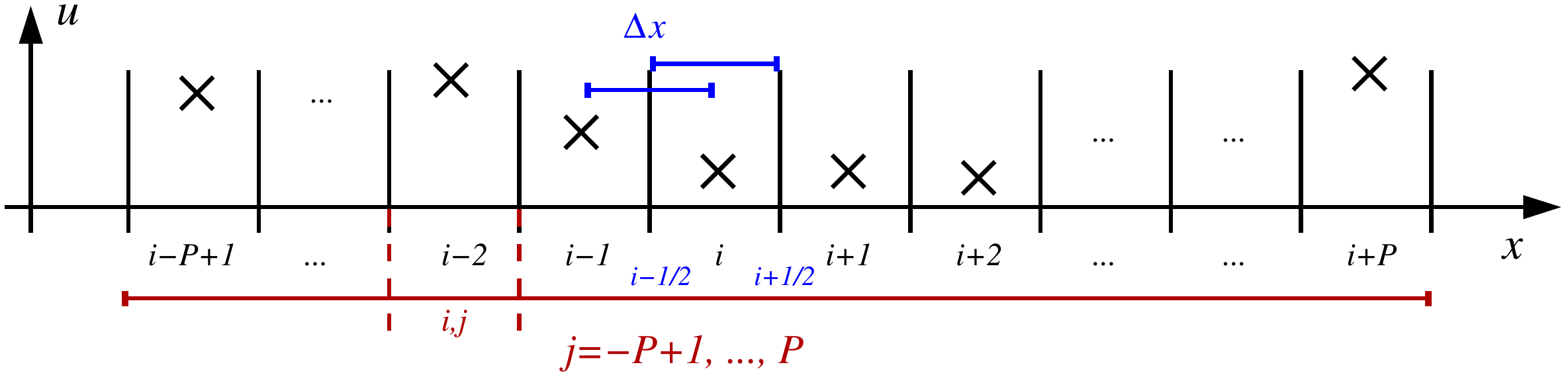}
 \caption{Left: Logical rectangular grid and notation used in this paper.
 Right: illustration of spatial stencil in 1D around cell $i$ with stencils $i,j$ with $-P+1\leq j \leq P$.}
 \label{fig:grid}
\end{figure}

\subsection{Notation}
In this article the scheme description is made mostly in 1D for the sake of clarity.
The following notations refer to
the different type of derivatives or approximations.
\begin{itemize}
    \item $u_{i,j}^{(k)}$ is the $k$-th time derivative of $u$ at time $t^n$ in position $x_{i+j},$ where $i$ refers to the cell and $j$ to the position in the stencil. In general, for a scheme of order $2P,$ $k = 1,\ldots,2P-1$ and $j = -P+1,\ldots,P.$

    \item $f_{i,j}^{(k)}$ is the $k$-th time derivative of $f(u)$ at time $t^n$ in position $x_{i+j}$, likewise for $u$. 
    In general, for a scheme of order $2P,$ $k = 0,\ldots,2P-1$ and $j = -P+1,\ldots,P,$ under the assumption that $f_{i,j}^{(0)} = f(u_{i+j}^{n})$.

    \item $u_{i,j}^{k,n+r}$ is the explicit Taylor expansion of function $u$ in time truncated to order $k$,  centered at time $t^n$ at distance $r\Delta t$ in time and at spatial location $x_{i+j}$. Again $i$ refers to the cell and $j$ to the position in the stencil. In general, for a scheme of order $2P,$ $k = 1,\ldots,2P-1,$ while $j,r = -P+1,\ldots,P.$

    \item $f_{i,j}^{k,n+r}$ refers to $f\Bigl(u_{i,j}^{k,n+r}\Bigr).$

    \item The symbol $*$ will be used to indicate to which index (space or time) the differentiation is applied as illustrated in the equations below for the approximation of space and time derivatives, respectively (see also ~\eqref{eq:interpolA})
    \begin{eqnarray}
            \label{not_1}
            \partial_x^kf(x_i+q\Delta x,t^n) &\approx& A^{k,q}_{p}(f_{i,*}^n, \Delta x)   =   \frac{1}{\Delta x^k} \sum_{j=-p+1}^{p} \gamma^{k,q}_{p,j} f^n_{i+j}, \\
            \label{not_2}
            \partial_t^kf(x_i,t^n) &\approx& A^{k,0}_{p}(f^{n,*}_i, \Delta t)   =   \frac{1}{\Delta t^k} \sum_{r=-p+1}^{p} \gamma^{k,q}_{p,r} f^{n,r}_{i}.
        \end{eqnarray}
\end{itemize}
From now on, since all the formulas are computed at time $t=t^n$, we avoid the extra index $n$ on equations like \eqref{not_1} and \eqref{not_2}.

%
\section{Compact Approximate Taylor (CAT) schemes} \label{sec:CAT} 

The focus of this chapter is the presentation of a family of numerical methods for non-linear systems of conservation law, named Compact Approximate Tayor (CAT) schemes\footnote{For simplicity we will summarize the numerical deduction of the method in the scalar case.}.
This family is based on an approximate Taylor procedure that constitutes a proper generalization of Lax-Wendroff (LW) method, in the sense that it reduces to the standard high-order LW method when the flux is linear. 
In this section we recall the LW and CAT procedures of second and fourth order.

\subsection{Lax Wendroff procedure} \label{ssec:LW}
A scheme of historical as well as practical importance is the celebrated \textbf{Lax Wendroff} scheme introduced by Peter Lax and Burton Wendroff in 1960 in \cite{LaxWendroff}, and \cite{Lax1957,Wendroff,LaxWendroffHigh,Hirsch22007,Toro2009}. It has been the most widely adopted scheme for aeronautical applications, up to the end of the 1980s under various forms.\footnote{Probably the most widely used variant is the two-step Mac Cormack scheme, which has properties similar to LW schemes, and avoids the computation of the second derivative. 
Published originally at a conference in 1969, the paper has been reproduced in \cite{maccormack2002effect}.}

%
The original derivation of Lax and Wendroff was based on a Taylor expansion in time of function $u$ at point $(x_i,t)$ up to second order of accuracy, thus
\begin{equation}
    \label{sec:LW_Taylor2}
    u(x_i,t+\Delta t)=u(x_i,t) + \Delta t \, u_t(x_i,t) + \frac{\Delta t^2}{2}u_{tt}(x_i,t) + O(\Delta t^3),
\end{equation}
where $\Delta t>0$ is a small increment in time.
The numerical scheme is then obtained by neglecting the higher order term in $\Delta t$, using the governing equation to replace time derivatives by spatial ones, and then substituting the obtained space derivatives with their finite difference approximations. 
For the \textit{linear case}, $u_t =-a u_x$ with $f(u)=au$, we obtain
\bea
    \label{eq:lw0}
    u(x_i,t^n+\Delta t)=u(x_i,t^n) - \Delta t \, a u_x(x_i,t^n)  + \frac{\Delta t^2}{2} a^2 u_{xx} (x_i,t) + O(\Delta t^3).
\eea
Using centred finite differences to approximate spatial derivatives, the numerical scheme follows as:
\begin{equation}
    \label{lw}
    u_i^{n+1}=u_i^n - \Delta t \, a \, \frac{u_{i+1}^n-u_{i-1}^n}{2\Delta x} + \frac{\Delta t^2}{2} a^2  \Frac{u_{i+1}^n - 2u_i^n + u_{i-1}^n}{\Delta x^2},
\end{equation}
where 
$u_i^m$ is an approximation of the point value of the solution at position $x_i$ at the time $t^m$. 
A useful alternative formulation written in conservative form yields
\begin{equation}
    \label{lw_conservative}
    u_i^{n+1}=u_i^n - \frac{\Delta t}{\Delta x}\left( F_{i+\ha}^{\rm{LW}}- F_{i-\ha}^{\rm{LW}} \right),
\end{equation}
where the so-called LW numerical flux, $F_{i+\ha}^{\rm{LW}}$, is given by
\bea \label{eq:LW_flux}  
F_{i+\ha}^{\rm{LW}} = \frac{a}{2}\left(u_{i+1}^n + u_i^n\right) - \frac{a^2\Delta t}{2\Delta x}\left(u_{i+1}^n - u_i^n\right).
\eea
%
The \textit{non-linear case}, $u_t =-f_x(u)$, yields
\begin{equation}
    \label{non-lin_lw}
    u_i^{n+1}=u_i^n - \frac{\Delta t}{2\Delta x}\left(f_{i+1}^n-f_{i-1}^n\right) + \frac{\Delta t^2}{2\Delta x^2}\left(\text{A}_{i+\ha}\left(f_{i+1}^n - f_{i}^n\right) - \text{A}_{i-\ha}\left(f_{i}^n - f_{i-1}^n\right)\right),
\end{equation}
where $f_{i+j}^n = f(u_{i+j}^n)$ for $j = -1,0,1$ and $\text{A}$ is an approximation of the derivative of $f$, i.e $A= \partial f/\partial u$. 
Hence, $\text{A}_{i\pm\ha}$ is the approximation derivative of $f$
evaluated at $u_{i+\ha}^n=\ha(u_{i}^n + u_{i\pm1}^n)$, that is $\text{A}_{i\pm\ha}\equiv\text{A}(u_{i+\ha}^n)$, or, the average between the cell-based derivative, that is $\text{A}_{i\pm\ha}\equiv \ha \left( \text{A}(u_i^n)+\text{A}(u_{i\pm1}^n) \right)$. 
Notice that they  depend non-linearly on variable $u$. Moreover they are always evaluated at time $t^n$, so we can omit this time dependency.
The alternative conservative formulation of the LW scheme is expressed as:
\begin{equation}
    \label{non-lin_lw_conservative}
    u_i^{n+1}=u_i^n - \frac{\Delta t}{\Delta x}\left(F_{i+\ha}^{\rm{LW}}- F_{i-\ha}^{\rm{LW}}\right),
\end{equation}
where the numerical flux
\bea \label{eq:flux_NL} 
    F_{i+\ha}^{\rm{LW}} = \underbrace{\frac{1}{2}\left( f_{i+1}^n + f_i^n \right)}_{\text{Physical flux}} \underbrace{- \frac{\Delta t}{2\Delta x}A_{i+\ha}\left(f_{i+1}^n - f_i^n\right)}_{\text{Dissipation}} ,
\eea
is composed of two parts: the average of the physical fluxes at cells $i$ and $i+1$, and, the numerical dissipation.

\subsection{Compact Approximate Taylor (CAT) procedure} \label{ssec:CAT}
The generalized Lax-Wendroff method is used to update the numerical solution:
\bea \label{eq:general_LW}
 u_i^{n+1} = u_i^n + \sum_{k=1}^{2P} \frac{(\Delta t)^k}{k!} u_i^{(k)},
\eea
where we recall that
$u_i^n$ is an approximation of the value of the exact solution $u(x,t)$ at time $t^n$ at position $x_i$ \cite{Hirsch12007}, and, $u_i^{(k)}$ is an approximation of $\partial_t^k u(x_i,t^n)$. 
The $k$-th derivative in time of $u$ is computed with a compact and numerical version of the Cauchy-Kovalesky procedure introduced by Carrillo and Par{\'e}s in \cite{Carrillo-Pares}.  \\ 
The final expression of the $2P$-order CAT method in conservative form is:
\begin{equation}
    \label{CAT2P}
    u_i^{n+1} = u_i^n + \frac{\Delta t}{\Delta x}\left(F_{i-\ha}^P - F_{i+\ha}^P\right).
\end{equation} 
Let us introduce the sets $\mathcal{S}_{i\pm\ha}^P$ of values $u_i^n$ on  
stencils centered around interface $i\pm\ha$ of size $2P$, 
that is
\bea \label{eq:stencil2P} 
\mathcal{S}_{i+\ha}^P = \left\{ u_{i-P+1}^n,\ldots, u_{i}^n, u_{i+1}^n, \ldots u_{i+P}^n \right\}.
\eea
The flux functions $F_{i\pm \ha}^P$ are then computed, respectively, on the sets $\mathcal{S}_{i\pm\ha}^P$, as 
\begin{equation}
    \label{FP}
    F_{i+\ha}^P = \sum_{k=1}^{2P}\frac{\Delta t^{k-1}}{k!}f^{(k-1)}_{i+\ha}, 
\end{equation}
and 
\bea \label{eq:interpolA} 
f^{(k-1)}_{i+\ha} = \mathcal{A}_{P}^{0,\ha}\left(f^{(k-1)}_{i,*},\Delta x\right), 
\quad 
\text{with}
\quad
\mathcal{A}_{P}^{0,\ha}\left(f^{(k-1)}_{i,*},\Delta x\right)  = \sum_{p=-P+1}^P\gamma_{P,p}^{0,\ha} \, f_{i+p}^{(k-1)},
\eea 
where $\mathcal{A}_{P}^{0,\ha}$ is an interpolation formula of order $2P-1$ based on $2P$-point stencil.
In the following we use the index $i$ for the cell global index, $j$ for the local position inside the stencil, $r$ for the Taylor expansion in time, and, $k$/$(k)$ to refer to the $k$-th time step/$(k)$th time derivative. 
For the sake of clarity, we detail the description of the second order ($P=1$) CAT2, in the next sub-section, and CAT4 in the Appendix \ref{sec:CAT4}.

\subsubsection{Second order version -- CAT2}
Let $p\in \mathbb{N}$ denote an integer such that $0 \leq p \leq P$. 
In the case $P=1$ then the relative stencil is simply $\mathcal{S}^1_{i+\ha} = \{u_i^n,u_{i+1}^n\}$ for interface $i+1/2$,
while the flux reconstructions are:
\bea \label{eq:fluxCAT2}
	F_{i+\ha}^1 = f^{(0)}_{i+\ha} + \frac{\Delta t}{2}f^{(1)}_{i+\ha}, \qquad
	F_{i-\ha}^1 = f^{(0)}_{i-\ha} + \frac{\Delta t}{2}f^{(1)}_{i-\ha},
\eea
where $f^{(0)}_{i+\ha} = \ha\left(f_i^n + f_{i+1}^n\right)$ 
is the interpolation of the flux at time $t^n$, while
$f^{(1)}_{i+\ha}$  
is the interpolation of the first time-derivative of the flux for any interface $i+1/2$.
This implies that $\gamma_{1,p}^{0,\ha}=\ha$ in (\ref{eq:interpolA}).
These are computed over stencils $\mathcal{S}_{i\pm\ha}^1$ for
$p\in \{ 0,1 \}$ as
\begin{align}
\nonumber
\mathcal{S}_{i-\ha}^1 &:&
f^{(1)}_{i-\ha} = \frac{f_{i-1}^{(1)}+f_i^{(1)} }{2}, \quad 
&f_{i-1,p}^{(1)} = \frac{f\left( u_{i-1+p}^n + \Delta t \, u_{i-1,p}^{(1)}\right) - f_{i-1+p}^n}{\Delta t}, 
 & u_{i-1,p}^{(1)} = -\frac{f_{i}^n - f_{i-1}^n}{\Delta x},\\
\nonumber
\mathcal{S}_{i+\ha}^1 &:& 
f^{(1)}_{i+\ha} = \frac{f_i^{(1)} + f_{i+1}^{(1)}}{2}, \quad 
&f_{i,p}^{(1)} = \frac{f\left( u_{i+p}^n + \Delta t \, u_{i,p}^{(1)}\right) - f_{i+p}^n}{\Delta t}, 
 &  u_{i,p}^{(1)} = -\frac{f_{i+1}^n - f_{i}^n}{\Delta x}.
\end{align}
Notice that the $\mathcal{S}_{i+\ha}^1$ expressions can be obtained from $\mathcal{S}_{i-\ha}^1$ by replacing $i$ by $i+1$.

Finally, the expanded form of the fluxes (\ref{eq:fluxCAT2}) is given by
\begin{align}
\label{F2-}
    F_{i-\ha}^1 &= \frac{1}{4}\left(f^n_{i-1} + f^n_i + f\left(u_{i-1}^n + \Delta t \, u_{i-1,0}^{(1)}\right) + f\left(u_{i}^n + \Delta t \, u_{i-1,1}^{(1)}\right) \right),\\
    \label{F2+}
    F_{i+\ha}^1 &= \frac{1}{4}\left(f^n_i + f^n_{i+1} + f\left(u_{i}^n + \Delta t \, u_{i,0}^{(1)}\right) + f\left(u_{i+1}^n + \Delta t \, u_{i,1}^{(1)}\right) \right),
\end{align}
and the solution is obtained by substituting these fluxes in formula (\ref{CAT2P}).\\ 
The computation of the interfacial flux $F_{i+\ha}^1$ can be recast into the algorithm:
\begin{itemize}
    \item[Step 1:] Compute $f_{i+\ha}^{(0)}$ adopting an interpolation formula over stencil $\mathcal{S}_{i+\ha}^1$ at time $t^n$.
    \item[Step 2:] Compute the first derivatives $u^{(1)}_{i,p}$ in time through the numerical compact Cauchy-Kovalesky procedure, using  $\partial_t u = -\partial_x f $, with data at time $t^n$. 
    \item[Step 3:] Compute the truncated Taylor expansions: $u^{1,n+1}_{i,p} = u^n_{i+p} + \Delta t\, u^{(1)}_{i,p}$ for $p=0$ and $1$.
    \item[Step 4:] Compute the first time derivatives of the flux using the first difference formulas: 
    \[
        f^{(1)}_{i,p} = \Frac{f \left(u^{1,n+1}_{i,p}\right) - f^n_{i+p}}{\Delta t}.
    \]
    \item[Step 5:] Compute $f_{i+\ha}^{(1)}$ through $f^{(1)}_{i,j}$ adopting an interpolation formula on stencil $\mathcal{S}_{i+\ha}^1;$
    \item[Step 6:] Compute $F_{i+\ha}^1$ as a Taylor expansion: $F_{i+\ha}^1 = f^{(0)}_{i+\ha} + \Frac{\Delta t}{2}f^{(1)}_{i+\ha}$ with (\ref{eq:fluxCAT2}).
\end{itemize}   

\subsubsection{CAT2P}
Generically, the $2P$ CAT scheme follows formulation \eqref{CAT2P} with the interface fluxes given by \eqref{FP}.
The expression of the right numerical flux of order $2P$ is obtained with formula:
\begin{equation}
    \label{FP+}
    F^P_{i+\ha} = \sum_{k=1}^{2P}\frac{\Delta t^{k-1}}{k!}\mathcal{A}^{0,\ha}_{P}\left(f_{i,*}^{(k-1)},\Delta x\right) = \sum_{k=1}^{2P}\frac{\Delta t^{k-1}}{k!}\sum_{j=-P+1}^P \gamma^{0,\ha}_{P,j}f_{i,j}^{(k-1)}, 
\end{equation}
where the high order time derivatives of the flux are computed following and extending the iterative algorithm presented for CAT2 in the previous section (see also \cite{Carrillo-Pares, CPZMR2020, Macca-Pares} for more details):
\begin{enumerate}
    \item Define $f^{(0)}_{i,j} := f(u_{i+j}^n)$ for all $j = -P+1,\ldots,P;$
    \item For every $k = 1,\ldots,2P-1$:
    \begin{enumerate}
        \item Compute the $k$-th derivative of $u$ at time step $t^n$ for each position $x_{i+j}$ with $j = -P+1,\ldots,P$ through the numerical compact version of the Cauchy-Kovalesky identity \eqref{CK}  as:
        $$u_{i,j}^{(k)} = -\mathcal{A}^{1,j}_{P}\left(f_{i,*}^{(k-1)},\Delta x\right); $$
        \item Compute the Taylor expansion of $u$ in time truncated to term $k$ for all positions $x_{i+j}$ with $j = -P+1,\ldots, P$ at time $t^{n+r}$ with $r = -P+1,\ldots,P$ as:
        $$u(x_{i+j},t^{n+r}) \approx u_{i,j}^{k,n+r} = u^n_{i+j} +\sum_{m=1}^k\frac{(r\Delta t)^m}{m!} u^{(m)}_{i,j}; $$
        \item Compute the $k-$th time derivative of flux for each position $x_{i+j}$ with $j = -P+1,\ldots, P$ at time $t^{n}$ as:
        $$f^{(k)}_{i,j} = \mathcal{A}^{k,j}_{P}\left(f_{i,j}^{k,*},\Delta t\right), $$ where $f_{i,j}^{k,*}$ means that we are applying the $\mathcal{A}$ operator in time and in particular we apply the differentiation formula to the set of flux approximations 
        $$ f_{i,j}^{k,n-P+1},\ldots,f_{i,j}^{k,n+P},$$ 
        in which $f_{i,j}^{k,n+r}=f\left(u_{i,j}^{k,n+r}\right)$ for all $j,r =-P+1,\ldots,P$.
    \end{enumerate}
\end{enumerate}

\begin{rem}
    Observe that the computation of the numerical flux  $F_{i+1/2}^{P}$ requires the approximation of $u$ at the nodes of a space-time grid of size $2P \times 2P$, represented by $u_{i,j}^{k,n+r}$, for $-P+1 \leq j,r =\leq P$ (see Figure~\ref{1D_grid}). 
    The approximations of the solution $u$ at successive times $(n-P+1)\Delta t$, \dots, $(n-1)\Delta t$ are different from the ones  already computed in the previous time steps $t^{n-P},\ldots, t^{n-1}$ which are $u_{i+j}^{n-P}$, \ldots, $u_{i+j}^{n-1}$.  
    In other words, the discretization in time is not based on a multi-step method but on a one-step one. In fact, it can be re-interpreted as a Runge-Kutta method whose stages are $u^{n+r}_{i,j}$, $r = -P+1, \dots, P$, see \cite{Carrillo-Pares,CPZMR2020,Macca-Pares}.
\end{rem}
\begin{rem} \label{rem:local}
    These approximations are local. Indeed, suppose that $i_1+j_1 = i_2+j_2 = \ell,$ i.e. $x_{\ell}>0$ belongs to $\mathcal{S}^P_{i_1+1/2}$ and $\mathcal{S}^P_{i_2+1/2}$ with local coordinates $j_1$ and $j_2$ respectively. Then, $f^{(k)}_{i_1,j_1}$ and $f^{(k)}_{i_2,j_2}$ are, in general, two different approximations of $\partial_t^k f(u)(x_{\ell},t^n).$
\end{rem}

\subsubsection{Computational complexity}\label{sec:Comput_Comp}
In this paragraph we estimate the  computational complexity for the CAT$2P$ scheme. 
The details about the complexity of the algorithm can be found in Appendix~\ref{sec:comp_compl}.
As a summary, the operation count per cell per time step for CAT$2P,$ $P>1,$ applied to the scalar case is $3.5(2P)^3 - 1.5(2P)^2 + (2P)$ flop plus $(2P)^3 - 2(2P)^2 + (2P) + 1$ function evaluations. 
For CAT2 this gives $24$ flop and $3$ function evaluations, while for CAT4 we get $204$ flop and $37$ function evaluations.
%
%
%
%
\subsection{Adaptive limiter for CAT schemes - ACAT schemes} \label{ssec:beyondCAT}
Although the Compact Approximate Taylor (CAT) schemes are linearly stable in the $L^2$-sense under the usual CFL-$1$ condition, they may produce bounded oscillations close to the discontinuity of the solution. 
To avoid these spurious phenomena, an Adaptive \apriori shock-capturing technique has been developed in \cite{CPZMR2020} and called ACAT. There, the order of the method is locally adapted to the smoothness of the numerical solution by means of indicators which check the regularity of the data for each temporal step. More specifically,  once the approximations  of the solution $u$ at time $t^n$ have been computed, the stencil of the data adopted to actually compute the right flux $F_{i+1/2}^P$ are set to belong to 
 $$  \mathcal{S}_p = \{ u_{i-p +1}^n, \dots, u_{i+p}^n\}, \quad p = 1, \dots, P.  $$
 The selected stencil is the one with maximal length among those in which the solution at time $t^n$ is 'smooth'.
 The smoothness is assessed according to some smoothness indicators: $\psi^p_{i+1/2}$, for any $p = 1, \ldots, P$, which are defined as:
\begin{equation}\label{condiconesS0}
\psi^p_{i+1/2} \approx \left\{
\begin{array}{cl}
1 & \mbox { if $u$ is 'smooth' in $\mathcal{S}_p $,}\\
0 & \mbox {otherwise.} 
\end{array}\right.
\end{equation} 
For this strategy one needs to define a robust first-order flux reconstruction, for instance  Rusanov-, HLL- or HLLC-based flux reconstruction \cite{Toro2009,LeVeque2002book}.
Next, one can employ a TVD flux-limiter for the second-order flux reconstruction to be combined with CAT2, such as, for instance minmod, van Alabada or superbee \cite{LeVeque2002book}.

\paragraph{ACAT2}
The expression of the ACAT2 numerical method (for $P=1$) based on a flux limiter (see \cite{LeVeque2002book,LeVeque2007book,Toro2009}) is given by
\begin{equation}\label{ACAT2meth}
     u_i^{n+1} = u_i^n + \frac{\Delta t}{\Delta x}\left(F^{*}_{i-1/2} - F^{*}_{i+1/2} \right) ,
\end{equation}
where the fluxes $F^{*}_{i+1/2}$ are blended as   
\begin{eqnarray}\label{ACAT2_flux} 
    {F}^*_{i\pm1/2} & = & \varphi^1_{i\pm1/2} \, {F}^1_{i\pm1/2} +(1-\varphi^1_{i\pm1/2}) \, {F}^{\text{low}}_{i\pm1/2}.
\end{eqnarray}
${F}^1_{i\pm1/2}$ is the CAT2 flux given by \eqref{F2+}-\eqref{F2-}, while $F^{\text{low}}_{i\pm1/2}$ is the first order flux reconstruction, and, $\varphi^1_{i\pm 1/2}$ is a switch computed by a flux limiter which verifies
\begin{equation} \nonumber
    \varphi^1_{i-1/2} \approx \begin{cases}
        1 &\text{if $\{ u_{i-2}^n, \dots, u_{i+1}^n \}$ is 'smooth'}, \\
        0 & \text{otherwise,}
    \end{cases}
    \quad
    \varphi^1_{i+1/2} \approx \begin{cases}
        1 &\text{if $\{ u_{i-1}^n, \dots, u_{i+2}^n \}$ is 'smooth'}, \\
        0 & \text{otherwise.}
    \end{cases}
\end{equation}
For scalar problems, standard flux limiter functions, $\varphi^1(r)$, such as minmod, superbee, van Leer \cite{Sweby,Kemm2010}, may  be used:
\begin{equation}\label{indicador1}
    \varphi_{i+1/2}^1  =  \varphi^1(r_{i+1/2}), 
\end{equation}
where
\begin{equation}\label{wavefuction}
        r_{i+1/2}=
         \left\{
        \begin{array}{cl}
            \displaystyle r^L_{i+1/2}= 
            \frac{u^n_i-u^n_{i-1}}{u^n_{i+1}-u^n_{i}} & \mbox {if } a_{i+1/2}  >0, \\
            \displaystyle  r^R_{i+1/2}= 
            \frac{u^n_{i+2}-u^n_{i+1}}{u^n_{i+1}-u^n_{i}} & \mbox {if } a_{i+1/2}  \leq 0,
        \end{array}\right.
        \quad 
        a_{i+1/2} = \begin{cases}  \displaystyle \frac{f(u^n_{i+1}) - f(u^n_i)}{u^n_{i+1} - u^n_i} & \text{ if $|u^n_i - u^n_{i+1}|> \varepsilon$} ,\\
    f'(u^n_i) & \text{otherwise,}
    \end{cases}
\end{equation}
where $a_{i+1/2}$ is an approximation of the wave speed such as Roe's intermediate speed and 
$\varepsilon$ is a small number (but not too small to avoid numerical cancellation). 
An alternative procedure introduced in \cite{Toro2009} avoids the computation of an intermediate speed by defining $\varphi^1_{i+1/2}=\min( \varphi^1(r_{i+1/2}^R), \varphi^1(r_{i+1/2}^L))$.
This strategy is easily extended to systems by computing the flux limiter component by component.

\paragraph{Smoothness indicators}
The smoothness indicators used in this work  have been introduced in \cite{CPZMR2020}. Nonetheless we briefly recall their construction for the sake of completeness. 
Given a set of the point values $f_j$ of a function $f$ in  the nodes of the stencil $\mathcal{S}^p_{i+1/2}$, $p \geq 2$ we define $\psi_{i+1/2}^p$ as follows. 
First define the lateral left ($^L$) and right ($^R$) weights $w_{i + 1/2}^{p,L/R}$ as
\begin{equation}
        w_{i+1/2}^{p,L}:=\sum_{j=-p+1}^{-1}(f_{i+1+j}-f_{i+j})^2+\varepsilon,\quad
        w_{i + 1/2}^{p,R}:=\sum_{j=1}^{p-1}(f_{i+1+j}-f_{i+j})^2+\varepsilon,
\end{equation}	
where $\varepsilon=10^{-8}$ is a small quantity only used to prevent the weights to vanish. 
Next, using the half harmonic mean,
$w_{i + 1/2}^p = \frac{w_{i + 1/2}^{p,L} w_{i +1/2}^{p,R}}{w_{i + 1/2}^{p,L}+w_{i+1/2}^{p,R}}$,
one defines the high order smoothness indicator over stencil $\mathcal{S}_p$ by 
\begin{equation}\label{indicadores_locales}
        \psi_{i+1/2}^p:=\left( \frac{w_{i + 1/2}^p}{w_{i + 1/2}^p+\tau_{i + 1/2}^p}\right),   
\quad
\text{with}
\quad
\tau_{i + 1/2}^p = \left( (2p - 1)! \sum^{p}_{j=-p+1}\,\gamma^{2p-1, 1/2}_{p,j} \, f^n_{i+j} \right)^2. 
\end{equation}
These indicators are such that
\begin{equation}\label{condiconesS}
    \psi^p_{i+1/2} \approx \left\{
    \begin{array}{cl}
    1 & \mbox { if $\{f_j\}$ are 'smooth' in $\mathcal{S}^p_i$},\\
    0 & \mbox {otherwise}, 
    \end{array}\right.
\end{equation}
see  \cite{CPZMR2020} for a precise statement of this property and its proof. 
\begin{rem}
Observe that, if data in the stencil $\mathcal{S}^p_i$ are smooth, then
$w_{i + 1/2}^{p,L} = O(\Delta x^2)$, 
$w_{i + 1/2}^{p,R} = O(\Delta x^2)$, and 
$\tau_{i + 1/2}^p = O(\Delta x^{4p})$.
Since we adopt the harmonic mean
\[
    \frac{1}{w_{i + 1/2}^p} = \frac{1}{w_{i + 1/2}^{p,L}} + \frac{1}{w_{i + 1/2}^{p,R}},
\]
then $w_{i + 1/2}^p = O(\Delta x^2)$. 
As such
\[
    \psi_{i+1/2}^p =  \frac{w_{i + 1/2}^p}{w_{i + 1/2}^p+\tau_{i + 1/2}^p} =  \frac{O(\Delta x^2)}{O(\Delta x^2)+ O(\Delta x^{4p})},
\]
so that $\psi_{i+1/2}^p$ is expected to be close to $1$. 
On the other hand, if there is an isolated discontinuity in the stencil then 
$
\tau_{i + 1/2}^p = O(1),
$
therefore one of the lateral weights is $O(1)$ and the other $O(\Delta x^2)$ so that the harmonic mean implies that  $w_{i + 1/2}^p = O(\Delta x^2)$ and thus:
\[
    \psi_{i+1/2}^p =  \frac{w_{i + 1/2}^p}{w_{i + 1/2}^p+\tau_{i + 1/2}^p}  =  \frac{O(\Delta x^2)}{O(\Delta x^2)+ O(1)}, 
\]
so that $\psi_{i+1/2}^p$ is expected to be close to 0. 
\end{rem}

\paragraph{ACAT2P schemes}
Using these ingredients, the final expression of the ACAT2$P$ for $P>1$ scheme is of the form
\begin{equation}
    \label{ACAT2P}
u_i^{n+1} = u_i^n + \frac{\Delta t}{\Delta x} \left({F}^{\mathcal{P}_i}_{i-\ha} - 
      {F}^{\mathcal{P}_i}_{i+\ha} \right),
\end{equation}  
where
\begin{equation}\label{numfluxacat}
    {F}^{\mathcal{P}_i}_{i\pm1/2} = \begin{cases}
    {F}^{*}_{i\pm1/2} & \text{if $\mathcal{P}_i = \emptyset$;}\\
    {F}^{p_{\max}}_{i\pm1/2} & \text{otherwise.}
    \end{cases}
\end{equation}
Here, $\mathcal{P}_i$ is the set of consecutive indices 
    \begin{equation}
    \mathcal{P}_i= \{  p \in \{2, \dots, P \} \ \text{s.t.}\  \psi^p_{i+ 1/2} \approx 1 \},
    \quad \text{and}    \quad
    p_{\max} = \max(\mathcal{P}_i).
    \label{eq:Ai}
\end{equation}
Moreover $F^*_{i+1/2}$ is the ACAT2 second-order numerical flux given by \eqref{ACAT2_flux}, $F^{p_{\max}}_{i+1/2}$ is the ACAT2$p_{\max}$ numerical fluxes defined in \eqref{FP+}. 
If $\mathcal{P}_i \not= \emptyset$, the order of the flux $F^{p_{\max}}_{i+1/2}$ can range from $2P$ to $4$.
Throughout the paper we clip $\psi$ to $1$ if $\psi\geq 0.95$.
\begin{rem}
Notice that in \eqref{eq:Ai} index $p$ starts from $2$ since it is not possible to determine the smoothness of the data in the two-point stencil $\mathcal{S}^1_{i+1/2}$.
\end{rem}

\subsection{Extension to 2D} \label{ssec:CAT_2D}
In this section we focus on the extension of CAT methods to non-linear two-dimensional systems of hyperbolic conservation laws
\begin{equation}
\label{2Dequ}
{u}_t + {f}({u})_x + {g}({u})_y=0.
\end{equation}
The following multi-index notation will be used:
$$
\bi = (i_1, i_2) \in \mathbb{Z} \times \mathbb{Z},
$$
and
$$
\bzero = (0,0),\quad \bone = (1,1), \quad \bhalf= (\ha, \ha), \quad \bef = (1,0), \quad \bes = (0,1).
$$
We consider Cartesian meshes with nodes
$$\mathbf{x}_{\bi} = (i_1 \Delta x, i_2 \Delta y).$$
Using this notation, the general form of the CAT$2P$ method will be as follows:
\begin{equation}
\label{2Dscheme}
u_\bi^{n+1}=u_\bi^n + \frac{\Delta t}{\Delta x}\left[ {F}_{\bi -  \frac{1}{2}\bef}^P-{F}_{\bi +  \frac{1}{2}\bef}^P\right] +
\frac{\Delta t}{\Delta y}\left[ {G}_{\bi -  \frac{1}{2}\bes}^P-{G}_{\bi +  \frac{1}{2}\bes}^P\right],  
\end{equation}
where the numerical fluxes ${F}_{\bi +  \frac{1}{2}\bef}^P$,  ${G}_{\bi +  \frac{1}{2}\bes}^P$ will be computed using the values of the numerical solution $U_{\bi}^n$ in the $P^2$-point stencil centered at $\mathbf{x}_{\bi + \bhalf} = 
((i_1 + \ha)\Delta x, (i_2 + \ha)\Delta y)$ 
$$
S_{\bi+\bhalf}^P = \{ \mathbf{x}_{\bi + \bj}, \quad \bj \in \mathcal{I}_P  \},
$$
where 
$$
\mathcal{I}_P =\{ \bj = (j_1, j_2) \in \Z \times \Z, \quad -P +1 \leq j_k \leq P, \quad k = 1,2 \}.
$$

\subsubsection{2D CAT2}
In order to show the extension of CAT2$P$ procedure let us start with the expression of the CAT2. The numerical fluxes are constructed as follows:
\begin{align}
\label{2Dfluxes}
F_{\bi +  \frac{1}{2}\bef}^1 = & \frac{1}{4}\left( {f}^{1,n+1}_{\bi,\bzero} + {f}^{1,n+1}_{\bi, \bef} + f^{n}_{\bi} +f^{n}_{\bi + \bef}\right), \\
G_{\bi +  \frac{1}{2}\bes}^1 = & \frac{1}{4}\left( {g}^{1,n+1}_{\bi,\bzero} + {g}^{1,n+1}_{\bi, \bes} + g^{n}_{\bi} +g^{n}_{\bi + \bes}\right), 
\end{align}
where 
\begin{align*}
f^{1,n+1}_{\bi, \bj} &= f\left(u_{\bi + \bj}^n +\Delta t \, u^{(1)}_{\bi, \bj}\right), \\
g^{1,n+1}_{\bi, \bj} &= g\left(u_{\bi + \bj}^n +\Delta t \,  u^{(1)}_{\bi, \bj}\right),
\end{align*}
for $\bj= \bzero, \bef$ in the $x$ direction, and $\bj = \bzero, \bes$ in the $y$ direction.
\begin{rem}
Despite what happens for the 1D reconstruction, the first time derivative of $u,$ $u^{(1)}_{\bi, \bj},$ does not coincide in the 2D-grid points. Indeed, observe  that $u^{(1)}_{\bi, \bzero} \neq u^{(1)}_{\bi, \bef}$ and $u^{(1)}_{\bi, \bzero} \neq u^{(1)}_{\bi, \bes}.$ 

Note that, in the 1D case, $u^{(1)}_{i,0 }= u^{(1)}_{i, 1} $.
\end{rem}
Hence, the first time derivatives $u^{(1)}_{\bi,\bj}$ are so defined: 
\begin{align*}
u^{(1)}_{\bi, \bzero} &= -\frac{1}{\Delta x}\left( f_{\bi + \bef}^n - f_{\bi}^n\right) 
-\frac{1}{\Delta y}\left( g_{\bi + \bes}^n - g_{\bi}^n\right),  \\
u^{(1)}_{\bi, \bef} &= -\frac{1}{\Delta x}\left( f_{\bi + \bef}^n - f_{\bi}^n\right)
-\frac{1}{\Delta y}\left( g_{\bi + \bone}^n - g_{\bi + \bef}^n\right), \\ 
u^{(1)}_{\bi, \bes} &= -\frac{1}{\Delta x}\left( f_{\bi + \bone}^n - f_{\bi + \bes}^n\right) 
-\frac{1}{\Delta y}\left( g_{\bi + \bes}^n - g_{\bi}^n\right), 
\end{align*} 
where 
$$ f_{\bi+\bj}^n = f(u_{\bi+\bj}^n), \quad g_{\bi+\bj}^n =g(u_{\bi+\bj}^n), \quad \forall \bj.$$
Finally, the 2D CAT2 method is so defined:
\begin{equation}
\label{2D_CAT2_cons}
u_\bi^{n+1}=u_\bi^n + \frac{\Delta t}{\Delta x}\left[ {F}_{\bi -  \frac{1}{2}\bef}^1-{F}_{\bi +  \frac{1}{2}\bef}^1\right] +
\frac{\Delta t}{\Delta y}\left[ {G}_{\bi -  \frac{1}{2}\bes}^1-{G}_{\bi +  \frac{1}{2}\bes}^1\right],  
\end{equation}

\subsubsection{2D CAT2P}
The high order CAT2$P$ iterative procedure are computed as follows:

\begin{enumerate}

\item  {Define}
$$
f^{(0)}_{\bi,\bj}=f^n_{\bi + \bj}, \quad  g^{(0)}_{\bi,\bj}=g^n_{\bi + \bj}, \quad \bj \in \mathcal{I}_P.
$$

\item  {For $k = 2, \dots, 2P$:}

\begin{enumerate}
\item Compute
\begin{equation*}
 u^{(k-1)}_{\bi,\bj} = - A^{1,j_1}_{P}(f^{(k-2)}_{\bi,(\bcdot, j_2)}, \Delta x) 
 - A^{1,j_2}_{P}(g^{(k-2)}_{\bi,(j_1, \bcdot)}, \Delta y), \quad \bj \in \mathcal{I}_P.
\end{equation*}
\item Compute 
$$
f^{k-1,n+r}_{\bi,\bj} = f \left(  u^n_{\bi + \bj} + 
\sum_{l=1}^{k-1} \frac{(r \Delta t)^l}{l!} u^{(l)}_{\bi,\bj} \right), \quad  \bj\in \mathcal{I}_P,\;\, r = -P+1, \dots, P.
$$
\item Compute
$$
f^{(k-1)}_{\bi,\bj} =   A^{k-1,0}_{P}( f^{k-1, \bcdot}_{\bi,\bj}, \Delta t),\quad  \bj \in \mathcal{I}_P.
$$
\end{enumerate}

\item Compute
\begin{eqnarray}\label{cat22Dx}
F^P_{\bi+\frac{1}{2}\bef}  &= &\sum_{k=1}^{2P} 
\frac{\Delta t^{k-1}}{k!}A^{0, 1/2}_{P}(\tilde{f}_{\bi,(\bcdot, 0)}^{(k-1)}, \Delta x), 
\\\label{cat22Dy}
G^p_{\bi+\frac{1}{2}\bes}  &= &\sum_{k=1}^{2P} 
\frac{\Delta t^{k-1}}{k!}A^{0, 1/2}_{P}(\tilde{g}_{\bi, (0, \bcdot)}^{(k-1)}, \Delta y).
\end{eqnarray}

\end{enumerate}
The notation used for the approximation of the spatial partial derivatives is the following:
\begin{eqnarray*}
 A^{k,q}_{P}(f_{\bi, (\bcdot, j_2)}, \Delta x) & =&  \frac{1}{\Delta x^k} \sum_{l = -P + 1}^P \gamma^{k,q}_{P,l} f_{\bi, (l,j_2)}\\
 A^{k,q}_{P}(g_{\bi, (j_1,\bcdot{} )}, \Delta y) & =&  \frac{1}{\Delta y^k} \sum_{l = -P + 1}^P \gamma^{k,q}_{P,l} g_{\bi, (j_1,l)}
 \end{eqnarray*}
 
 \begin{rem}
 In the last step of the algorithm above the set $\mathcal{I}_P$ can be replaced by its $(2P -1)$-point subset
 $$
 \mathcal{I}^0_P = \{ \bj = (j_1, j_2)\; \text{ \rm{such that} } \;j_1 = 0 \text{ or } j_2 = 0 \}
 $$
 since only the corresponding values of $\tilde{f}^{(k-1)}_{\bi,\bj}$ are  used to compute the numerical fluxes \eqref{cat22Dx} and
 \eqref{cat22Dy}.
 \end{rem}

\subsection{Discussion} \label{ssec:beyondCAT_discusssion}
The ACAT numerical method presented in previous sections is exhaustively described in \cite{CPZMR2020,Macca-Pares} where numerical results are provided.
The extension to systems of conservation laws is done component-by-component, and the extension to multi-dimensions relies on a direction-by-direction splitting.
However we can point out several defects in the design of such an approach.
First of all, the impact of the \apriori smoothness indicators is tremendous. Indeed, they must effectively detect the presence of a discontinuity, but also the occurrence of new ones. While the former is achievable with \apriori limiters, the latter is usually more difficult.
Moreover the smoothness indicators can detect discontinuity only if the mesh is fine enough, otherwise it is hard to distinguish between a discontinuity and a smooth region with a large gradient.
Secondly, the smoothness indicator cost drastically increases with $P$. Moreover, it becomes increasingly more difficult to determine the smoothness of a numerical solution when higher and higher order of accuracy is required.
Thirdly, just like any \apriori limiter, the current one has a natural tendency to over-estimate and over-react to possible spurious troubles. As such, the nominal order of accuracy on smooth solution is not always achieved, unless the grid is really fine. \\
These facts have restricted the effective use of ACAT schemes for orders not greater than $2P=6$. 
Also the cost of the current version of ACAT2P constitutes a genuine limitation in 2D. \\
In this work we present an alternative way to couple CAT scheme with an \aposteriori limiting technique, called MOOD, and remove
some of the previously described defects.

%
%
\section{CAT-MOOD} \label{sec:CATMOOD}
The main objective of this work is to combine the \aposteriori shock capturing technique (MOOD) \cite{CDL1} to this family of one-step spatial and temporal reconstructions of high order CAT schemes \cite{CATMOOD_1D}.  

The MOOD algorithm evaluates \aposteriori the solution of the high-order numerical method using a class of criteria that detect a variety of oscillations, even very small ones. It is therefore possible to ensure that the numerical solution preserves some properties of the exact solutions of the PDE system, such as, for example, positivity, monotonicity and increase of the physical entropy, even in complex cases. 
In addition, this procedure allows to reduce and eliminate the numerical oscillations introduced by the high order methods in the presence of shocks or large gradients as is the case for the CAT methods. 

The basic idea of the MOOD procedure is to apply a high-order method over the entire domain for a time step, then check locally, for each cell $i$, the behavior of the solution via admissibility criteria. If the solution computed in cell $i$ at time $t^{n+1}$ is in accordance with the criteria considered, it is kept, otherwise, it is recomputed with a new numerical method of order lower than the previous one. This operation is repeated until acceptability, or, when the (last) first order scheme is used. This last case occurs when the admissibility criteria fail with any previous reconstruction.

Therefore, the object of this work is to design a cascade of CAT methods in which the order is locally adapted according to \aposteriori admissibility criteria thus creating a new family of adaptive CAT methods called CAT-MOOD schemes.

\subsection{MOOD admissibility criteria} \label{ssec:MOODdetection}
\begin{figure}[!ht]
    \centering    
    \includegraphics[width=0.7\textwidth]{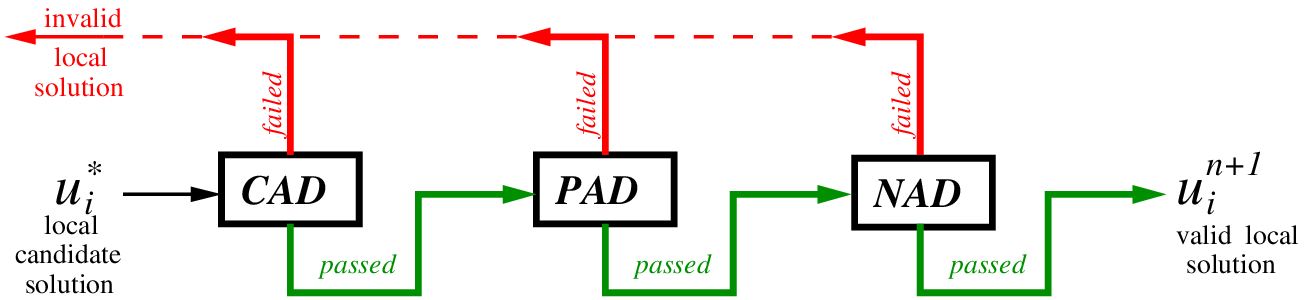}
    \includegraphics[width=0.28\textwidth]{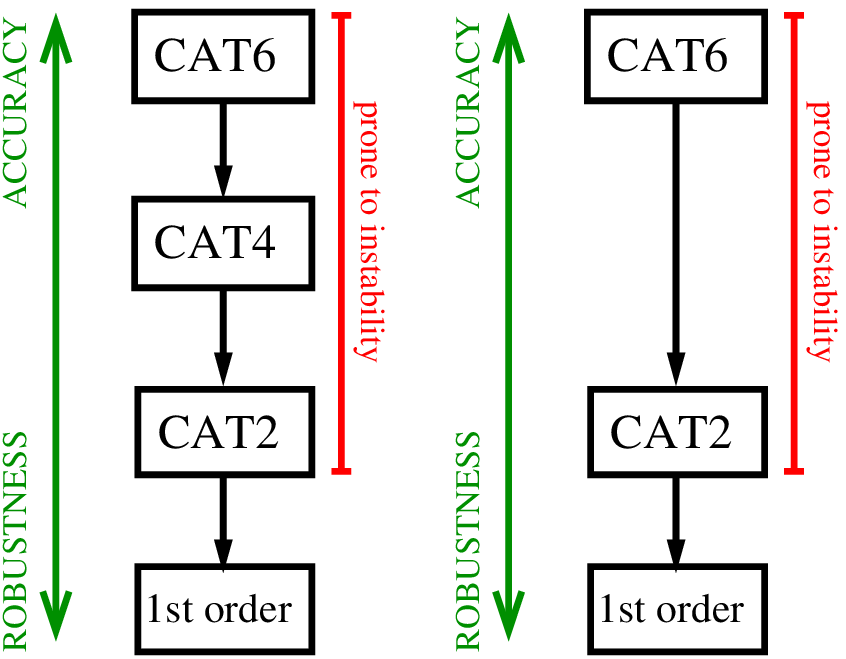}
    \caption{Left: Detection criteria of the MOOD procedure for a candidate solution $u_i^*$. \textit{Computer Admissible Detector (CAD)}, \textit{Physical Admissible Detector (PAD)}  and \textit{Numerical Admissible Detector (PAD)} ---
    Right: Cascades of CAT schemes used in the MOOD procedure. Starting from the most accurate one, CAT6, downgrading to lower order schemes, and, at last to a $1$st order accurate scheme employed to ensure robustness.}
   \label{fig:MOOD_chain_cascade}
\end{figure}

In this work we select 3 different admissibility criteria 
\cite{CDL2,CDL3} which are invoked onto the candidate numerical solution $\left\{ u_i^{n+1} \right\}_{1\leq i \leq I}$:
\begin{enumerate}
    \item 
    \textit{Computer Admissible Detector (CAD)}: This criterion is responsible to detect 
    undefined or un-representable quantities, usually
    not-a-number \texttt{NaN} or infinity quantity due to some division by zero for instance.
    \item 
    \textit{Physical Admissible Detector (PAD)}: The second detector is responsible for ensuring the physical validity of the candidate solution. The detector reacts to every negative pressure $p$ or density $\rho$ in the computational domain, 
    in compliance with (\ref{eq:physical_states}), since otherwise the solution will create non-physical sound speeds, imaginary time steps and so on. 
    The physicality here is assessed from the point of view of a fluid flow, which limits the generality of the criteria as far as predicted pressures are concerned. 
    Said differently, this physical admissibility criteria must be adapted to the model of PDEs which is solved.
    \item 
    \textit{Numerical Admissible Detector (NAD)}: This criterion corresponds to a relaxed variant of the Discrete Maximum Principle (see \cite{Ciarlet,CDL1})
    $$\min_{c\in \mathcal{C}_i^P}(w^n_c) - \delta_i \le w^*_i \le \max_{c\in \mathcal{C}_i^P}(w^n_c) + \delta_i, $$
    where $\mathcal{C}_i^P = \{-P,\ldots,P\}$ is the local centered stencil of order $2P$ and $\delta_i$ is a relaxation term to avoid problems in flat region. $w^*_i$ is the numerical solution obtained with the scheme of order $2P$ and $\delta_i$ is set as:
    \begin{equation}
        \label{eq:delta}
        \delta_i = \max\Big(\varepsilon_1,\varepsilon_2\left[\max_{c\in \mathcal{C}_i^P}w^n_c - \min_{c\in \mathcal{C}_i^P}w^n_c\right]\Big). 
    \end{equation}
    Here $\varepsilon_1$ and $\varepsilon_2$ are small dimensional constant. 
    This criterion is responsible for guaranteeing the essentially non oscillatory (ENO) character of the solution; that is, no large and spurious minima or maxima are introduced locally in the solution\footnote{In the numerical tests of Sec.~\ref{sec:numerics} we compute the relaxed discrete maximum principle only for density $\rho$ and pressure $p$. No limiting, on the velocity components, $u$, $v$, has been considered.}. 
\end{enumerate}    
If a \texttt{NaN} number is detected by CAD in the candidate solution $w_i^*$, then this cell is sent back for recomputation right away. 
Next, if the candidate solution has some positivity issues, then the PAD test is not passed, and, the cell is also invalid.
At last the cell enters the NAD criteria to test for possible numerical oscillation.
If spurious oscillations have contaminated the candidate solution, $w_i^*$, then it will fail NAD.
The candidate solution is then recomputed if at least one of the previous criteria ordered into a chain, see figure~\ref{fig:MOOD_chain_cascade}-left,
has failed.
As a consequence, the MOOD loop drives the code to locally downgrade the order of accuracy by using an auxiliary scheme of lower accuracy.

\subsection{CAT scheme with MOOD limiting} \label{ssec:CATMOOD} 
In this work we target to reach a maximal $6$th order of accuracy on the part of the domain where a smooth solution is present. 
This can be achieved with CAT6 scheme that we would like to employ as often as possible.
On the contrary, for cells presenting a discontinuous solution we plan to rely on a $1$st order low accurate but robust scheme, for instance using Rusanov or HLL fluxes, that we would employ only when and where necessary. 
In between these two extremes, the CAT4 schemes of $4$th order of accuracy is inserted and tried when CAT6 fails. If CAT4 scheme also fails, then CAT2 scheme ($2$nd order of accuracy) is further tried.
As such we build several cascades of schemes of decreasing orders but possibly increasing robustness, see figure~\ref{fig:MOOD_chain_cascade}-right.
Notice that, apart from CAT6 and the 1st order parachute scheme, the user can decide if the intermediate schemes (CAT4 and CAT2 here) should be included or not. 
In the numerical results in section~\ref{sec:numerics} only the CAT2 scheme is employed to spare computing resources.
This scheme is referred to as CATMOOD6.

\subsection{CATMOOD algorithm} \label{ssec:MOODalgo}
Practically, for the time-step $[t^n,t^{n+1}]$ and for each cell $i$, we define a 'mask', $M_i^{n}\in\{-1,0,1\}$, such that
\bea \label{eq:mask}
M_i^n = \begin{cases}
           1 \quad & \text{if} \; u_i^* \; \text{fails at least one criterion},  \\ 
           -1\quad & \text{if} \; \exists  j \in \mathcal{N}_i \; \text{s.t} \; M_j^n=1,\\ 
           0 \quad & \text{otherwise}.
         \end{cases}
\eea
where $\mathcal{N}_i$ is the set of direct neighbor cells of cell $\omega_i$. $M_i^n=-1$ means that cell $i$ is the neighbor of an invalid cell. \\
The algorithm designed for CATMOOD6 scheme is:
\begin{enumerate}
    \item[Init] Let  $\{u_i^n\}_{1\leq c\leq N_c}$ be the numerical solution at $t = t^n$ over the whole domain $\Omega.$
    \item[CAT6] Let $\{u_i^*\}_{1\leq c\leq N_c}$ be the numerical solution at $t = t^{n+1}$ of order 6 obtained from CAT6 scheme. Set $M_i^n=0$.\\
    For all cell $i$, check if $u_i^*$ satisfies all detection criteria. 
    In this case, then $u_i^{n+1} = u_i^*$ and $M_i^n=0$.
    Otherwise $M_i^n=1$ and, if $M_j^n=0$ then set $M_j^n=-1$ for all $j \in \mathcal{N}_i$.
    \item[CAT4] Only for the troubled cell $i$, i.e $M_i^n\not= 0$, recompute $\{u_i^*\}$ the numerical solution at time $t = t^{n+1}$ of order 4 obtained with CAT4 scheme. \\
   Check if $u_i^*$ satisfies all detection criteria.   
    In this case, then $u_i^{n+1} = u_i^*$ and $M_i^n=0$.
    Otherwise $M_i^n=1$ and, if $M_j^n=0$ then set $M_j^n=-1$ for all $j \in \mathcal{N}_i$.
    \item[CAT2] Only for the troubled cell $i$, i.e $M_i^n\not= 0$, recompute $\{u_i^*\}$ the numerical solution at time $t = t^{n+1}$ of order 2 obtained with CAT2 scheme. \\
   Check if $u_i^*$ satisfies all detection criteria.   
    In this case, then $u_i^{n+1} = u_i^*$ and $M_i^n=0$.
    Otherwise $M_i^n=1$ and, if $M_j^n=0$ then set $M_j^n=-1$ for all $j \in \mathcal{N}_i$.
    \item[1st-ord] Only for the remaining troubled cell $i$, i.e $M_i^n\not= 0$, recompute $\{u_i^*\}$ the numerical solution at time $t = t^{n+1}$ of order 1 obtained with a first order scheme, and set $u_i^{n+1} = u_i^*$.
\end{enumerate}
For efficiency purposes we may remove the CAT4 step in the cascade.

\subsection{Complexity, cost, convergence, implementation}
The unavoidable extra-cost when using a MOOD procedure is the detection of troubled cells with the admissible criteria from section~\ref{ssec:MOODdetection}.
However the detection criteria are particularly inexpensive to compute compared to smoothness indicators or sensors used for instance in ACAT scheme. Usually these admissibility criteria have a negligible cost. \\
The MOOD limiting procedure always converges because the number of cells and schemes are finite. In the worst case scenario the entire solution is computed successively with all CAT2P schemes up to the 1st order accurate solution.
In the best case scenario the solution from the unlimited CAT6 scheme is accepted without any correction.
Any situation in-between is possible, and generally only few cells need to be recomputed. \\
The detection criteria are fundamental, they must be designed to ensure that, if the mesh is fine enough, a smooth solution computed by CAT6 scheme does not produce any troubled cell. 
They must also ensure that, in the vicinity of strong discontinuity, the robust 1st order scheme is regularly employed to avoid spurious oscillations. \\
Concerning the complexity of CATMOOD, it becomes impossible to estimate \textit{a priori} the cost because the limiting adapts to the underlying flow and the computed numerical solution. In the best case scenario mentioned above, the cost of CATMOOD6 is the one of CAT6, in the worst, the cost of CATMOOD6 is the sum of the cost of all schemes in the cascade.
Generally one observes that the amount of troubled cells is of the order of $0-20\%$ of the total number of cells, which makes the MOOD procedure genuinely competitive compared to existing \textit{a priori} limiters, see \cite{CDL2,ADER_MOOD_14} and the numerical section in this paper. 

\subsection{CATMOOD vs LATMOOD}
In the context of the LAT methods proposed by D. Zorío et al. \cite{ZBM2017}, it is noteworthy that the number of operations per time step is significantly reduced compared to CAT$2P$ methods. This reduction is achieved by performing the Taylor expansions only once per point, while CAT$2P$ requires $2P$ expansions per point. However, it is important to acknowledge that CAT$2P$ exhibits superior stability properties compared to LAT.

Consequently, an intriguing question arises: can LATMOOD yield a more efficient approach than CATMOOD? To explore this question, it is important to consider two crucial points: 
\begin{enumerate}
    \item The local computation of numerical fluxes makes CAT more suitable for implementing MOOD. Indeed, unlike CAT schemes, LAT methods compute temporal derivatives in a global manner, which contradicts the local approach of MOOD. Specifically, in CAT methods, once the solution $u^*$ has been computed, if any of the detectors fail, it is sufficient to recompute all local time derivatives of the fluxes with a lower order of accuracy, without causing any issues in the MOOD approach. Contrarily in the case of LAT, since all the approximations are computed in a non local sense (see Remark~\ref{rem:local}) and since LAT does not use a compact stencil, $4P$ computations have to be updated for each bad cell and approximation $u_{i+p}^{(1)},$ $u_{i+p}^{1,n+1},$ etc., making the scheme computationally inefficient.
    \item From our preliminaries tests we notice that the proportion of low order cells substantially increases when LAT is employed instead of CAT. 
\end{enumerate}


%
\section{Numerical test cases}  \label{sec:numerics}
In this paper only numerical tests related to the two-dimensional Euler equations are considered. Obviously, CATMOOD schemes could be applied to any systems of conservation laws without any restriction.
Our methodology of testing relies on several classical and demanding test cases: 
\begin{enumerate}
    \item \textit{Isentropic vortex in motion}. This test measures the ability of the MOOD procedure combined with CAT schemes to achieve the optimal high order for a smooth solution. We also compared CATMOOD6 with unlimited CAT schemes, limited ACAT ones and some first order scheme using Rusanov, HLL and HLLC fluxes.
    \item \textit{Sedov Blast wave.} This problem has an exact solution presenting a single cylindrical shock wave followed by an exponential decay \cite{Sedov}. This test is used to check the behaviour of the CATMOOD6 scheme against shocks and to compare CATMOOD6 versus ACAT6.
    \item \textit{2D Riemann problems}. We simulate four versions of the four-state 2D Riemann problems \cite{RP_Schultz93}. These problems present large smooth regions, unstable shear layers, unaligned contact discontinuities and shock waves, along with complex interaction patches for which no exact solution has yet been derived.
    \item \textit{Astrophysical jet}. This intense and demanding test case challenges the positivity and robustness of the CATMOOD6 scheme as it presents an extremely violent jet at Mach 2000 generating a bow shock and unstable shear layers. 
\end{enumerate}
The time step is chosen according to 
$$
\Delta t^n = \mathrm{CFL}\, \min\left( \frac{\Delta x}{\lambda_{ {\rm max}_x}^n}, \frac{\Delta y}{\lambda_{ {\rm max}_y}^n}\right),   $$
where $\lambda_{{\rm max}_x}^n$ and $\lambda_{{\rm max}_y}^n$ represent the maximum, over cells, of the spectral radius of the Jacobian matrices $\partial \mathbf{F}/\partial u$ and $\partial \mathbf{G}/\partial u$. 
Empirically we found that the CFL number should be less or equal to $0.5$.
In our calculations we choose CFL$=0.4$.
In all the numerical examples below (except the one in Section \ref{ssec:Mach}), the units are chosen in such a way that the order of magnitude of all field quantities is one, therefore we always use the same value for the constants appearing in \eqref{eq:delta}, i.e.\ $\varepsilon_1=10^{-4}$ and $\varepsilon_2=10^{-3}$, while the value $\varepsilon = 10^{-8}$ has been used in \eqref{wavefuction}.

\subsection{Isentropic vortex in motion} \label{ssec:vortex}
The isentropic vortex problem \cite{Shu1} challenges the accuracy of numerical methods since an exact, smooth and analytic solution exists. The computational domain is set to $\Omega = [-10,10]\times [-10,10]$.
The ambient flow is characterized by $\rho_\infty = 1.0$, $u_\infty = 1.0$, $v_\infty = 1.0$ and $p_\infty = 1.0$,
with a normalized ambient temperature $T^*_\infty =1.0$. 
At the initial time, $t=0$, onto this ambient flow is superimposed a vortex centered at $(0,0)$ with the following state:
$u = u_\infty + \delta u$, $ v = v_\infty + \delta v$, 
$ T^* = T^*_\infty + \delta T^*$, where the increments are given by
\begin{equation}
    \begin{array}{c}
      \nonumber
      \delta u   = -y' {\dfrac {\beta} {2 \pi}} \exp \left( {\dfrac {1-r^2} {2}} \right),  \quad 
      \delta v   = x' {\dfrac {\beta} {2 \pi}} \exp \left( {\dfrac {1-r^2} {2}} \right), \quad
      \delta T = - { \dfrac {(\gamma - 1 ) \beta^2} {8 \gamma \pi^2}} \exp \left( {1-r^2} \right),
    \end{array}
\end{equation}
with $r = \sqrt{x^2 + y^2}$. The so-called strength of the vortex is set to $\beta=5.0$ and the initial density is
given by  $\rho =   \rho_\infty \left( {{T}/{T_\infty} } \right)^{\frac{1}{\gamma-1} }$. Periodic boundary conditions are prescribed. 
At final time $t=t_{\text{final}}=20$ the vortex is back to its original position, and, the final exact solution matches the initial one. 
Since the solution is smooth, it should be simulated with optimal high accuracy, in other words, the limiting/stabilization procedure employed in the scheme should not have any effect. 
\begin{table}[!ht]
\numerikNine
    \centering
    \begin{tabular}{|c||cc|cc|cc|cc|}
    \hline \multicolumn{9}{|c|}{\textbf{2D Isentropic Vortex in motion - Rate of convergence}} \\
\hline   
\hline   
& \multicolumn{2}{c|}{\textbf{Rusanov-flux}} & \multicolumn{2}{c|}{\textbf{HLL}} & \multicolumn{2}{c|}{\textbf{HLLC}} & \multicolumn{2}{c|}{\textbf{CATMOOD6}}\\
 $N$ &  $L^1$ error &  order & $L^1$ error &  order & $L^1$ error &  order & $L^1$ error &  order \\ \hline
50 $\times$ 50  & 8.44$\times 10^{-3}$  &  ---    & 8.44$\times 10^{-3}$ & ---   &7.91$\times 10^{-3}$  &  -& 8.48$\times 10^{-3}$ & ---     \\
100 $\times$ 100 & 8.04$\times 10^{-3}$ &  0.07 & 8.04$\times 10^{-3}$ & 0.07&6.86$\times 10^{-3}$ &  0.21& 3.77$\times 10^{-3}$ & 1.17 \\
200 $\times$ 200 & 6.68$\times 10^{-3}$  & 0.27 & 6.67$\times 10^{-3}$ & 0.27&5.31$\times 10^{-3}$ &  0.37& 2.40$\times 10^{-7}$ & 13.94  \\
300 $\times$ 300 & 5.71$\times 10^{-3}$  & 0.36 & 5.71$\times 10^{-3}$ & 0.36&4.53$\times 10^{-3}$ &  0.39& 2.06$\times 10^{-8}$ & 6.05   \\
400 $\times$ 400 & 4.98$\times 10^{-3}$ &  0.47 & 4.98$\times 10^{-3}$ & 0.47&3.86$\times 10^{-3}$ &  0.55& 3.52$\times 10^{-9}$ & 6.14   \\
& \text{Expected} & 1 
& \text{Expected} & 1 
& \text{Expected} & 1 
& \text{Expected} & 6
\\
\hline  \hline  
& \multicolumn{2}{c|}{\textbf{CAT2}} & \multicolumn{2}{c|}{\textbf{CAT4}} & \multicolumn{2}{c|}{\textbf{CAT6}} &\multicolumn{2}{c|}{\textbf{ACAT6}} \\ 
     $N$ &  $L^1$ error &  order & $L^1$ error &  order & $L^1$ error &  order & $L^1$ error &  order  \\ \hline
50 $\times$ 50  & 7.94$\times 10^{-3}$  &  ---    & 2.03$\times 10^{-3}$ & ---   &8.46$\times 10^{-4}$  &  -& 8.95$\times 10^{-3}$ & --- \\
100 $\times$ 100 & 2.55$\times 10^{-3}$ &  1.64 & 1.42$\times 10^{-4}$ & 3.83&1.56$\times 10^{-5}$ &  5.76& 8.28$\times 10^{-3}$ & 0.11\\
200 $\times$ 200 & 6.12$\times 10^{-4}$ &  2.06 & 8.34$\times 10^{-6}$ & 4.09&2.41$\times 10^{-7}$ &  6.02& 8.34$\times 10^{-5}$ & 9.95\\
300 $\times$ 300 & 2.69$\times 10^{-4}$ &  2.02 & 1.64$\times 10^{-6}$ & 4.02&2.09$\times 10^{-8}$ &  6.03& 1.05$\times 10^{-5}$ & 5.14\\
400 $\times$ 400 & 1.52$\times 10^{-4}$ &  1.99 & 5.16$\times 10^{-7}$ & 4.01&3.68$\times 10^{-9}$ &  6.03& 2.48$\times 10^{-6}$ & 4.93\\
& \text{Expected} & 2
& \text{Expected} & 4 
& \text{Expected} & 6 
& \text{Expected} & 6
\\
    \hline 
    \end{tabular}
    \caption{Isentropic vortex in motion  $L^1-$norm errors on density $\rho$ between the numerical solution and the exact solution of the isentropic vortex in motion problem at $t_{\text{final}} = 20$ on uniform Cartesian mesh. }
    \label{tab:vortex_convergence}
\end{table}

\begin{table}[!ht]
    \centering
    \numerikNine
 \resizebox{\columnwidth}{!}{
    \begin{tabular}{|c||c| c| c| c|c|c|c|c|}
    \hline 
    & \multicolumn{3}{c|}{1st order methods} & 2nd ord. & 4th ord. & \multicolumn{3}{c|}{6th order methods} \\
    \cline{2-9}
    & \textbf{Rusanov} & \textbf{HLL} &  \textbf{HLLC}   &  \textbf{CAT2} & 
    \textbf{CAT4} & \textbf{CAT6} &  \textbf{CATMOOD6}   &  \textbf{ACAT6} \\
    \hline \hline
    \textbf{CPU(s)} & 48.22 & 60.58   & 98.41 & 424.72  & 4375.68 & 16112.38   & 20848.26 & 23597.67    \\ \hline
    \textbf{Ratio} & 1 & 1.26  & 2.04 &  8.81 & 90.74 & 334.14  & 432.36 &  489.38\\ 
    \hline
    \end{tabular} 
    }
    \caption{Isentropic vortex in motion with $300\times300$ cells. First row: computational costs expressed in seconds. Second row: ratio with respect to the Rusanov cost.}
    \label{tab:vortex_cpu}
\end{table}
\begin{figure}[!ht]
    \centering    
    \includegraphics[width=0.495\textwidth]{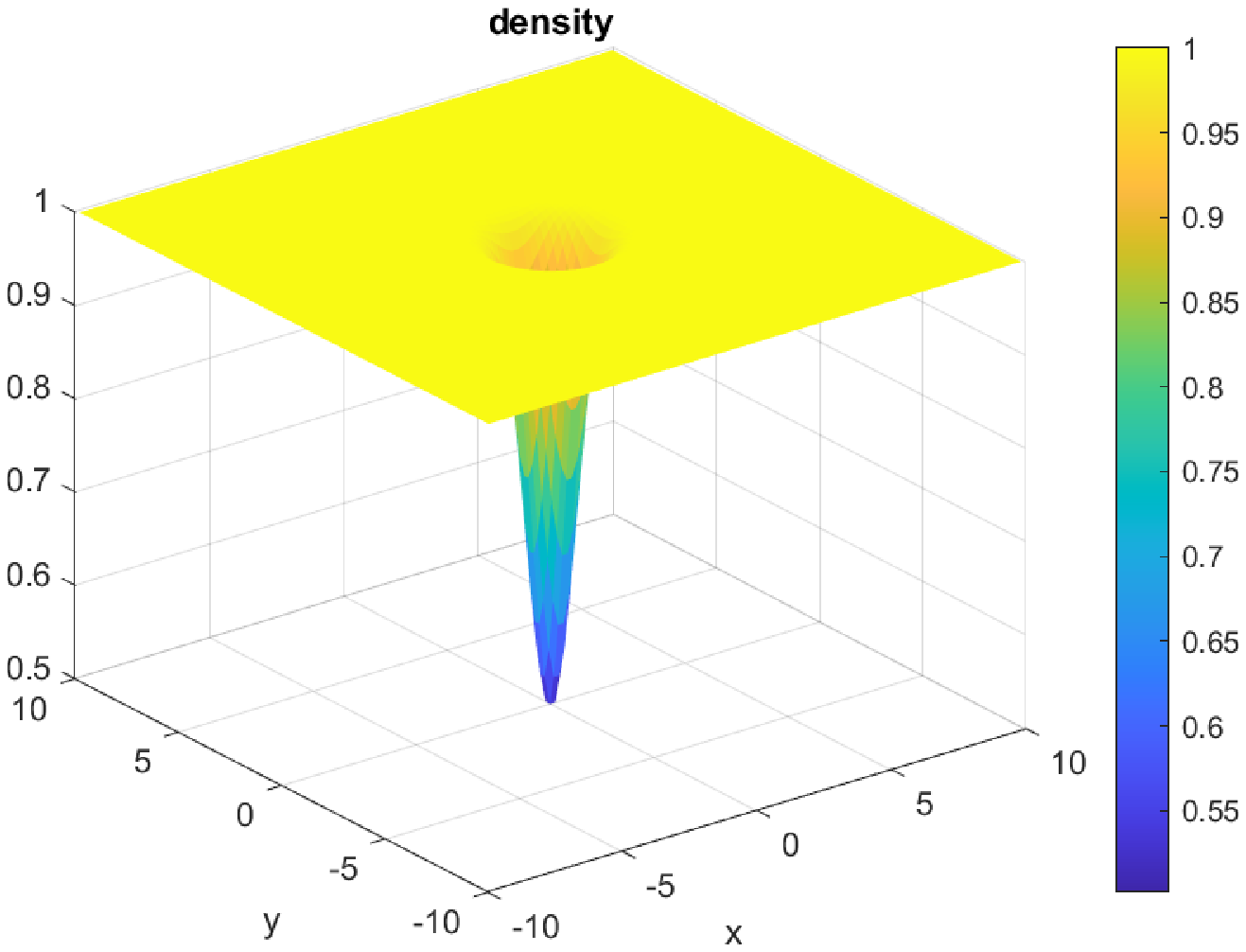}
    \includegraphics[width=0.495\textwidth]{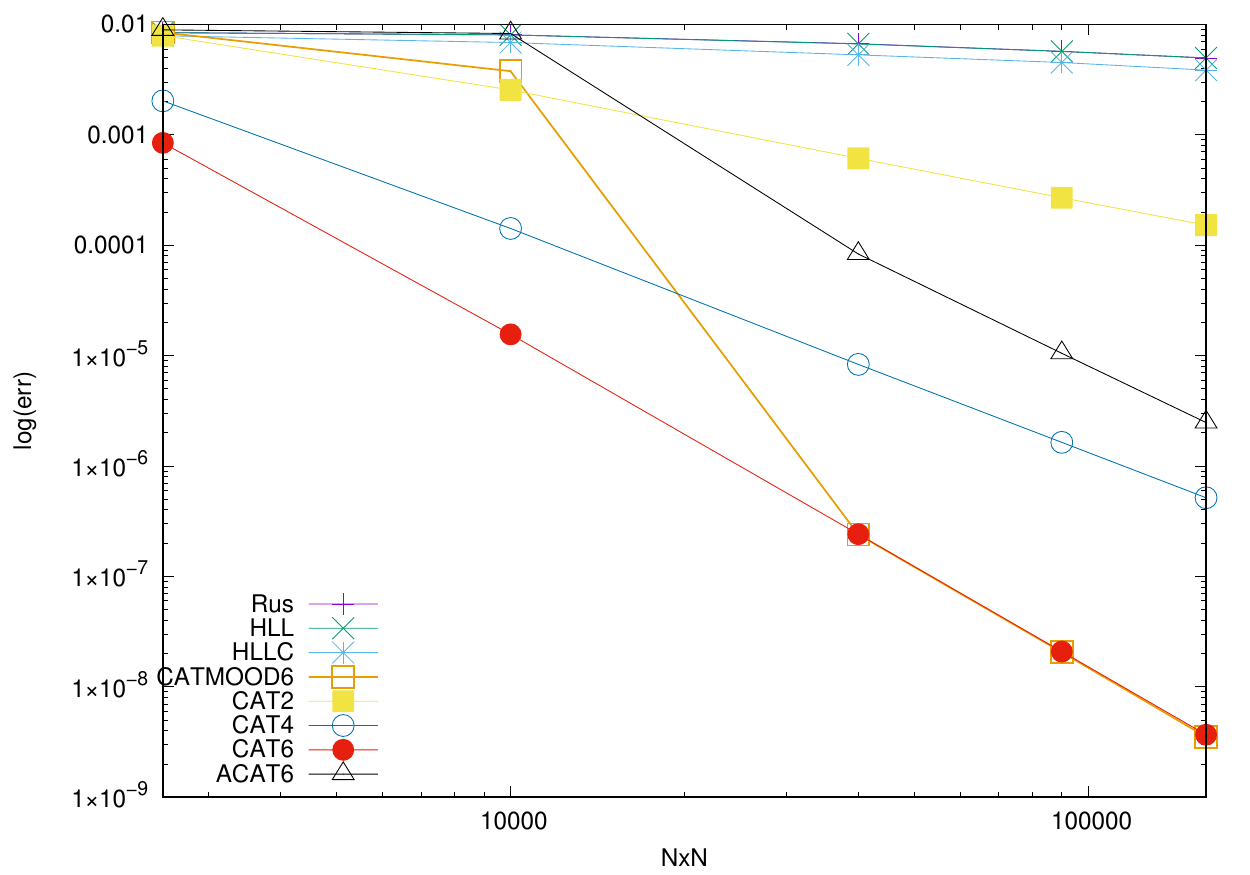}
    \caption{2D isentropic vortex in motion ---
    Left: Numerical solution at final time for density with CATMOOD6 with HLLC flux scheme used for the first order method on $100 \times 100$ uniform mesh and CFL$=0.4$ ---
    Right: Errors in $L_1$ norms given by all schemes.}
   \label{fig:vortex_density}
\end{figure}
We run the isentropic vortex test case with first order Rusanov-flux, HLL, HLLC and  sixth order CATMOOD method on successive refined Cartesian meshes going from $50 \times 50$  up to  $400 \times 400$ cells. 
The density at final time is plotted in figure~\ref{fig:vortex_density}-left.
The results of the convergence analysis for all the schemes are displayed in table~\ref{tab:vortex_convergence} and in figure~\ref{fig:vortex_density}-right.
The expected rate of convergence is reached for CAT schemes, while for the first-order schemes the convergence is below the expected order.
Ideally the limited $6$th-order ACAT6 scheme should produce the same errors as the ones given by CAT6 scheme, but ACAT6 errors are two orders of magnitude greater.
Contrarily, starting at mesh $200\times 200$, CATMOOD6 errors match the optimal ones from CAT6, and, the nominal $6$th order is retrieved. This proves that the \aposteriori MOOD limiting avoids spurious intervention if the mesh is fine enough.
Table~\ref{tab:vortex_cpu} presents the computational costs of various methods for  the isentropic vortex test case on a uniform grid with $300\times300$ cells, CFL$ = 0.4,$ final time $t_{\text{final}} = 20$, and periodic boundary conditions. The first row displays the CPU computational costs in seconds for the eight methods, while the second row shows the ratio between the computational cost of each scheme and the one of the Rusanov method. This table emphasizes the superiority of the \aposteriori approach over the \apriori one. Indeed ACAT6 scheme appears to be not only less accurate than CATMOOD6, but also its execution time is higher.


\subsection{Sedov blast wave}\label{ssec:Sedov}
\begin{figure}[!ht]
     \centering
     \begin{subfigure}[b]{0.49\textwidth}
         \centering
         \includegraphics[width=\textwidth]{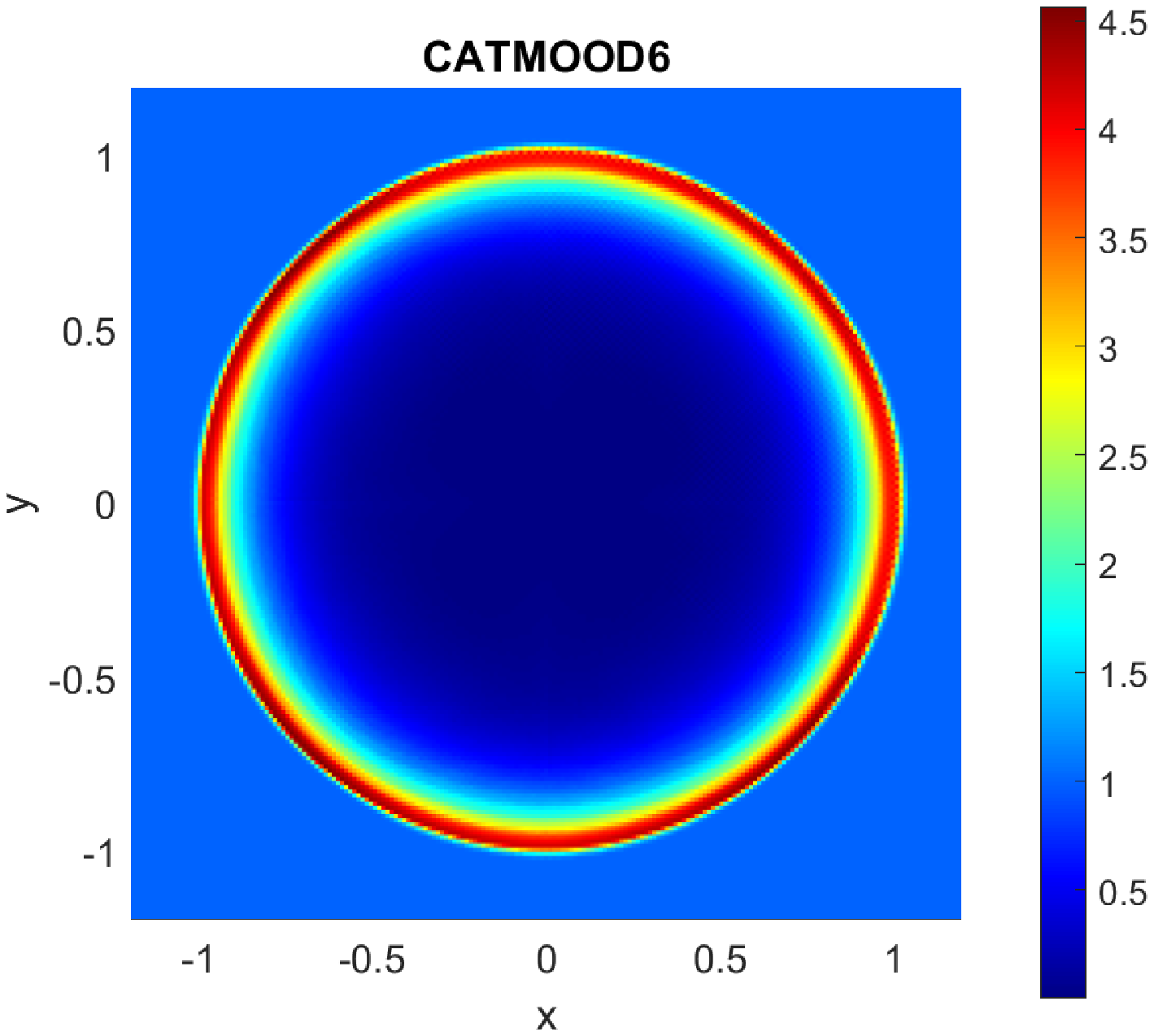}
         \caption{Projection of the numerical density for CATMOOD6 with HLLC as first order scheme}
         \label{Sed:CATMOOD6}
     \end{subfigure}
     \hfill
     \begin{subfigure}[b]{0.49\textwidth}
         \centering
         \includegraphics[width=\textwidth]{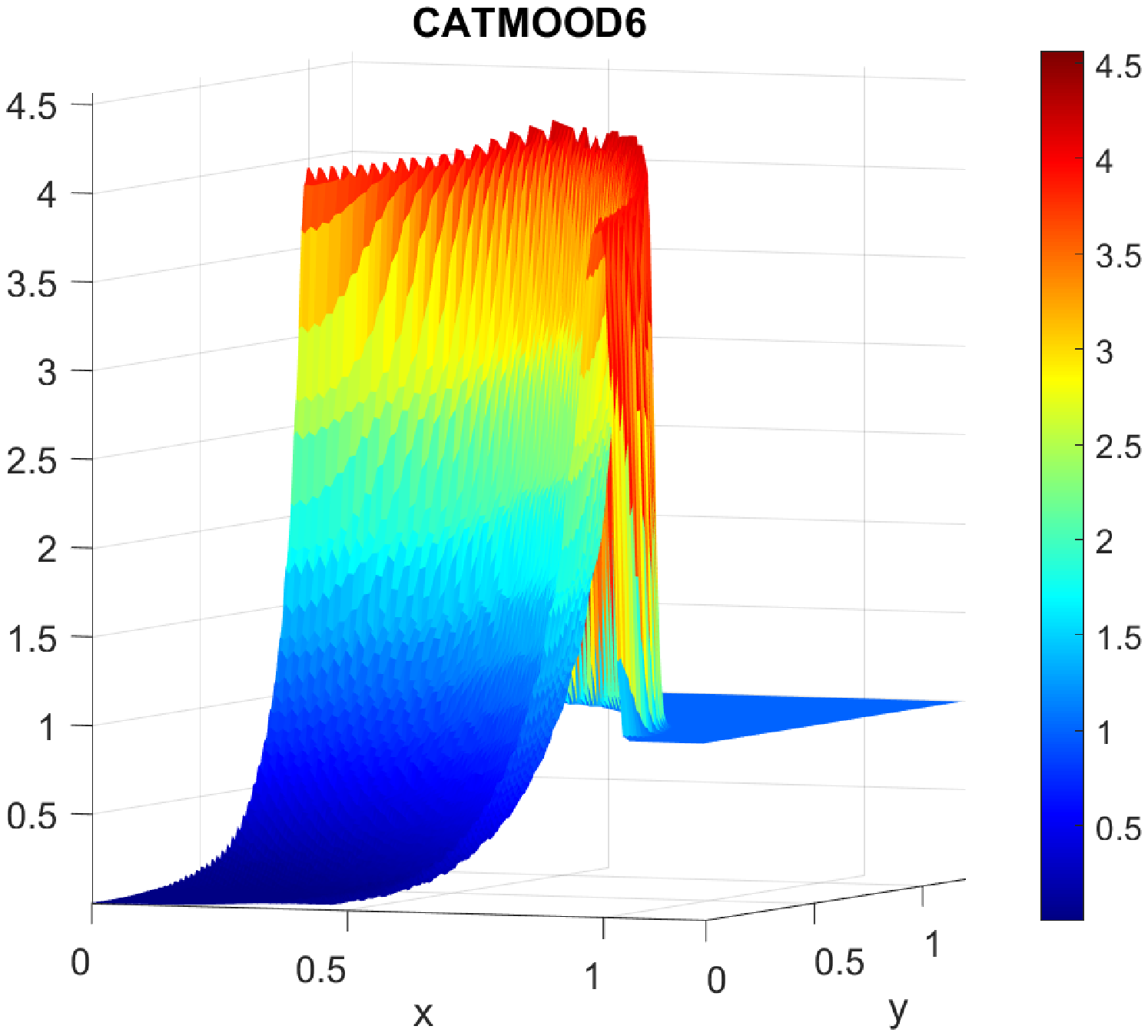}
         \vspace{0.2cm}
         \caption{Zoom on $[0,1.2]\times[0,1.2]$ of the numerical density for CATMOOD6.}
         \label{Sed:CATMOOD6_cut}
     \end{subfigure}
     \caption{Sedov blast wave \ref{ssec:Sedov}. Numerical density obtained with CATMOOD6 with HLLC as first order scheme on the interval $[-1.2,1.2]\times[-1.2,1.2]$ adopting a $200\times200$ mesh and CFL$=0.4$ at time $t=1$. Projection on the plane O-$xy$ (left); zoom on the interval $[0,1.2]\times[0,1.2]$ (right).}
     \label{Se:CATMOOD6}
\end{figure}

\begin{figure}[!ht]
     \centering
\begin{tabular}{cccc}
\hspace{-0.5cm}
\rotatebox{90}{\hspace{0.2cm} $100\times 100$ unif. cells} &
\hspace{-0.5cm}
\includegraphics[width=0.35\textwidth]{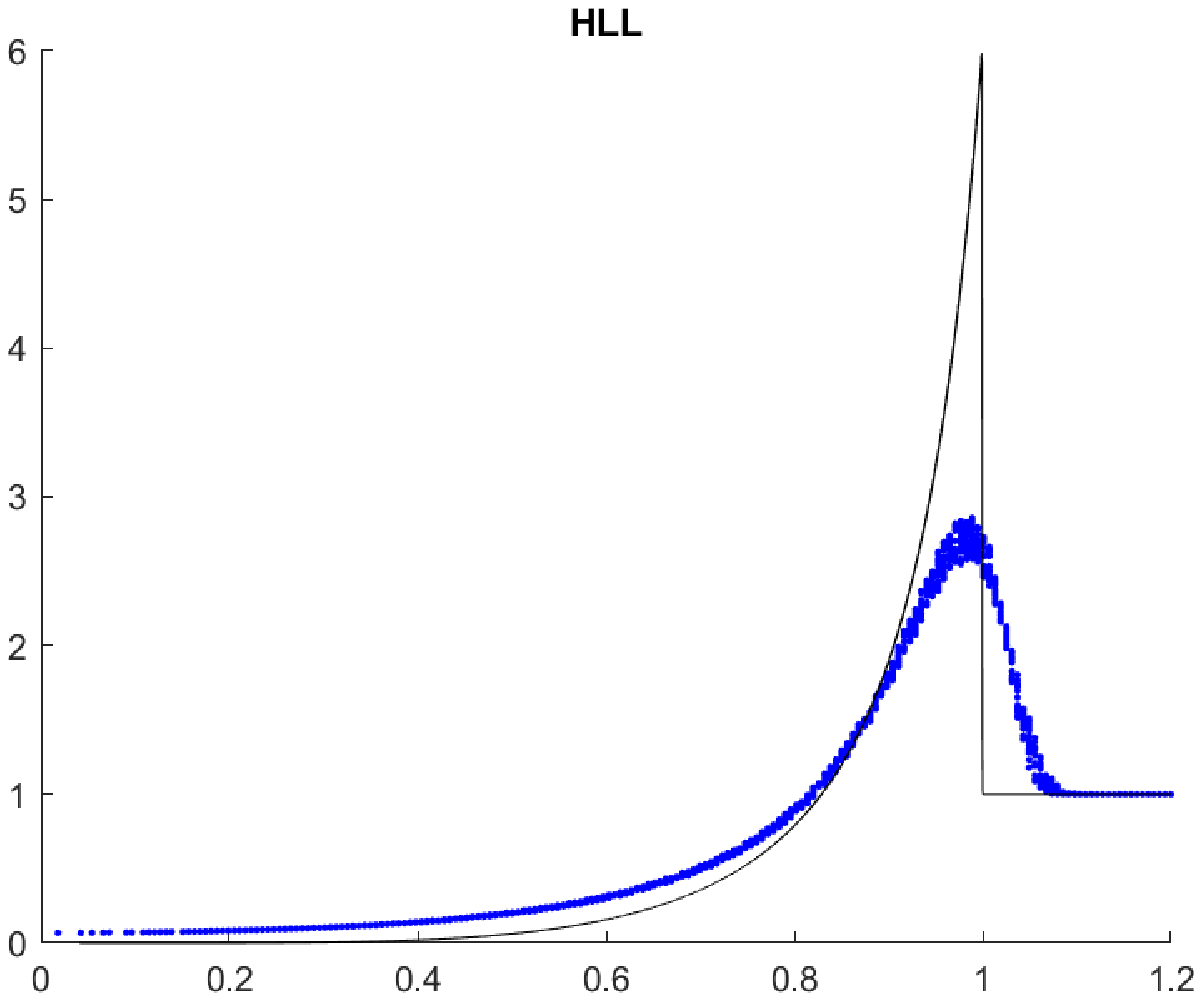}&
\hspace{-0.65cm}
\includegraphics[width=0.35\textwidth]{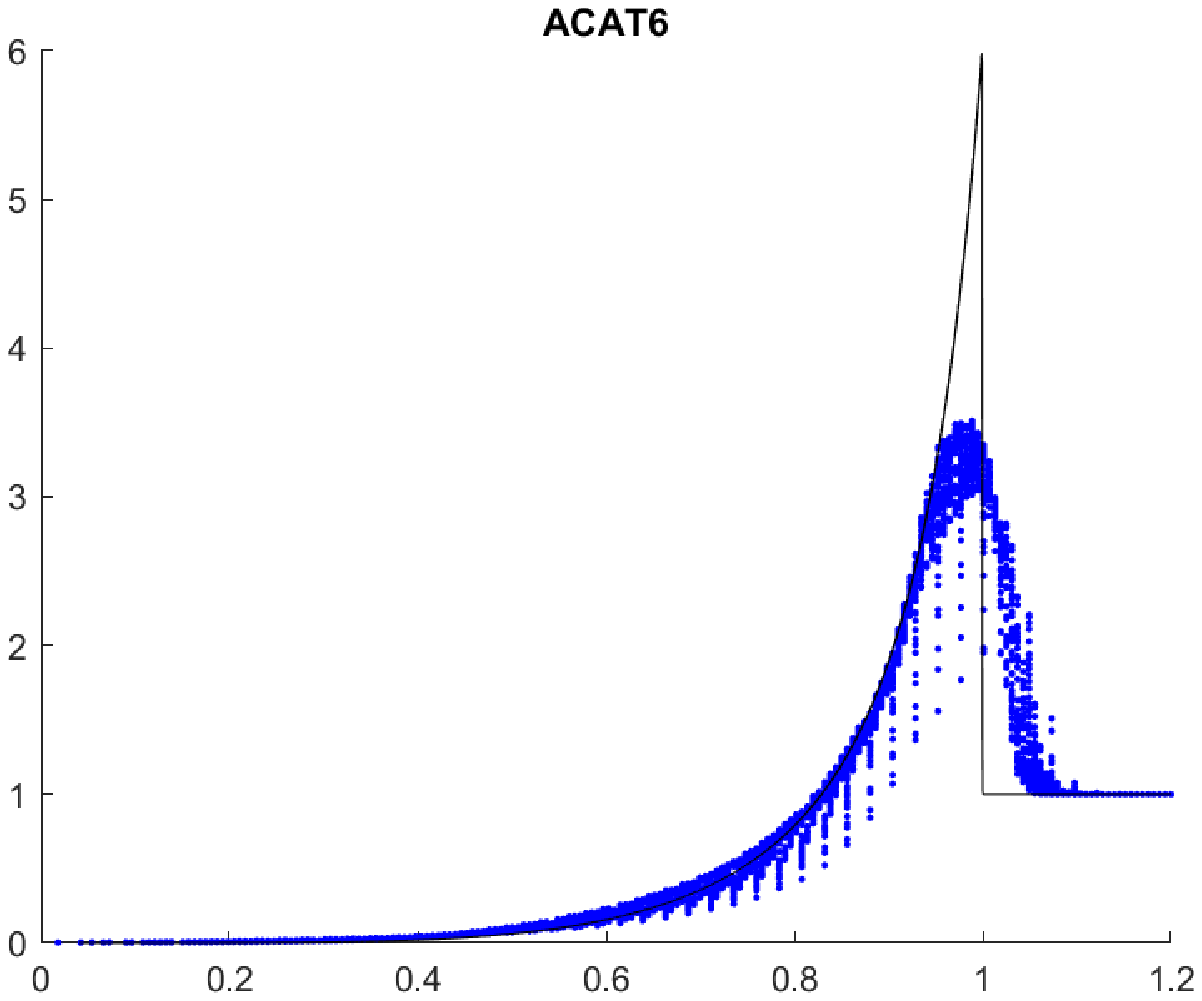}&
\hspace{-0.65cm}
\includegraphics[width=0.35\textwidth]{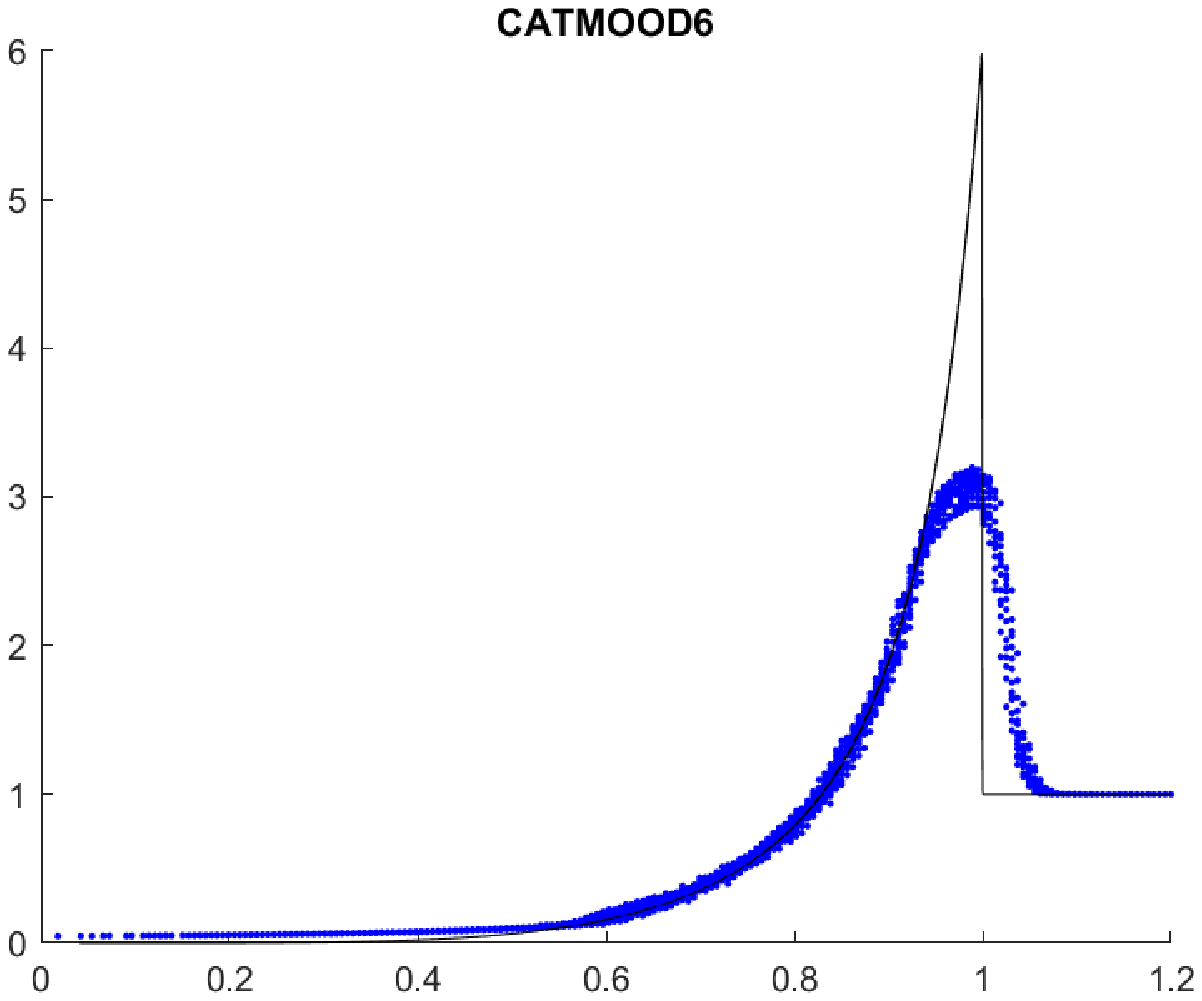} \\
\hspace{-0.5cm}
\rotatebox{90}{\hspace{0.2cm} $200\times 200$ unif. cells} &
\hspace{-0.5cm}
\includegraphics[width=0.35\textwidth]{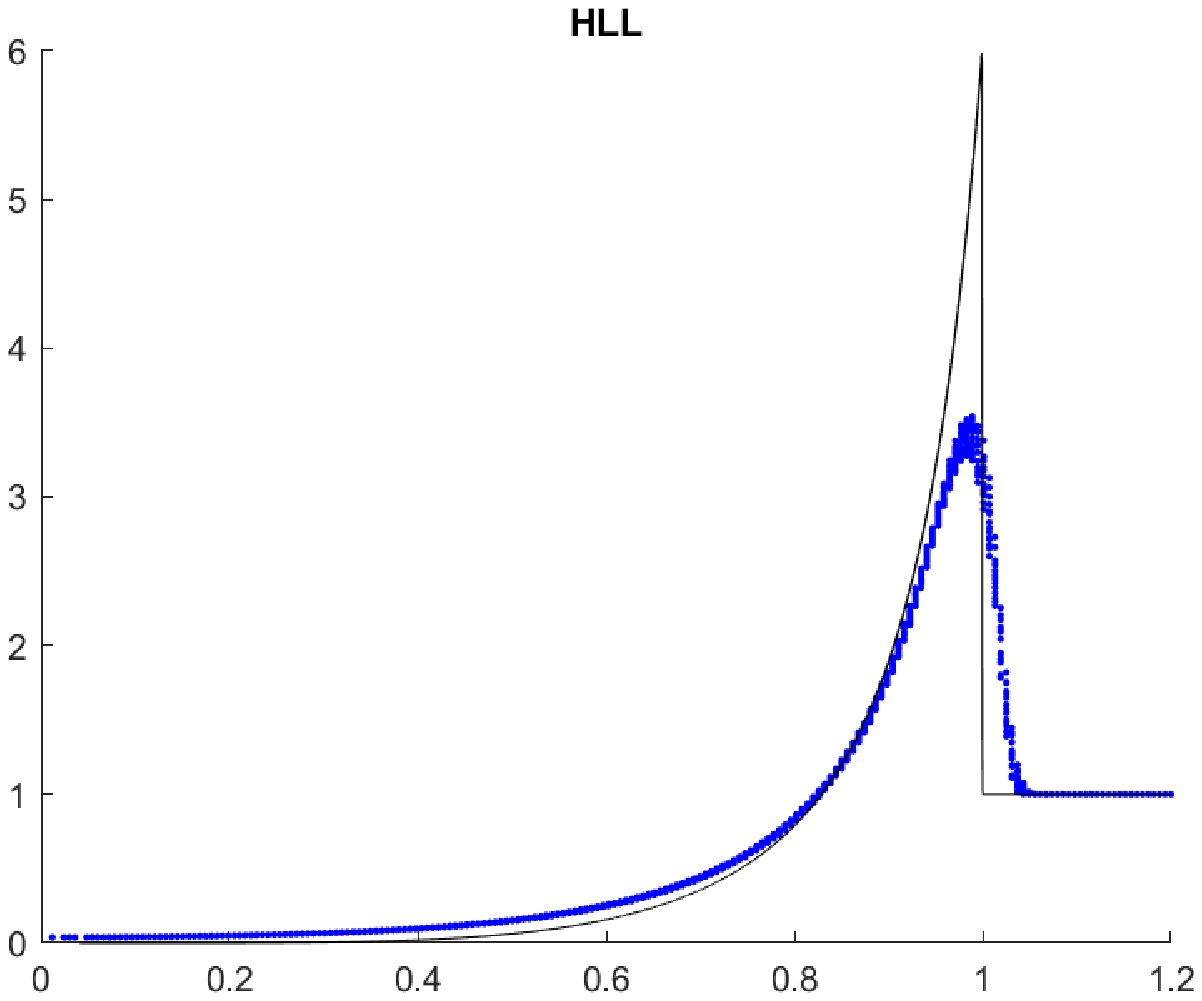}&
\hspace{-0.65cm}
\includegraphics[width=0.35\textwidth]{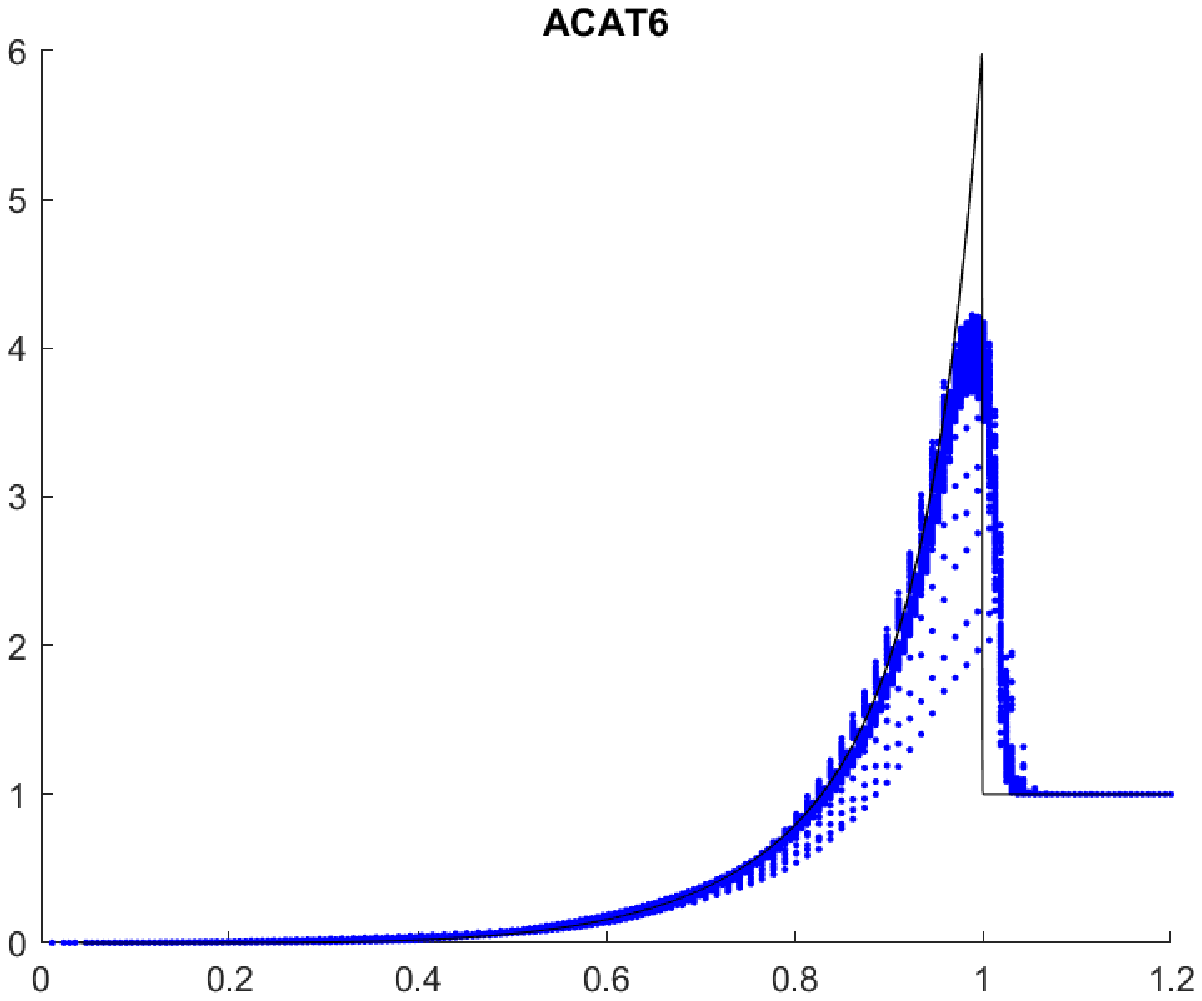}&
\hspace{-0.65cm}
\includegraphics[width=0.35\textwidth]{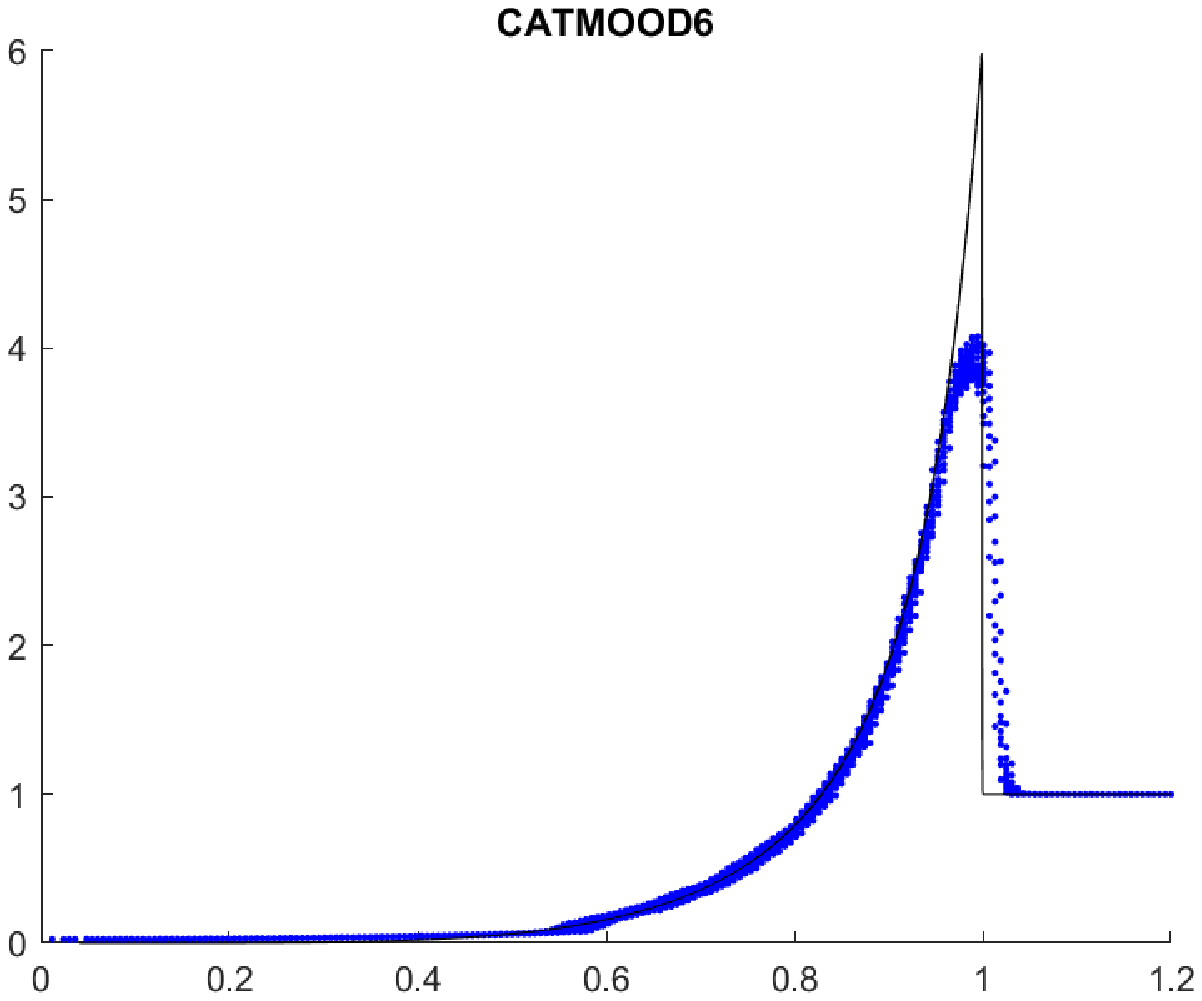}
\end{tabular}
     \caption{Sedov blast wave from section~\ref{ssec:Sedov} --- Scatter plots of the numerical density obtained with a 1st order scheme HLL, and 6th order ACAT6 and CATMOOD6 with HLL flux on the interval $[-1.2,1.2]\times[-1.2,1.2]$.
     Top: results for $100\times100$ uniform cells.
     Bottom: results for $200\times200$ uniform cells. In black is the exact solution.
     }
     \label{Se:CM}
\end{figure}
The 2D Sedov problem~\cite{Sedov} is a cylindrical symmetric explosion. The domain is given by $(x,y) \in [-1.2,1.2]^2 $ initially filled with perfect gas at rest such as $(\rho^0, u^0, v^0, p^0, \gamma) = (1, 0, 0, 10^{-13}, 1.4)$. 
A total energy of $E_{total}=0.244816$ is concentrated at the origin \cite{LOUBERE2005105}. 
This configuration corresponds to a point-wise symmetric explosion, for which a cylindrical shock front reaches radius $r=\sqrt{x^2+y^2}=1$ at $t_\text{final}=1$ with a density peak of $\rho=6$, see figure~\ref{Sed:CATMOOD6_cut} for an example of a quasi symmetric numerical solution. \\
Figure~\ref{Se:CM} presents the scatter plots of numerical density as a function of cell radius using, respectively for the first and second columns either $100\times100$ or $200\times200$ uniform cells. The tested schemes are  
the 1st order, ACAT6 and CATMOOD6 using HLL numerical flux. 
Notice that the unlimited CAT schemes fails in producing a numerical solution due to the presence of a strong shock wave.
%
Because the exact solution has a cylindrical symmetry (see the black line in figure~\ref{Se:CM}), all cells at the same radius should share the same numerical density, if the scheme exactly preserves the cylindrical symmetry.
The width of the spread/variance of numerical data somehow measures how good the scheme can preserve this symmetry.
As is visible in the figures, the refinement of the grid implies a reduction of the variance for any of the methods.
Moreover, the $6$th order schemes (ACAT6 and CATMOOD6) seem to produce a sharper shock wave and more accurate results. 
Importantly CATMOOD6 has a better cylindrical symmetry than ACAT6.
Of course no positivity issue is reported for any of these schemes.

\subsection{2D Riemann problems} \label{ssec:2DRiemann}
Here, we consider 2D Riemann problems on the interval $[-1,1]\times[-1,1]$ where zero Neumann conditions are applied, and  four initial conditions  are set respectively as in Table~\ref{IC:Riemann} \cite{RP_Schultz93}.
A $400 \times 400$ mesh is adopted for all schemes and configurations. The final time is set to $t_\text{final} = 0.3$ and the CFL to $0.4$.
We run four simulations using CATMOOD6 and the three first order schemes respectively with Rusanov, HLL and HLLC flux functions.
The 1st order scheme with Rusanov flux has been adopted as parachute scheme in the CATMOOD6 cascade.
In these 2D figures the schemes are respectively ordered from top-left to bottom-right.
In addition, in a different set of figures for CATMOOD6 schemes, we plot the percentage of cells updated with CAT6 (top), CAT2 (center) or 1st order Rusanov (bottom) schemes as a function of time-step.
\begin{table}[!ht]
    \centering
    \begin{tabular}{|l l | l l|}
    \hline 
    \multicolumn{4}{|c|}{\textbf{Configuration 3}} \\
    \hline 
   $\rho_2 = 0.5323$  &  $u_2 = 1.206$   &  $\rho_1 = 1.5$    & $u_1 = 0$     \\
   $v_2 = 0$          &  $p_2 = 0.3$     &  $v_1 = 0$         & $p_1 = 1.5$   \\  
    \hline 
    $\rho_3 = 0.138$  &  $u_2 = 1.206$   &  $\rho_4 = 0.5323$ & $u_ 4= 0$     \\
    $v_3 = 1.206$      &  $p_2 = 0.029$   &  $v_4 = 1.206$    & $p_4 = 0.3$  \\  \hline
    \hline 
    \multicolumn{4}{|c|}{\textbf{Configuration 6}} \\
    \hline 
   $\rho_2 = 2     $  &  $u_2 = 0.75 $   &  $\rho_1 = 1.5$    & $u_1 = 0.75$     \\
   $v_2 = 0.5$        &  $p_2 = 1  $     &  $v_1 = -0.5$      & $p_1 = 1$   \\  
    \hline 
   $\rho_3 = 1$       &  $u_2 = -0.75$   &  $\rho_4 = 3$      & $u_4 = -0.75$     \\
   $v_3 = 0.5$        &  $p_2 = 1$       &  $v_4 = -0.5$      & $p_4 =1$  \\  \hline
    \hline 
    \multicolumn{4}{|c|}{\textbf{Configuration 11}} \\
    \hline 
   $\rho_2 = 0.5313$  &  $u_2 = 0.8276$  &  $\rho_1 = 1$      & $u_1 = 0.1$     \\
   $v_2 = 0  $        &  $p_2 = 0.4$     &  $v_1 = 0$         & $p_1 = 1$   \\  
    \hline 
   $\rho_3 = 0.8$     &  $u_2 = 0.1$     &  $\rho_4 = 0.5313$ & $u_4 = 0.1$     \\
   $v_3 = 0$          &  $p_2 = 0.4$     &  $v_4 = 0$         & $p_4 =0.4$  \\  \hline
    \hline 
    \multicolumn{4}{|c|}{\textbf{Configuration 17}} \\
    \hline 
   $\rho_2 = 2     $  &  $u_2 = 0.   $   &  $\rho_1 = 1$      & $u_1 = 0$     \\
   $v_2 = -0.3$       &  $p_2 = 1  $     &  $v_1 = -0.4$      & $p_1 = 1$   \\  
    \hline 
   $\rho_3 = 1.0625$  &  $u_2 = 0$       &  $\rho_4 = 0.5197$ & $u_4 = 0$     \\
   $v_3 = 0.2145$     &  $p_2 = 0.4$     &  $v_4 = -1.1259$   & $p_4 =0.4$  \\  \hline
    \end{tabular} 
    \caption{2D Riemann problem initial conditions. }
    \label{IC:Riemann}
\end{table}

\begin{figure}[!ht]
    \centering
    \begin{subfigure}[b]{0.328\textwidth}
         \centering
         \includegraphics[width=1.3\textwidth]{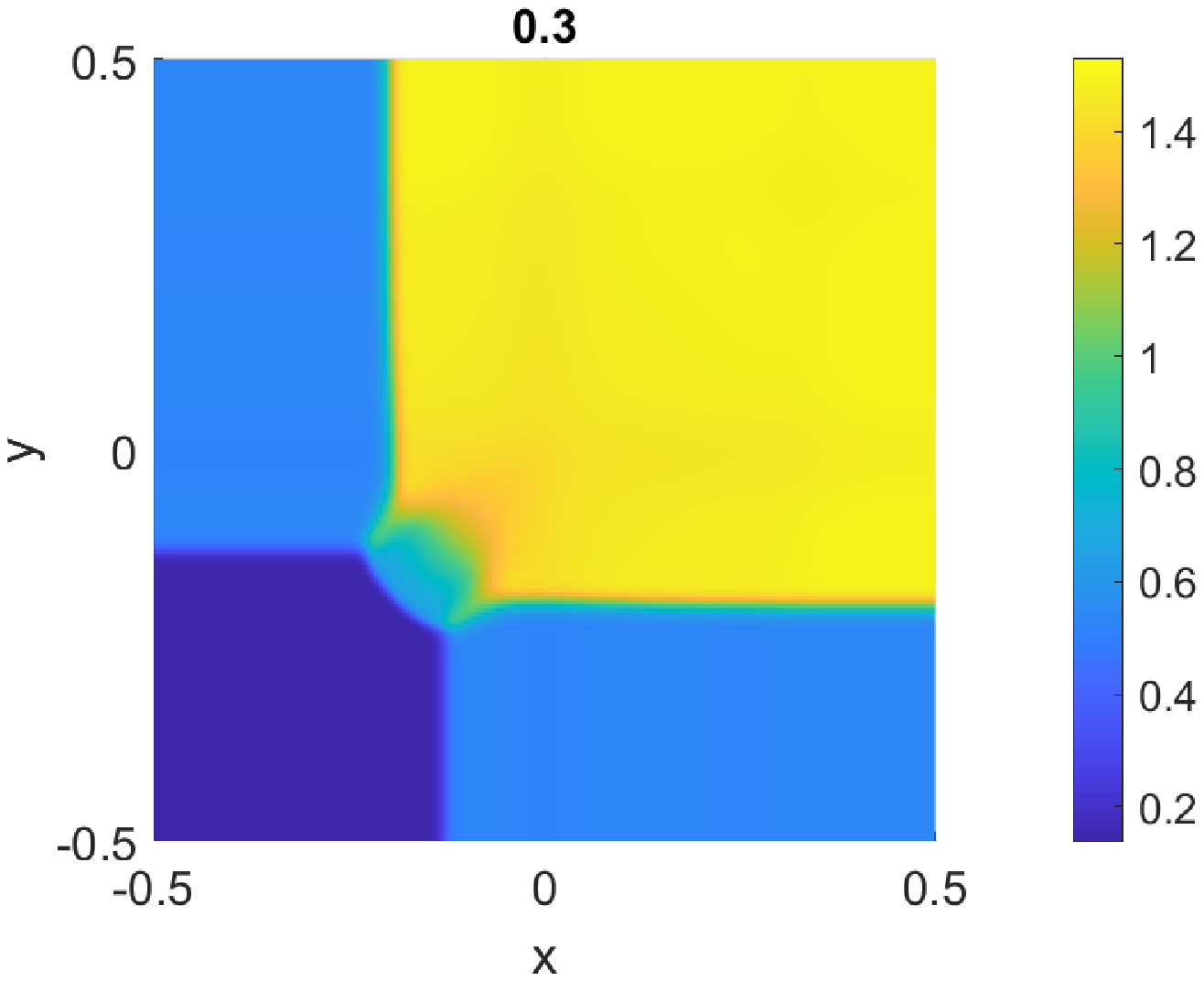}
         \caption{Rusanov flux}
         \label{fig:2D_R3_RU}
     \end{subfigure}
     \hfill
     \begin{subfigure}[b]{0.328\textwidth}
    \hspace{-1.cm}
         \centering
         \includegraphics[width=1.3\textwidth]{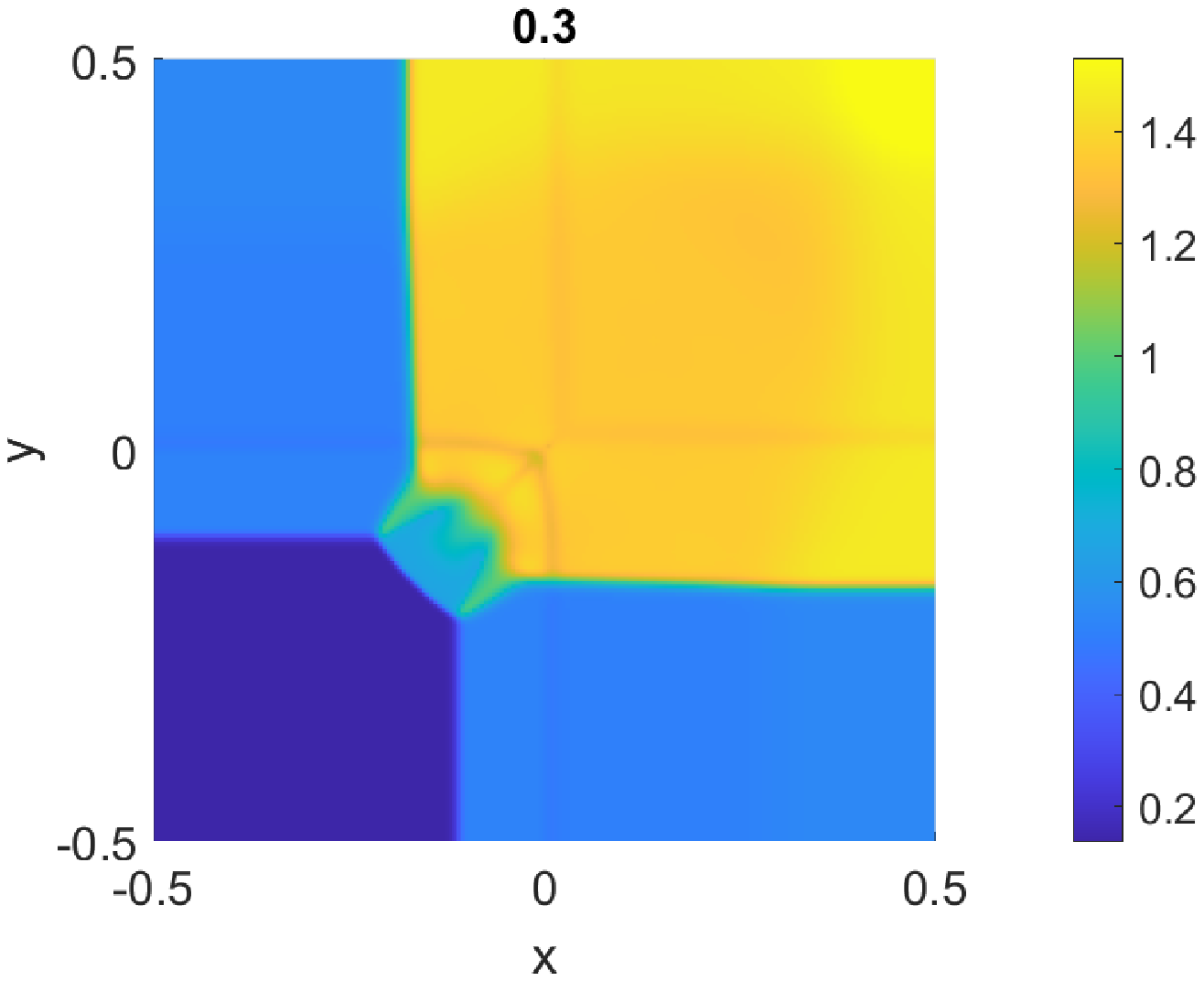}
         \caption{HLLC}
         \label{fig:2D_R3_HLLC}
     \end{subfigure}
     \hfill
     \begin{subfigure}[b]{0.328\textwidth}
    \hspace{-1.cm}
         \centering
         \includegraphics[width=1.3\textwidth]{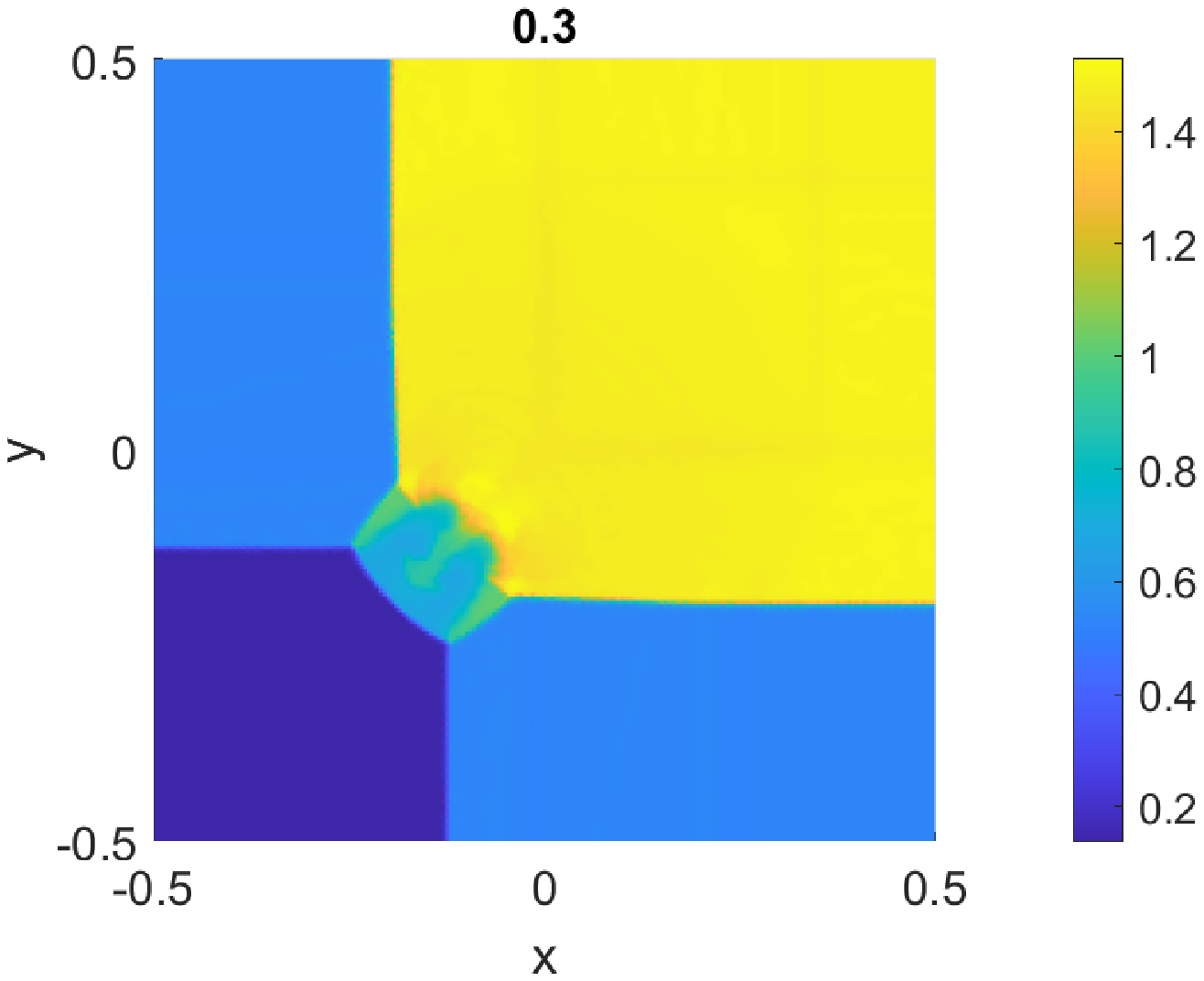}
         \caption{CATMOOD6}
         \label{fig:2D_R3_CATMOOD6}
     \end{subfigure}
    \caption{2D Riemann problem from section~\ref{ssec:2DRiemann} with initial condition configuration 3. Zoom of the numerical solution for density on the interval $[-1,1]\times[-1,1]$ adopting a mesh of $400 \times 400-$cells and CFL$=0.4$. Rusanov-flux (a); HLLC (b); and CATMOOD6 with Rusanov flux for the first order method (c).}
   \label{fig:2D_R_C3}
\end{figure}

\begin{figure}[!ht]
    \centering
    \includegraphics[width=0.99\textwidth]{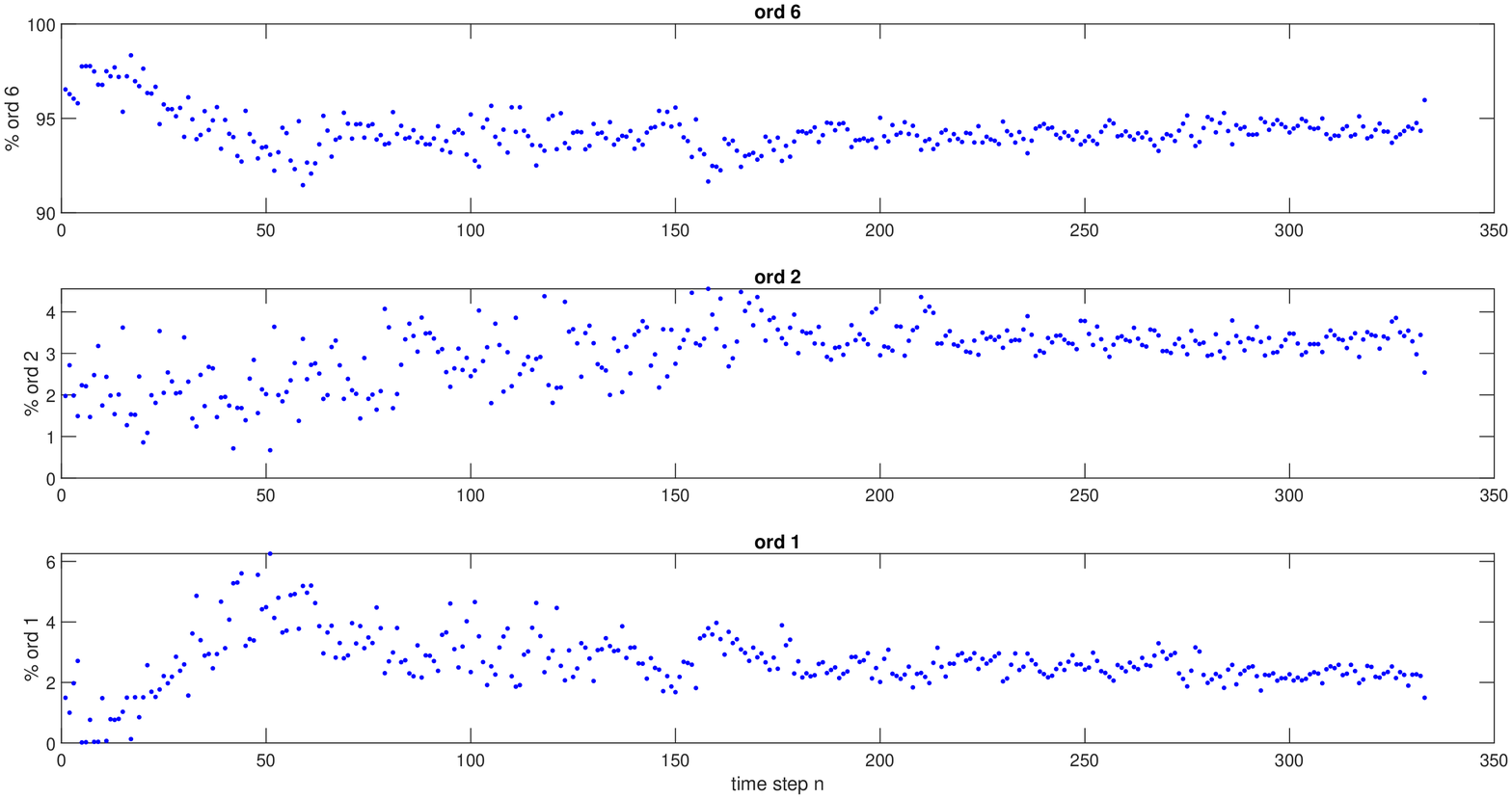}
    \vspace{-0.3cm}
    \caption{2D Riemann problem from section~\ref{ssec:2DRiemann} with initial condition configuration 3.Percentage of cells updated by CAT6 (top), by CAT2 (center), and by Rusanov (bottom).
    }
   \label{fig:2D_R_C3_ord}
\end{figure}

\paragraph{Configuration 3}
Figure~\ref{fig:2D_R_C3} shows a zoom of the numerical densities for the configuration 3 (Table~\ref{IC:Riemann}).
All the schemes capture the same global solution. The first order methods are  diffusive even if HLLC scheme performs better. Contrarily, the high order CATMOOD6 clearly improves the sharpness of the shear layers and contacts. 
Figure~\ref{fig:2D_R_C3_ord} exhibits the percentage of cells using CAT6 (top), CAT2 (center) and Rusanov scheme (bottom). We can observe that on average $95\%$ of the cells are updated with $6$th order of accuracy, about $3-4\%$ with $2$nd order or $1$st order. One observes that only few cells demand some limiting, and are sent back to $t^n$ for re-computation by the MOOD approach.
Consequently, the cost of CATMOOD6 is mainly the one of CAT6. 
\begin{figure}[!ht]
    \centering
    \begin{subfigure}[b]{0.328\textwidth}
         \centering
         \includegraphics[width=1.3\textwidth]{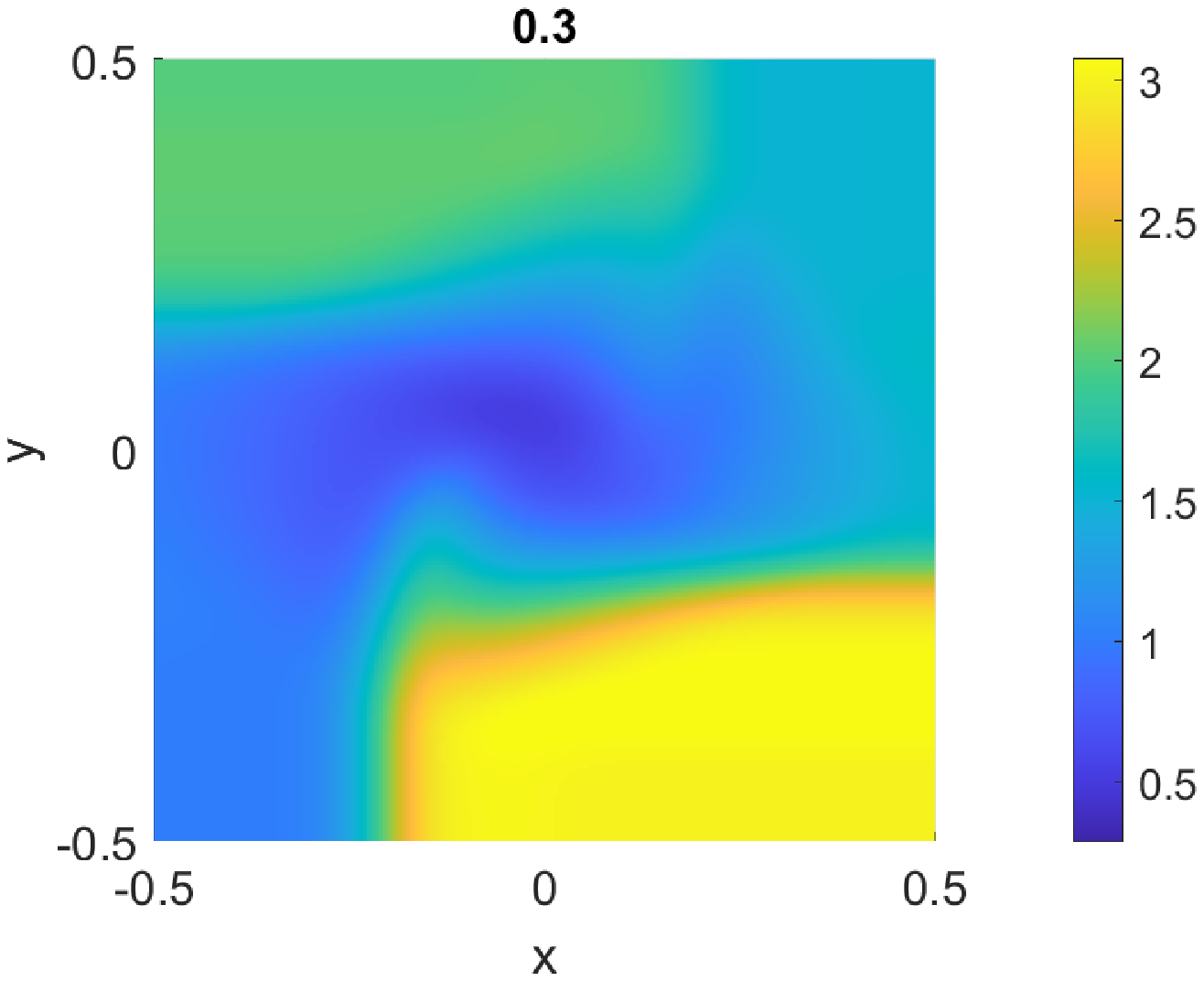}
         \caption{Rusanov flux}
         \label{fig:2D_R6_RU}
     \end{subfigure}
     \hfill
     \begin{subfigure}[b]{0.328\textwidth}
    \hspace{-1.cm}
         \centering
         \includegraphics[width=1.3\textwidth]{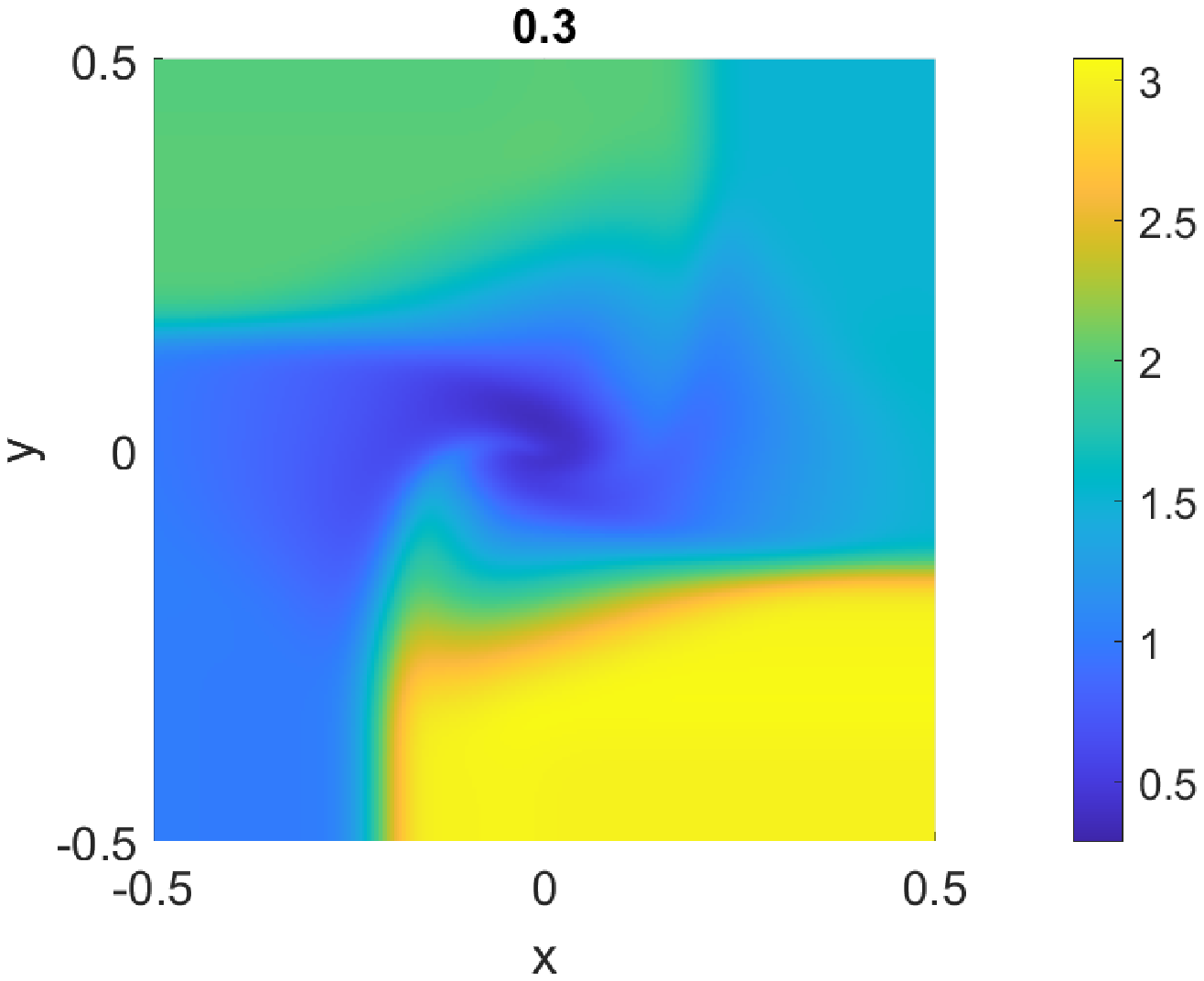}
         \caption{HLLC}
         \label{fig:2D_R6_HLLC}
     \end{subfigure}
     \hfill
     \begin{subfigure}[b]{0.328\textwidth}
    \hspace{-1.cm}
         \centering
         \includegraphics[width=1.3\textwidth]{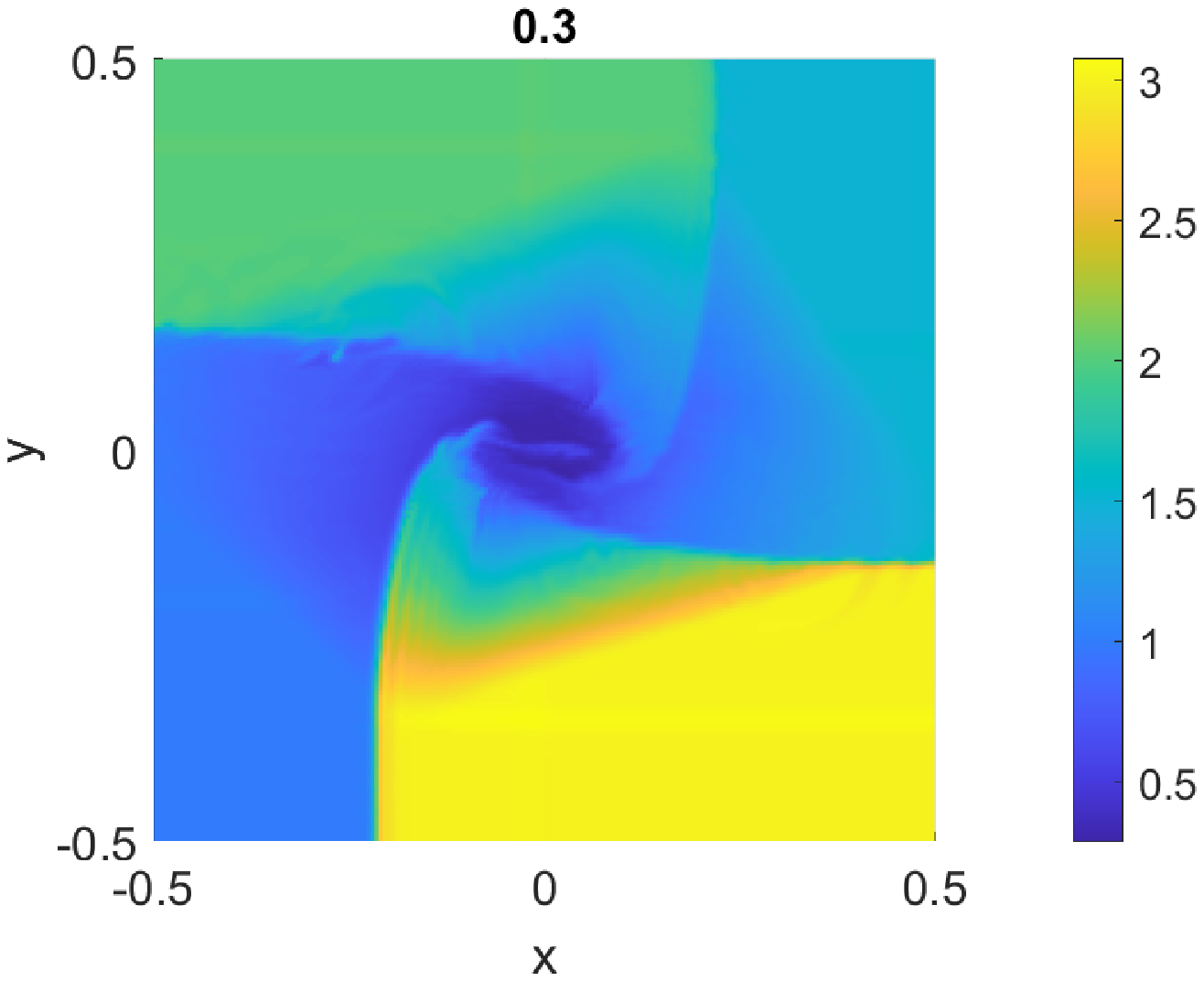}
         \caption{CATMOOD6}
         \label{fig:2D_R6_CATMOOD6}
     \end{subfigure}
    \caption{2D Riemann problem from section~\ref{ssec:2DRiemann} with initial condition configuration 6. Zoom of the numerical solution for density on the interval $[-1,1]\times[-1,1]$ adopting a mesh of $400 \times 400-$cells and CFL$=0.4$. Rusanov-flux (a); HLLC (b); and CATMOOD6 with HLLC for the first order method (c).}
   \label{fig:2D_R_C6}
\end{figure}
\begin{figure}[!ht]
    \centering
    \includegraphics[width=0.99\textwidth]{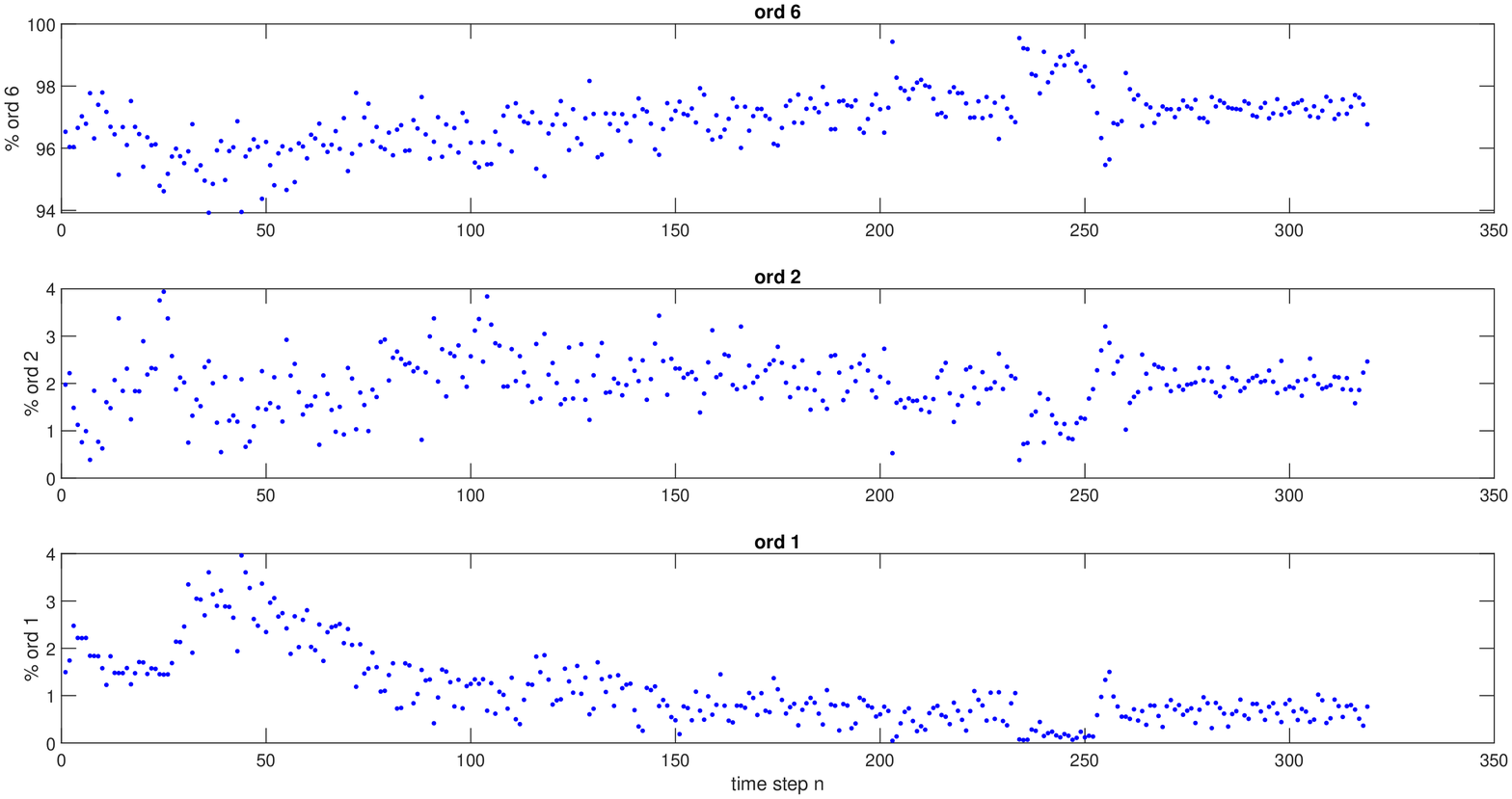}
    \vspace{-0.3cm}
    \caption{2D Riemann problem from section~\ref{ssec:2DRiemann} with initial condition configuration 6. Percentage of cells updated by CAT6 (top), by CAT2 (center), and by HLLC (bottom).}
   \label{fig:2D_R_C6_ord}
\end{figure}

\paragraph{Configuration 6}
Likewise for configuration 3, we plot the results in Figure~\ref{fig:2D_R_C6} and \ref{fig:2D_R_C6_ord}.
We again observe that the CATMOOD6 results are far more sharper than the 1st order schemes, with HLLC being the less dissipative scheme of the 1st order ones.
Notice that the shear layers along the curves are Kelvin-Helmholtz unstable, so that it seems normal to see the occurrences of small vortices. This is a numerical evidence of low dissipation of CATMOOD6.
The percentage of troubled cells is of the order of $3-4\%$ on averaged, again showing that CATMOOD6 has globally the same cost than CAT6 plus $10\%$ of the cost of a first order scheme.
\begin{figure}[!ht]
    \centering
    \begin{subfigure}[b]{0.328\textwidth}
         \centering
         \includegraphics[width=1.3\textwidth]{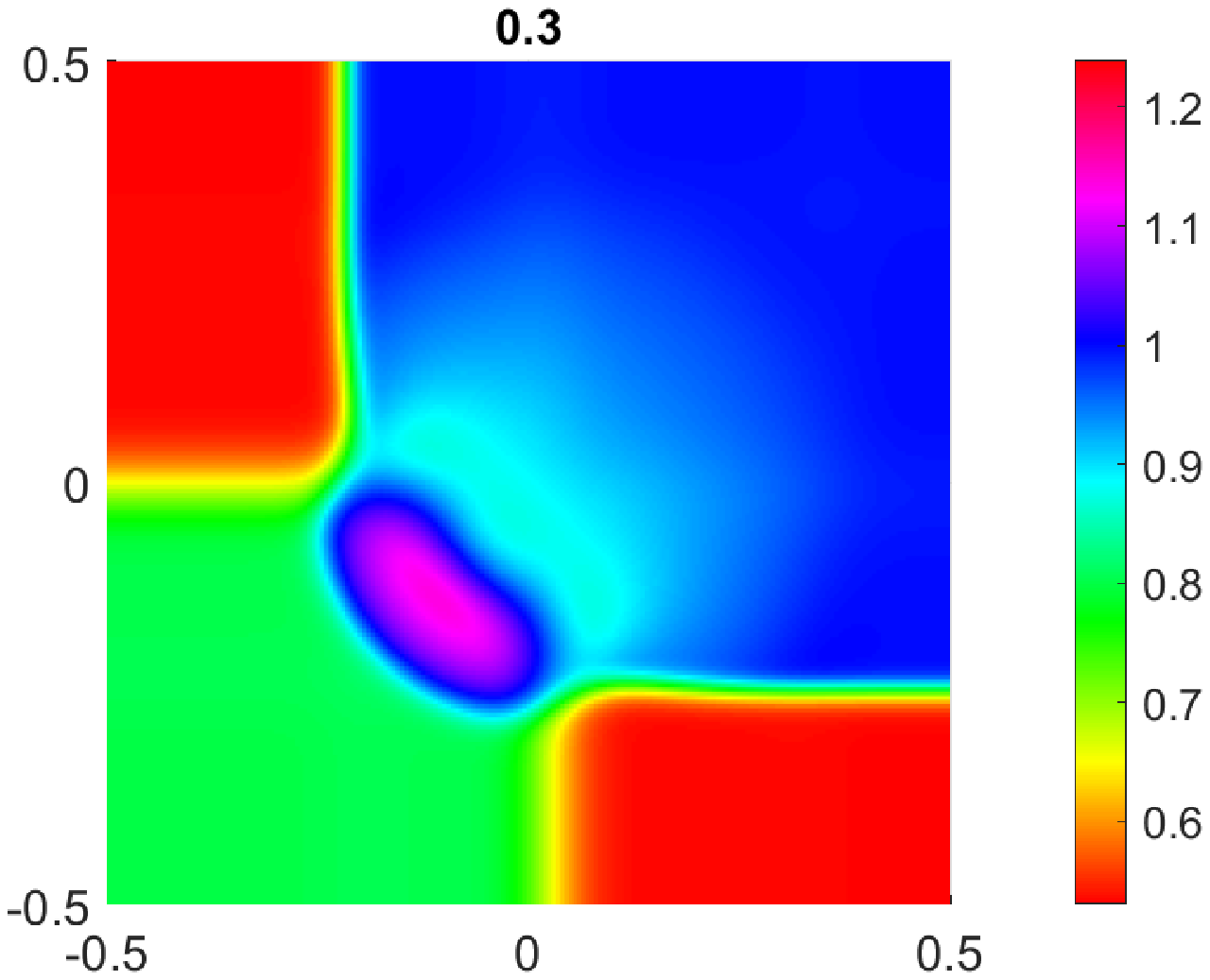}
         \caption{Rusanov flux}
         \label{fig:2D_R11_RU}
     \end{subfigure}
     \hfill
     \begin{subfigure}[b]{0.328\textwidth}
        \hspace{-1cm}
         \centering
         \includegraphics[width=1.3\textwidth]{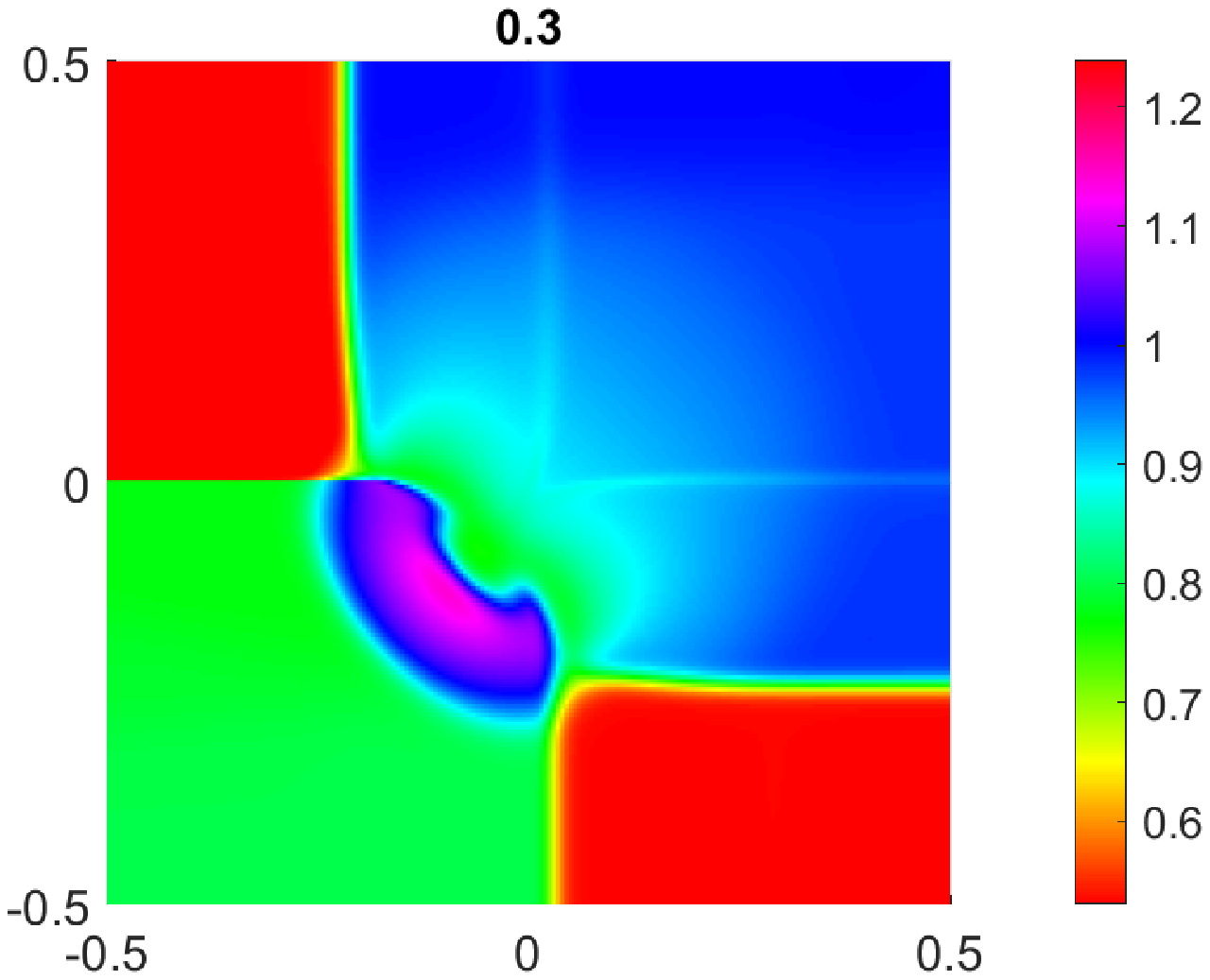}
         \caption{HLLC}
         \label{fig:2D_R11_HLLC}
     \end{subfigure}
     \hfill
     \begin{subfigure}[b]{0.328\textwidth}
        \hspace{-1cm}
         \centering
         \includegraphics[width=1.3\textwidth]{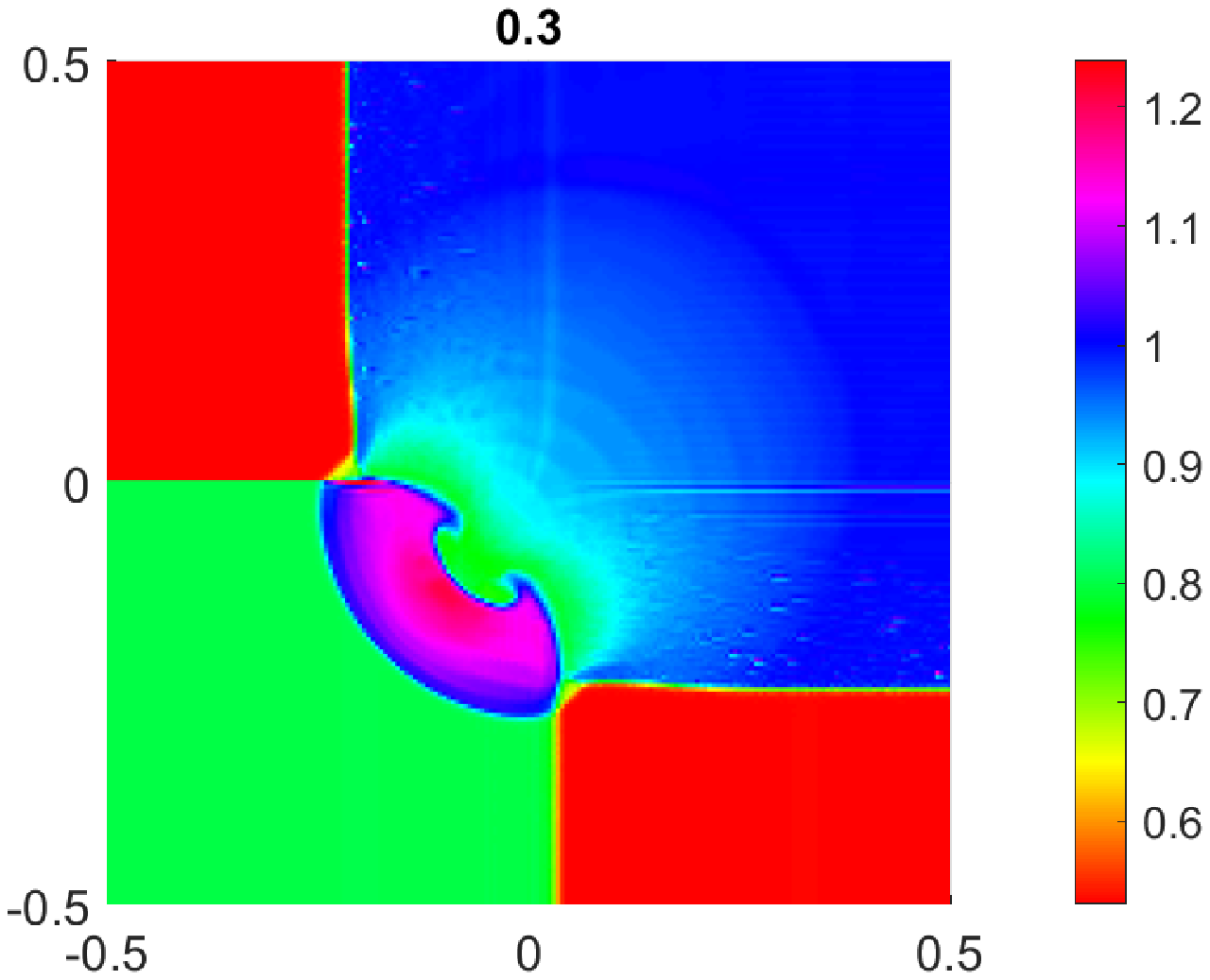}
         \caption{CATMOOD6}
         \label{fig:2D_R11_CATMOOD6}
     \end{subfigure}
    \caption{2D Riemann problem from section~\ref{ssec:2DRiemann} with initial condition configuration 11. Zoom of the numerical solution for density on the interval $[-1,1]\times[-1,1]$ adopting a mesh of $400 \times 400-$cells and CFL$=0.4$. Rusanov-flux (a); HLLC (b); and CATMOOD6 with HLLC for the first order method (c).}
   \label{fig:2D_R_C11}
\end{figure}
\begin{figure}[!ht]
    \centering
    \includegraphics[width=0.99\textwidth]{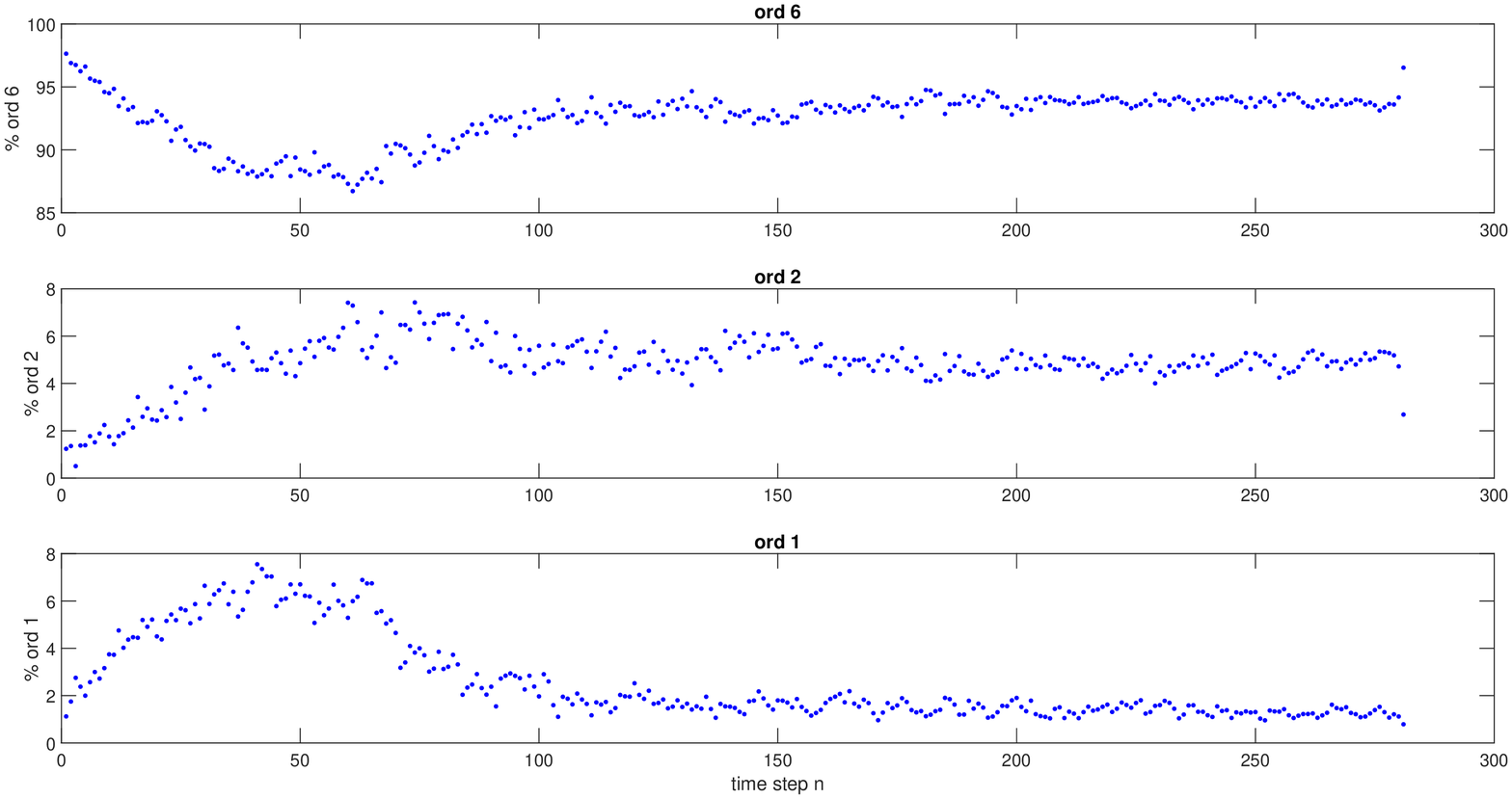}
    \vspace{-0.3cm}
    \caption{2D Riemann problem from section~\ref{ssec:2DRiemann} with initial condition configuration 11.  Percentage of cells updated by CAT6 (top), by CAT2 (center), and by HLLC (bottom).}
   \label{fig:2D_R_C11_ord}
\end{figure}

\paragraph{Configuration 11}
The numerical solutions for configuration 11 are gathered in Figure~\ref{fig:2D_R_C11} and the percentages of troubled cells in \ref{fig:2D_R_C11_ord}.
It is important to notice that the wave are visibly sharper with the high-order CATMOOD6 as expected. The complex internal pattern is also better captured, avoiding the excessive diffusion that can be seen when any 1st order scheme is employed. 

\begin{figure}[!ht]
    \centering
    \hspace{-0.5cm}
     \begin{subfigure}[b]{0.495\textwidth}
         \centering
        \includegraphics[width=1.0\textwidth]{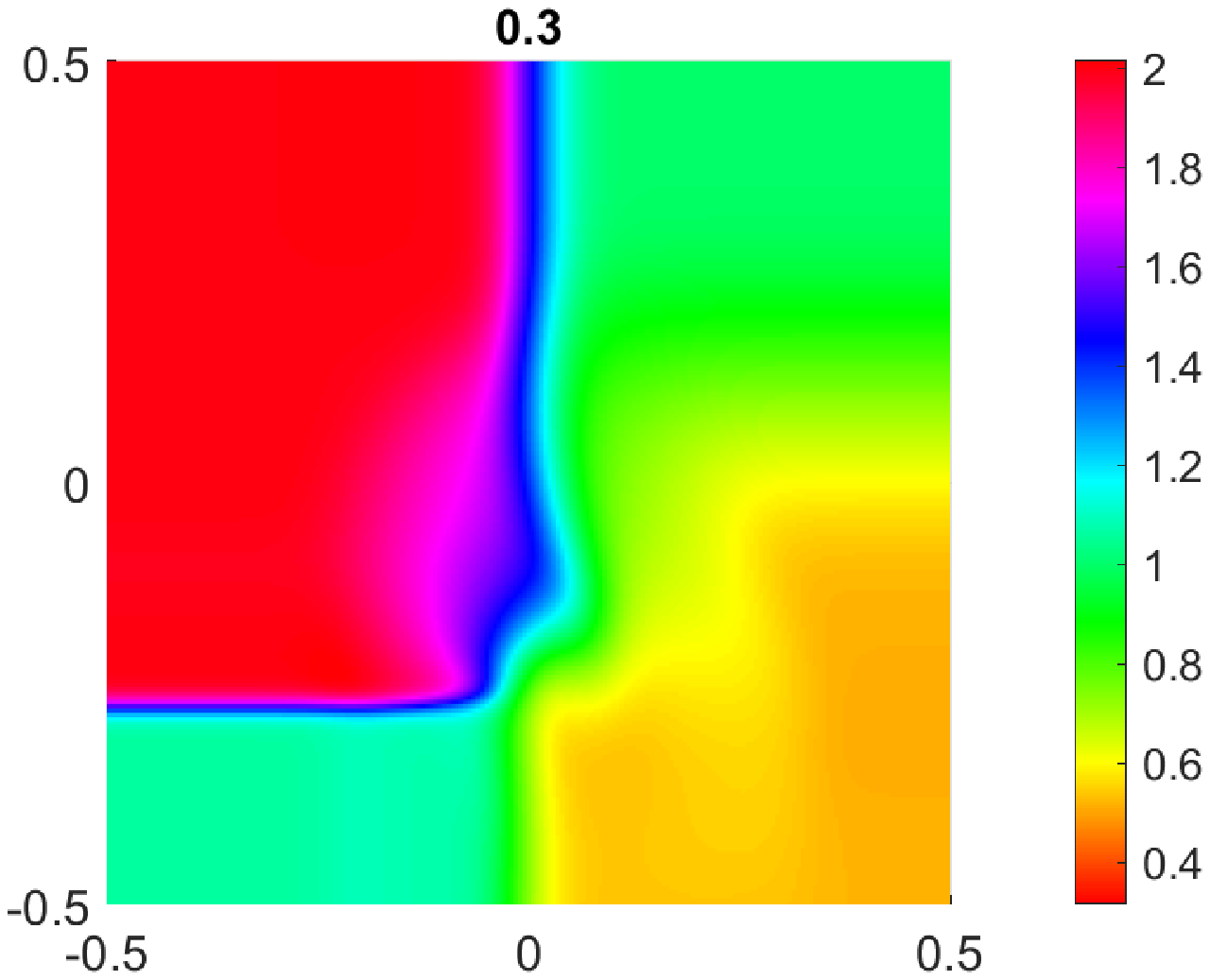}
         \caption{Rusanov flux}
         \label{fig:2D_R17_RU}
     \end{subfigure}
     \hfill
     \begin{subfigure}[b]{0.495\textwidth}
         \centering
         \includegraphics[width=1.0\textwidth]{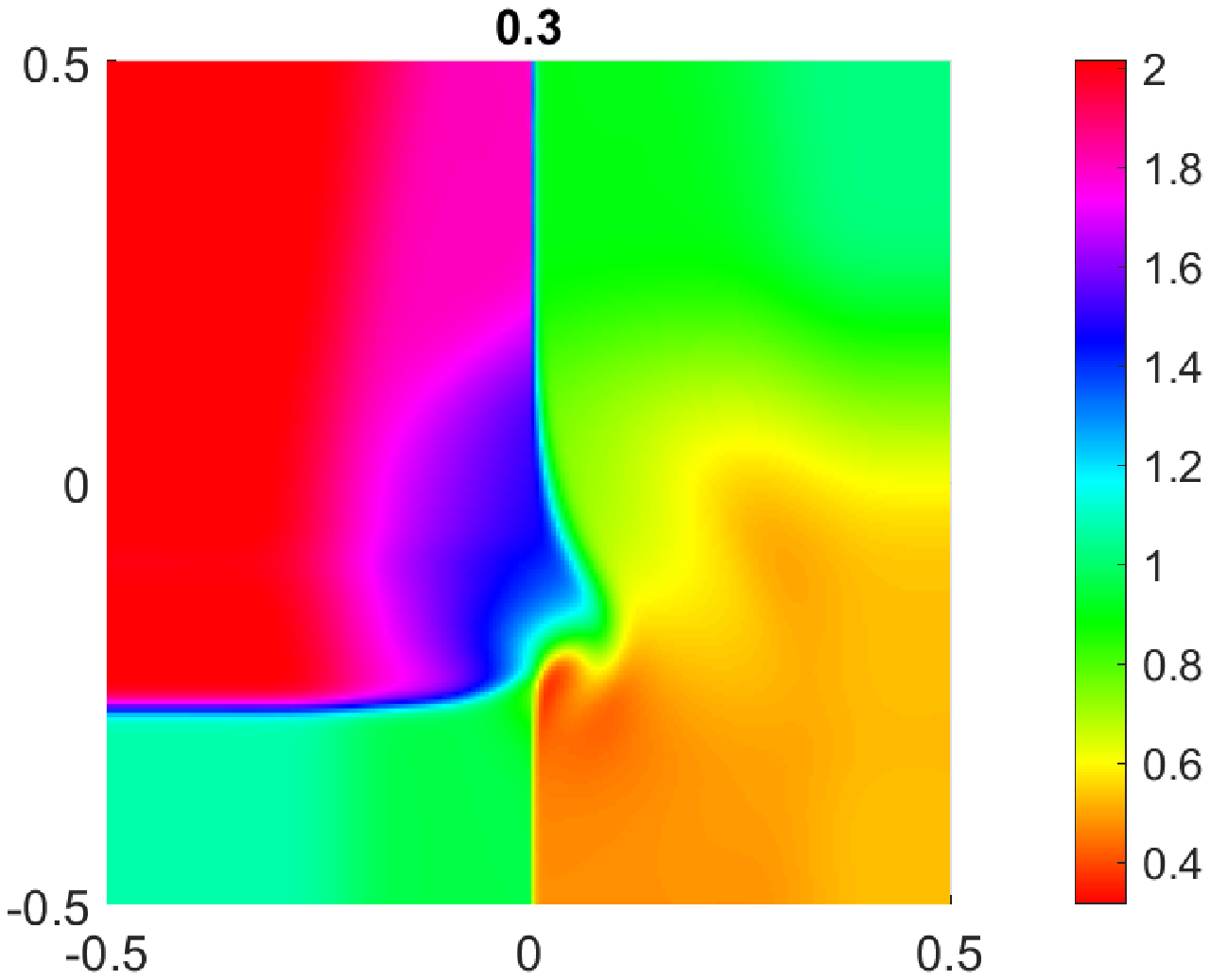}
         \caption{HLLC}
         \label{fig:2D_R17_HLLC}
     \end{subfigure}
     \\
    \hspace{-0.5cm}
     \begin{subfigure}[b]{0.495\textwidth}
         \centering
         \includegraphics[width=1.0\textwidth]{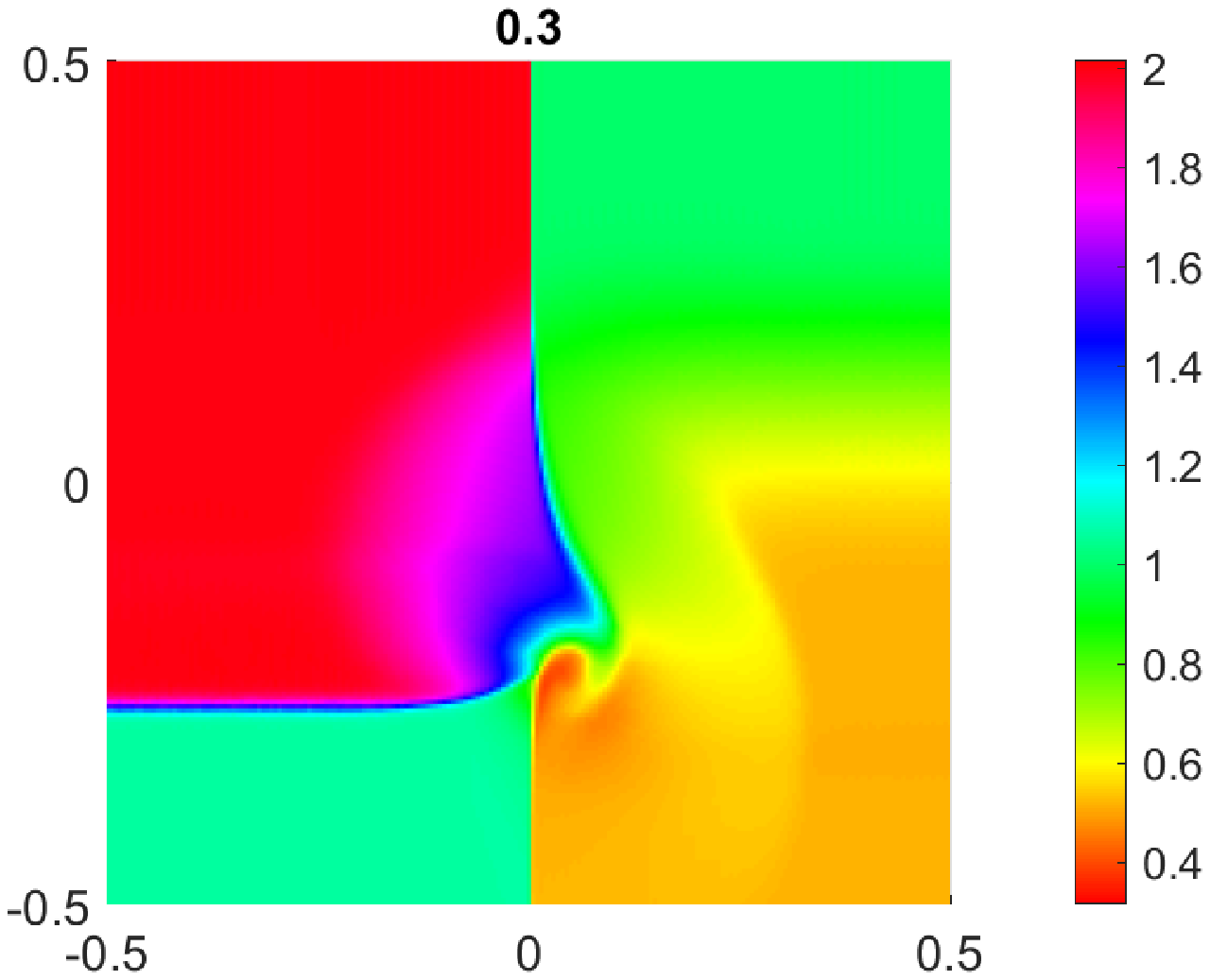}
         \caption{ACAT6}
         \label{fig:2D_R17_ACAT6}
     \end{subfigure}
     \hfill
     \begin{subfigure}[b]{0.495\textwidth}
         \centering
         \includegraphics[width=1.0\textwidth]{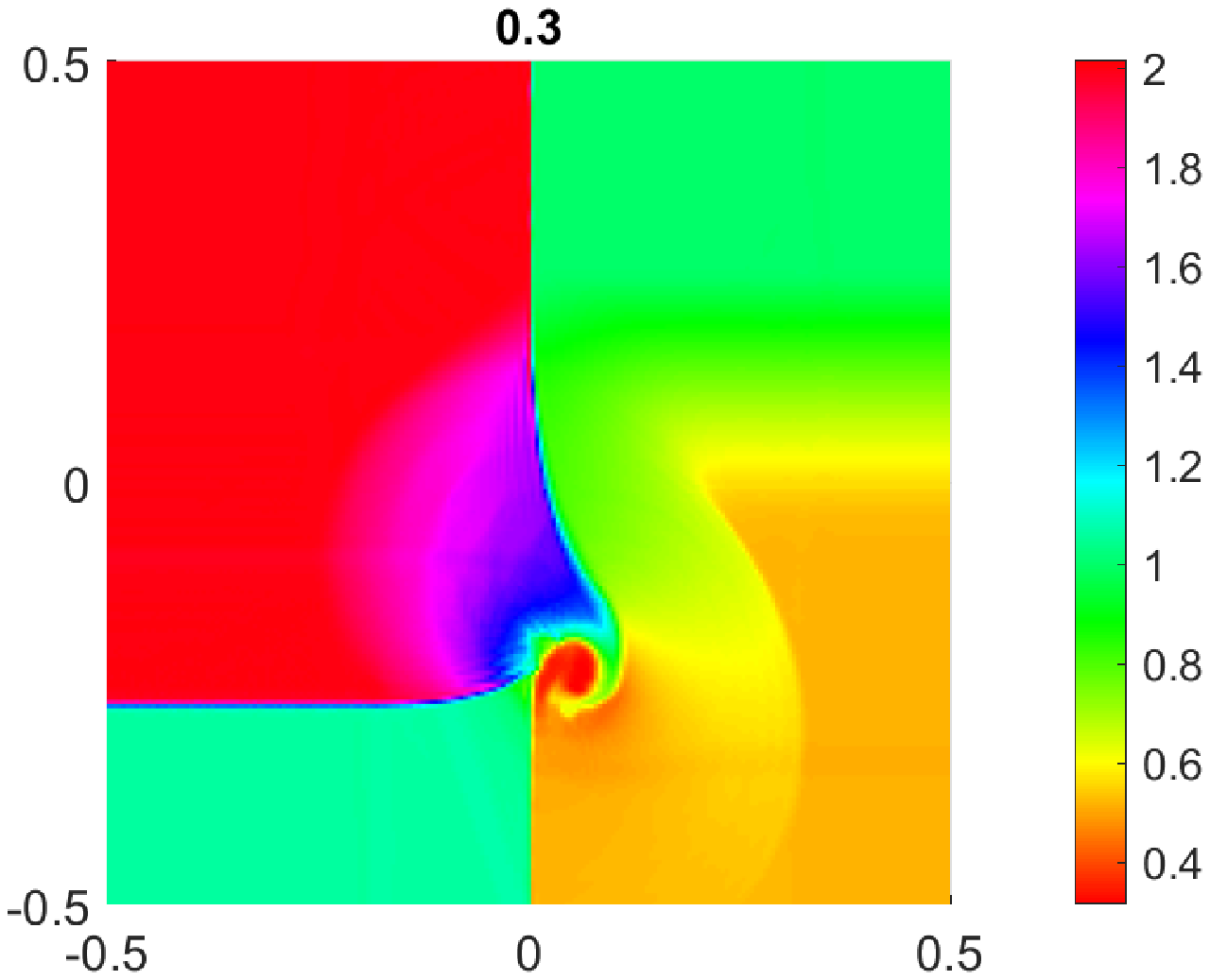}
         \caption{CATMOOD6}
         \label{fig:2D_R17_CATMOOD6}
     \end{subfigure}
    \caption{2D Riemann problem from section~\ref{ssec:2DRiemann} with initial condition configuration 17. Zoom of the numerical solution for density on the interval $[-1,1]\times[-1,1]$ adopting a mesh of $400 \times 400-$cells and CFL$=0.4$. Rusanov-flux (a); HLLC (b); ACAT6 (c); and CATMOOD6 (d).}
   \label{fig:2D_R_C17}
\end{figure}
\begin{figure}[!ht]
    \centering
    \includegraphics[width=0.99\textwidth]{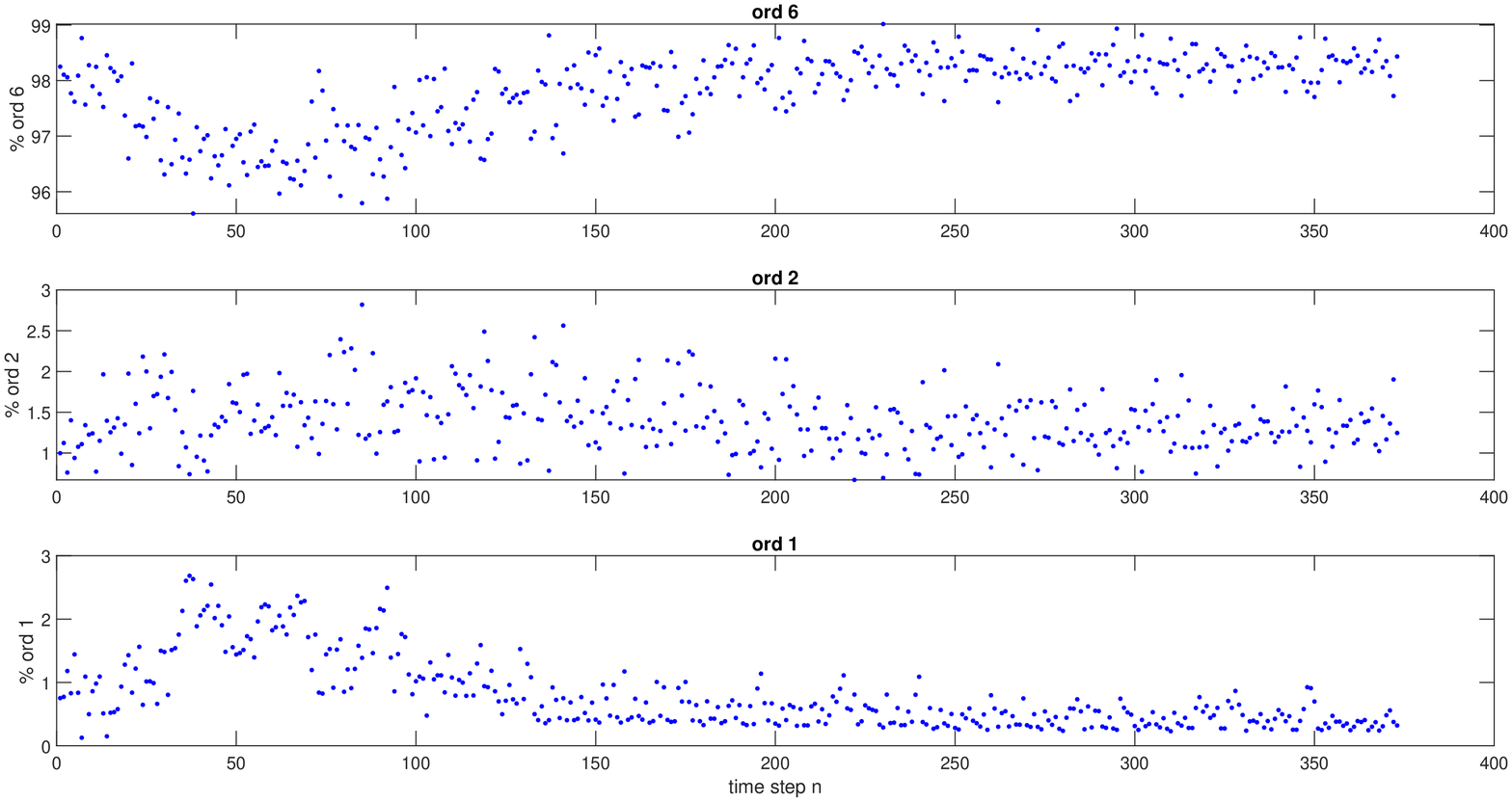}
    \vspace{-0.3cm}
    \caption{2D Riemann problem from section~\ref{ssec:2DRiemann} with initial condition configuration 17. Percentage of cells updated by CAT6 (top), by CAT2 (center), and by HLLC (bottom).}
   \label{fig:2D_R_C17_ord}
\end{figure}

\paragraph{Configuration 17}
Figure~\ref{fig:2D_R_C17} and 
Figure~\ref{fig:2D_R_C17_ord} presents the numerical results.
This test presents a small vortex type of structure along with first and secondary waves emanating from the quadruple point.
We can easily observe that CATMOOD6 can capture the vortex and the waves. Contrarily a diffusive scheme would require many more cells to get to this accuracy.
Here the HLLC flux is adopted in the first order scheme of the CATMOOD6 scheme.
The percentage of troubled cells are of the same order than previously, namely $97-98\%$ are updated with $6$th order of accuracy and $2\%$ and $1\%$ with $2$nd or $1$st orders for CATMOOD6. 
The extra cost compared to an unlimited CAT6 scheme is therefore acceptable. 
Table~\ref{CPU:2D_Riemann17} displays, in its first row, the computational cost expressed in seconds for the Configuration 17 concerning Rusanov, HLLC, ACAT6, and CATMOOD6 schemes. Meanwhile, in the second row are presented the ratio between the computational costs with respect to the Rusanov scheme cost. 
Figure~\ref{fig:2D_R_C17} clearly demonstrates the remarkable superiority of the adaptive \textit{a posteriori} approach CAT+MOOD over the \textit{a priori} technique ACAT. 
Despite exhibiting similar computational costs the CATMOOD6 numerical solution significantly outperforms the ACAT6 one in terms of accuracy. 
\begin{table}[!ht]
    \centering
    \numerikNine
    \begin{tabular}{|c||c| c| c| c|}
    \hline 
    & \textbf{Rusanov } & \textbf{HLLC} &  \textbf{ACAT6}   &  \textbf{CATMOOD6} \\
    \hline 
    \textbf{CPU (s)} & 25.28 & 34.11   & 14563.40  & 12871.25   \\ \hline
    \textbf{Ratio} & 1 & 1.35  & 576.08 &  509.15 \\
    \hline
    \end{tabular} 
    \caption{ 2D Riemann problem configuration 17. First row: CPU computational costs expressed in seconds. Second row: ratio between the computational costs with respect to the Rusanov scheme.}
    \label{CPU:2D_Riemann17}
\end{table}

%
%
\clearpage
\subsection{Mach 2000 astrophysical jet}\label{ssec:Mach}
Let us consider the high Mach number astrophysical jet problem \cite{HaShu,ZhangShu2010}. 
For this test in the interval $[0,1]\times[-0.25,0.25]$ we consider the following initial conditions:
\begin{equation}
\label{IC:Mach2000}
    (\rho^0,u^0,v^0,p^0) = \left\{ \begin{array}{ll}
        (5,800,0,0.4127) & \quad \text{if} \; x=0 \, \text{and} \, y\in [-0.05,0.05], \\ 
        (0.5,0,0,0.4127) & \quad \text{otherwise},
    \end{array} \right.
\end{equation} 
where $\gamma = 5/3$.

\begin{figure}[!ht]
    \centering
         \hspace{-0.76cm}
         \includegraphics[width=0.4\textwidth]{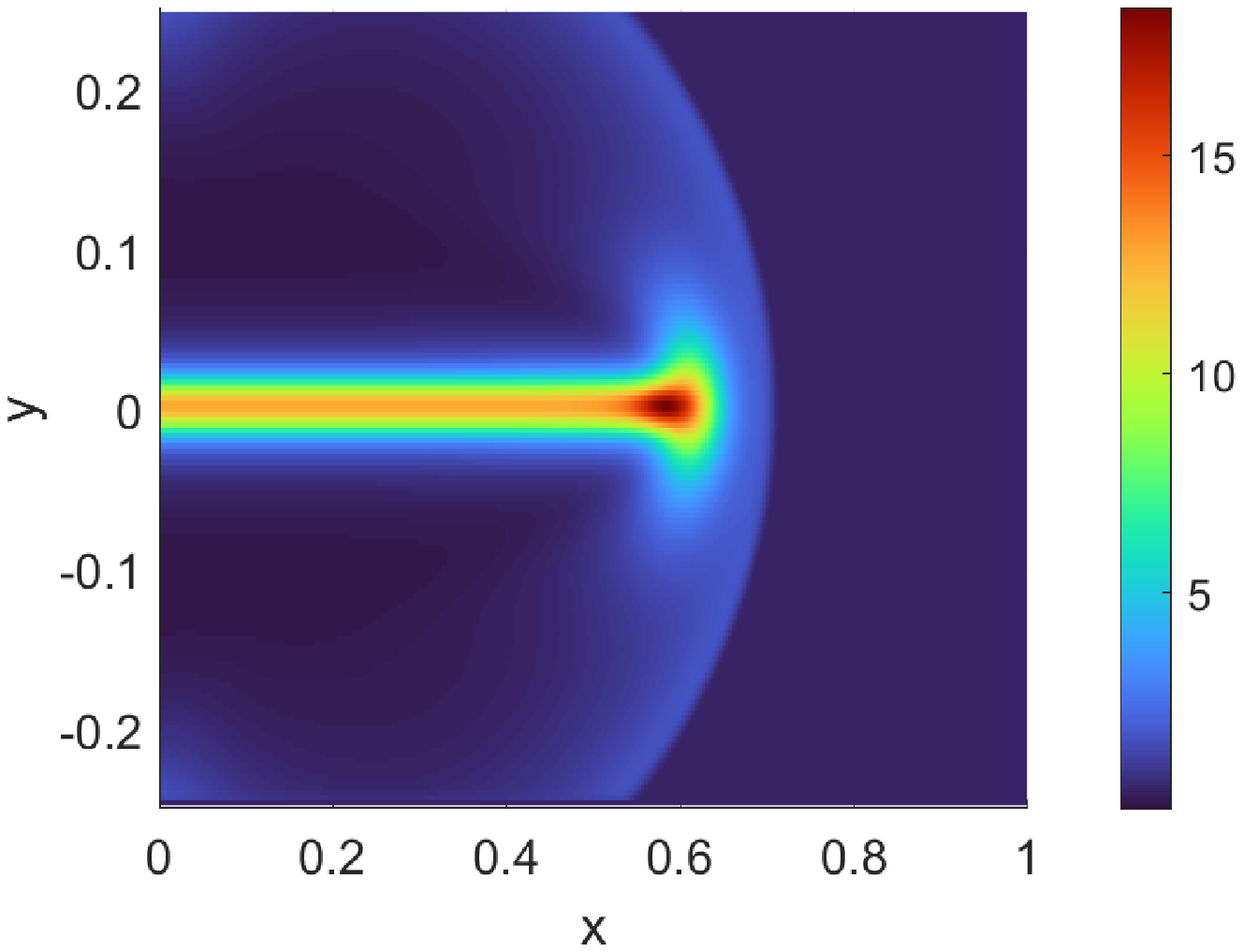}
         \hspace{-1.4cm}
         \includegraphics[width=0.4\textwidth]{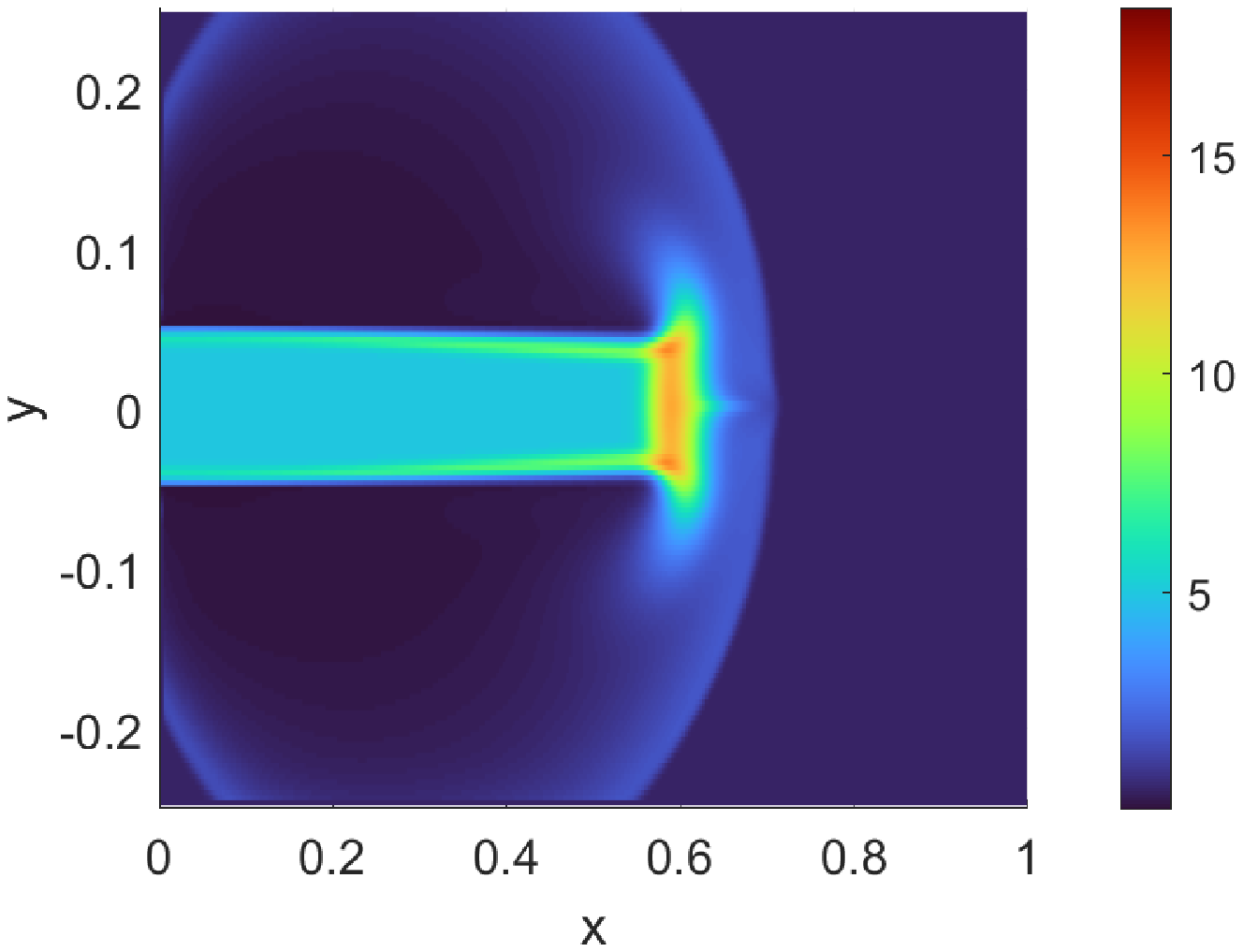}
         \hspace{-1.4cm}
         \includegraphics[width=0.4\textwidth]{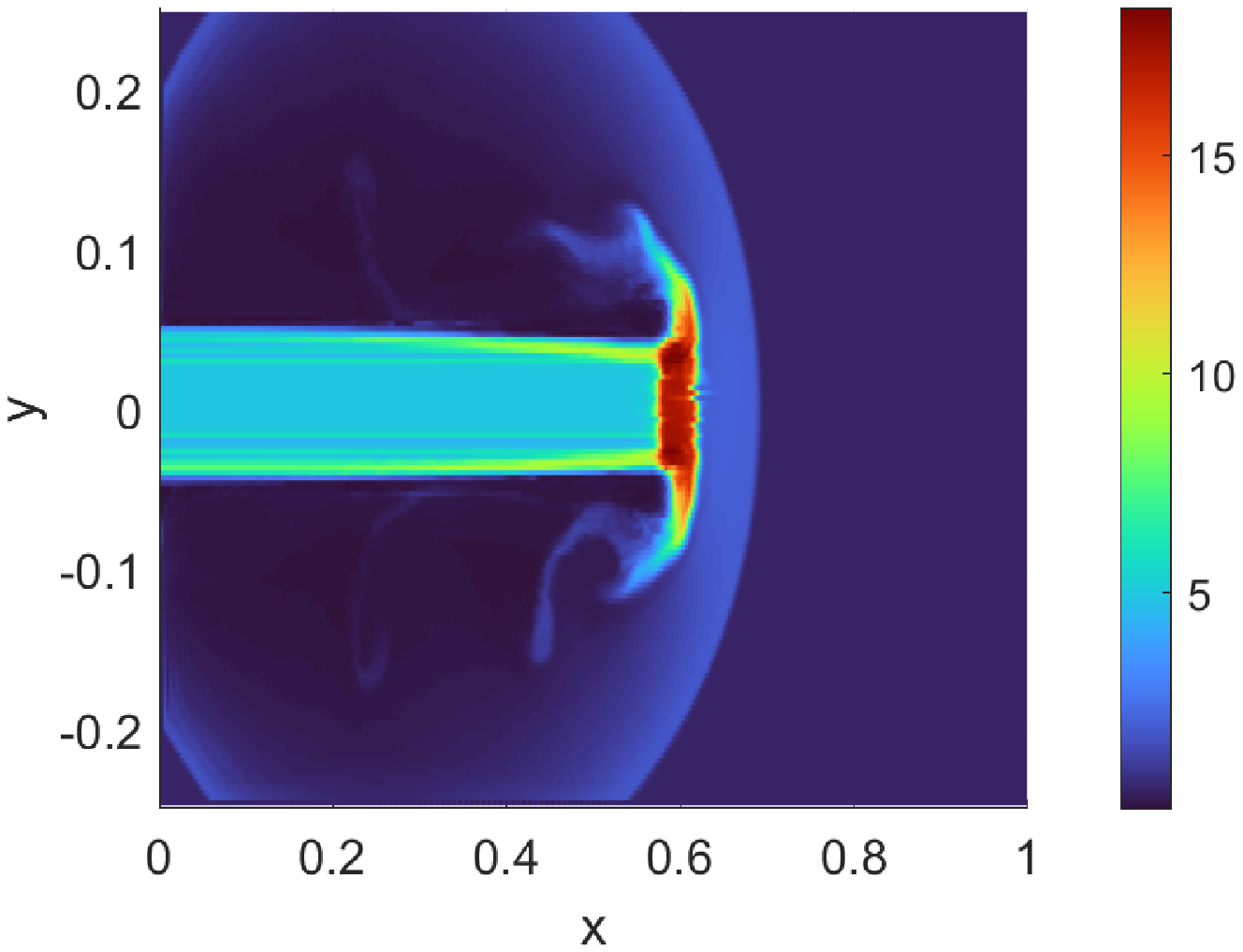}         
    \caption{Mach 2000 from section~\ref{ssec:Mach}. Numerical solution for density on the interval $[0,1]\times[-0.25,0.25]$ adopting a mesh of $300 \times 150-$ quadrangular cells and a CFL$=0.4$. The Rusanov flux (left); the HLLC (middle); and the CATMOOD6 with HLLC for the first order method (right).}
   \label{fig:Mach_2000_p}
\end{figure}

\begin{figure}[!ht]
    \centering
         \hspace{-0.76cm}
         \includegraphics[width=0.4\textwidth]{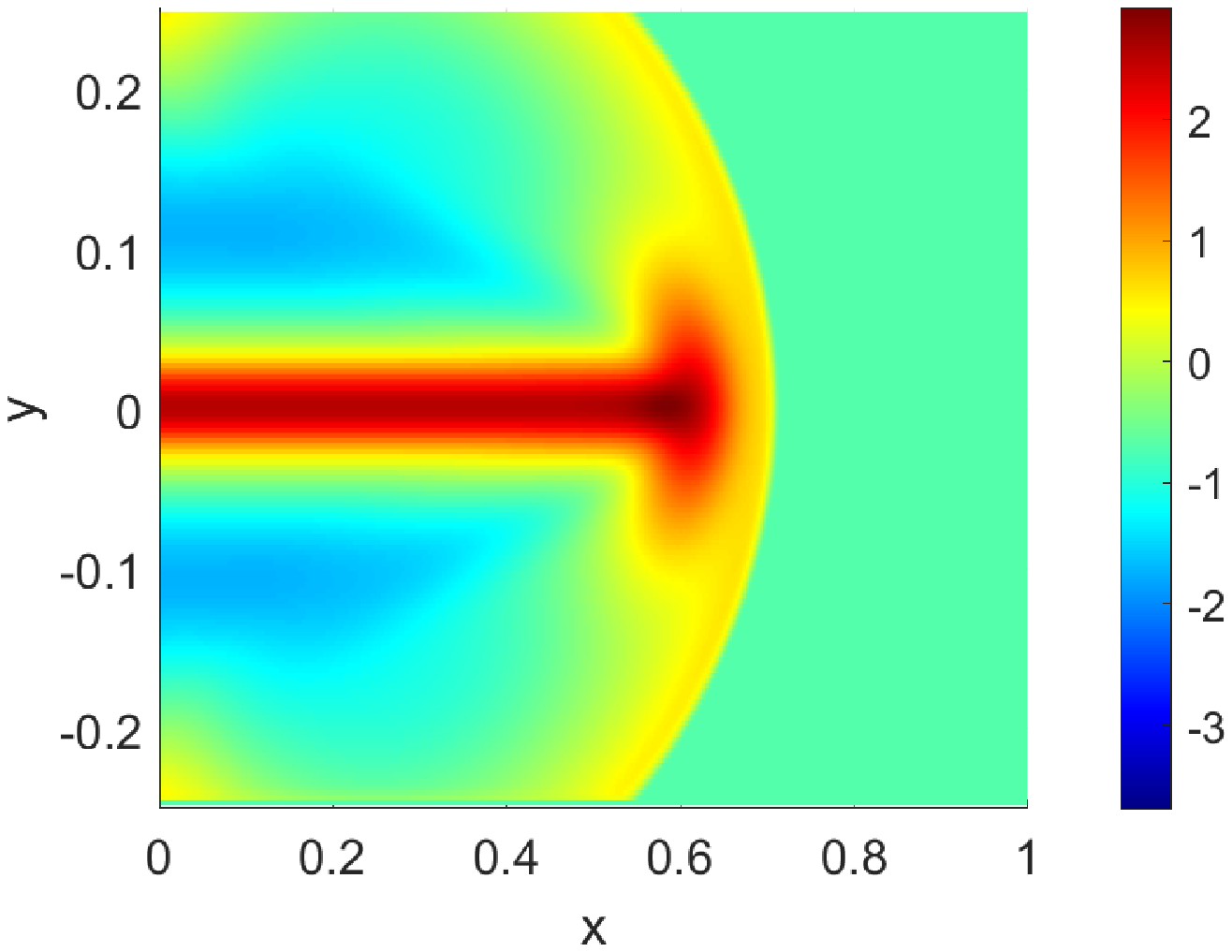}
         \hspace{-1.4cm}
         \includegraphics[width=0.4\textwidth]{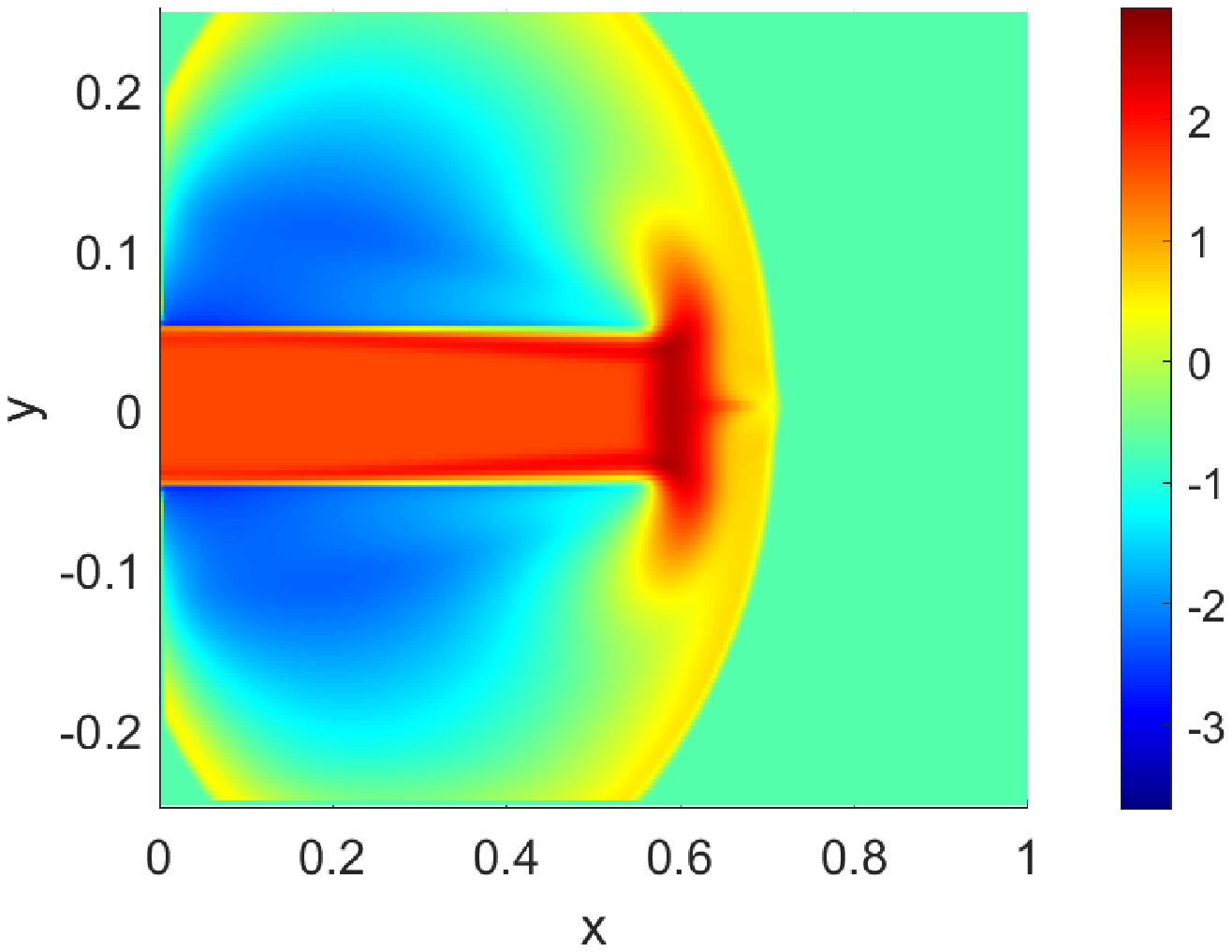}
         \hspace{-1.4cm}
         \includegraphics[width=0.4\textwidth]{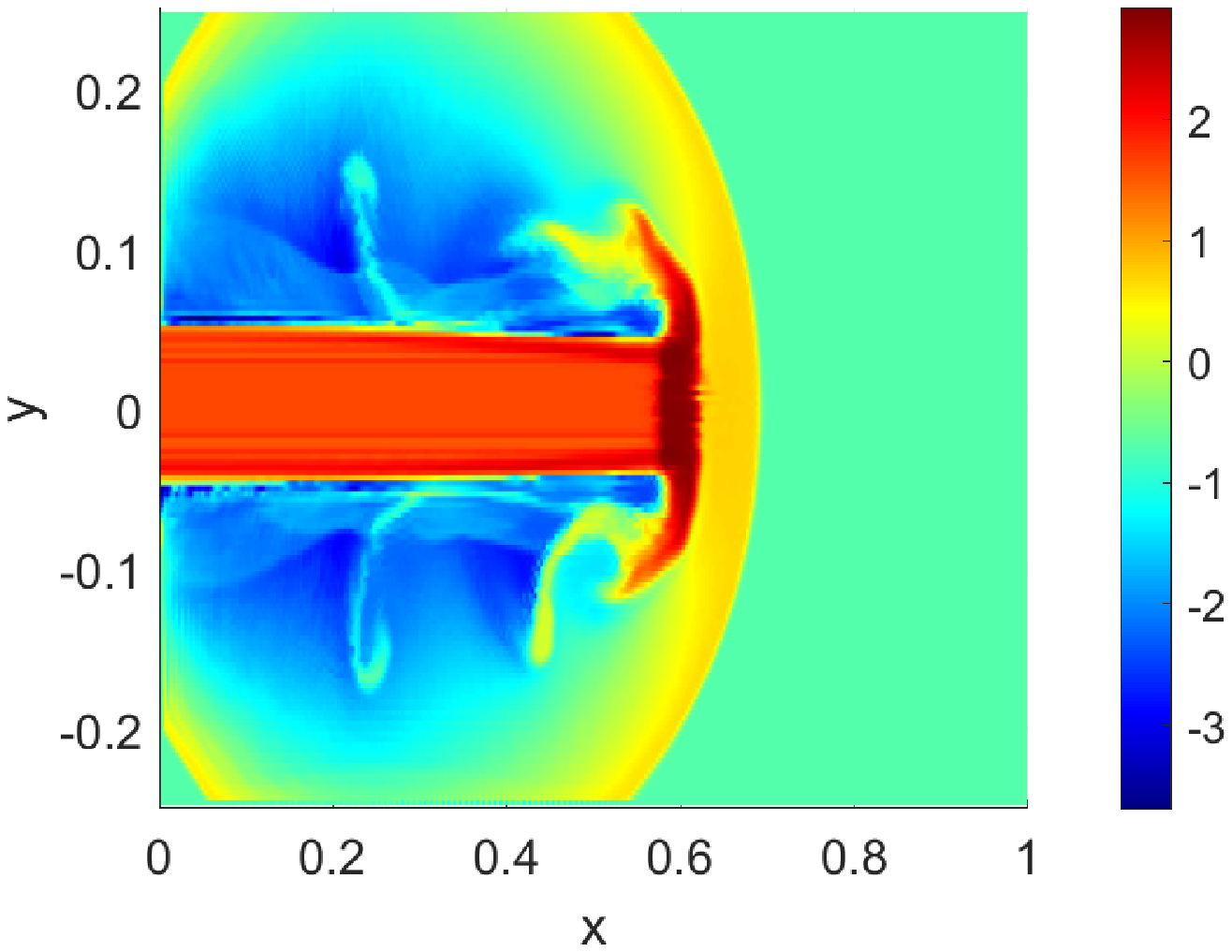}         
    \caption{Mach 2000 from section~\ref{ssec:Mach}. Numerical solution in the logarithm scale for density on the interval $[0,1]\times[-0.25,0.25]$ adopting a mesh of $300 \times 150-$cells and CFL$=0.4$. The Rusanov flux (Left); the HLLC (middle); and the CATMOOD6 with HLLC for the first order method (right).}
   \label{fig:Mach_2000_log_p}
\end{figure}
\begin{figure}[!ht]
    \centering
    \includegraphics[width=0.99\textwidth]{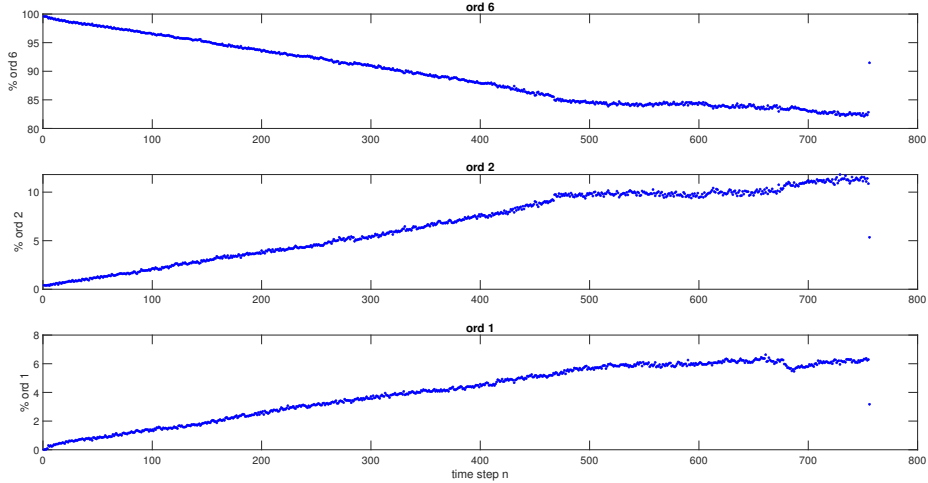}
    \vspace{-0.3cm}
    \caption{Mach 2000 from section~\ref{ssec:Mach} with initial condition \eqref{IC:Mach2000}. Percentage of cells updated by CAT6 (top), by CAT2 (center), and by HLLC (bottom).}
   \label{fig:Mach_2000_ord}
\end{figure}
Initially we have a constant ambient gas at rest except for a segment in the left boundary of the domain for which the gas is $10$ times denser and with a high $x$-component of velocity, corresponding to a Mach equals to $2000$. 
Then, we impose inflow boundary conditions in this left segment, and outflow conditions otherwise. 
As such this simulates the penetration of a 
dense jet at hyper-velocity from a portion of the left boundary. This jet generates a bow shock ahead of the jet and a complex shape of the tip of the jet.  Some reference
solutions are reported in \cite{HaShu,ZhangShu2010,BVD20} for instance.
Since the Mach number of the jet is extremely high, negative numerical pressure or density could easily appear during the computation, leading to the crash of the program. 
The simulations are run with CATMOOD6  and some 1st order schemes. The mesh is made of $300 \times 150$ uniform quadrangular cells, and a CFL equals to $0.4$ is adopted for a final time $t_\text{final} = 0.001$. 
Figures~\ref{fig:Mach_2000_p}-\ref{fig:Mach_2000_log_p} present the numerical densities, with a logarithmic scale in the latter. The simulations are performed with CATMOOD6, first order Rusanov flux and first order HLLC schemes. The first order scheme for the MOOD cascade uses the HLLC numerical flux. 
Notice that the ACAT6 scheme fails for this test due to the generation of nonphysical states (negative pressure). Obviously all unlimited CAP2P schemes for $P\geq1$ also fail.\\
The $1$st order scheme can capture the bow shock position but the tip and body of the jet are truly diffused, especially if HLLC flux is not employed. A better shape is gained with HLLC type of flux.
On the contrary CATMOOD6 is able to capture the complexity of the jet motion, its tip and the unstable lateral shear layers producing secondary waves and patterns into the post-shock region. These phenomena would be totally absent from low-accurate simulations for this mesh resolution, hence the need for truly accurate numerical methods. (They are indeed absent from the first order scheme results.)
This phenomenon is even more clear with the logarithmic scale on Figures~\ref{fig:Mach_2000_log_p}.
We would like to emphasize that CATMOOD6 does not have any issue related to positivity, because the parachute scheme is indeed one of the  the $1$st order schemes which is robust enough. As such CAT+MOOD coupling is an almost 'fail-safe' strategy. This is not an obvious property of high-order methods using  classical \apriori limiters. \\
Figure~\ref{fig:Mach_2000_ord} presents the percentage of troubled cells in CATMOOD6 as a function of time-steps.
We plot the cells updated by the schemes in the MOOD cascade, that is with the unlimited CAT6 (top), CAT2 (center) and $1$st order HLLC-base (bottom) schemes.
This is a more advanced unsteady test case compared to the 2D Riemann problems seen in section~\ref{ssec:2DRiemann} for which the solutions were self-similar.
Here, the troubled cell evolution presents two phases: first, for about $500$ time-steps the number of untroubled cells decreases linearly up to $85\%$, and, secondly, this number stagnates for about $250$ time-steps.
Interestingly, the number of cells updated with CAT2 scheme is of the order $10\%$, while truly demanding cells updated with the $1$st order scheme represent about $6\%$. 
These two schemes seem to be useful in our CATMOOD6. \\
This test case is a single example which validates the CATMOOD6 scheme for extremely demanding simulations. Here, both accuracy and robustness are required. 


\section{Conclusion and Perspectives} \label{sec:conclusion}
In this paper we have presented an \aposteriori way of 
limiting finite difference CAT2P schemes using MOOD paradigm.
We have focused on CAT2P schemes of even orders devoted to solve Euler equations on 2D Cartesian mesh with a maximal order of accuracy $6$.
CAT2P schemes are nominally of order $2P$ on smooth solutions. 
Some extra dissipative mechanism must be supplemented to deal with steep gradients or discontinuous solutions, likewise for any high order scheme.
Originally, CAT2P schemes were coupled with an automatic \apriori limiter which blends high- and low-order fluxes as in \cite{CPZMR2020}.
However the difficulty with \apriori limiters is the fact that they must (i) anticipate the occurrence of possible spurious oscillations from data at time $t^n$,
(ii) tailor the appropriate amount of dissipation to stabilize the scheme, and, (iii) ensures that a physically admissible numerical solution is always produced.
For $2$nd order schemes such \apriori limiters are available, but for higher orders they are not always performing well either on smooth solutions (lack of accuracy) or to ensure the physical admissibility (lack of robustness).
In this work we rely on an \aposteriori MOOD paradigm which computes an unlimited high-order candidate solution at time $t^{n+1}$, and, further detects troubled cells which are recomputed with a lower-order accurate scheme \cite{CDL2}.
The detection procedure marks troubled cells according to Physical, Numerical and Computer admissible criteria which are at the core of our definition of an acceptable numerical solution.
For a proof of concept, in this work, we have tested the so-called CATMOOD6 scheme based on the cascade of schemes: CAT6$\rightarrow$CAT2$\rightarrow$ 1st, where the last scheme is a first order robust scheme.  \\
We have tested this scheme on a test suite of smooth solutions (isentropic vortex), simple shock waves (cylindrical Sedov blastwave), complex self-similar solutions involving contact, shock and rarefaction interacting waves  (four state 2D Riemann problems) and, at last, on an extreme Mach $2000$ astrophysical-like jet.
CATMOOD6 has passed these tests. It has preserved the optimal accuracy on smooth parts of the solutions, an essentially-non-oscillatory behavior close to steep gradients, and, a physically valid solution. CATMOOD6 has been compared to  some $1$st order schemes to challenge its robustness and to unlimited CAT2P schemes to challenge its accuracy and cost.
We have observed that the percentage of cells updated with the 6th order accurate scheme is in the range $85-100\%$, leaving only few percents to be re-computed by the low order schemes. 
As such CATMOOD6 has a total cost $20\%$ superior to the unlimited CAT6 scheme on smooth solution, and, $10\%$ less expensive than limited ACAT6 scheme on discontinuous solution (for which ACAT6 does not fail).
From our test campaign we have observed that CATMOOD6 presents the robustness of its $1$st order scheme, and, the accuracy of its $6$th order one where appropriate. The detection procedure being able to sort out troubled cells from valid ones. \\
Concerning the perspectives, 
the extension to 3D is solely a question of implementation and testing in a parallel environment. 
A last, some extensions would be to consider different systems of PDEs, with source terms and well-balanced property, possibly stiff, or, more complex models such as Navier-Stokes equations.

\section*{Acknowledgments}
We would like to thank Simone Chiochetti for providing some reference for Sedov problem.

This research has received funding from the European Union’s NextGenerationUE – Project: Centro Nazionale HPC, Big Data e Quantum Computing, “Spoke 1” (No. CUP E63C22001000006). E. Macca was partially supported by GNCS No. CUP E53C22001930001 Research Project “Metodi numerici per problemi differenziali multiscala: schemi di alto ordine, ottimizzazione, controllo”. E. Macca and G. Russo are members of the INdAM Research group GNCS. The  research  of C. Parés  is part of the research project PDC2022-133663-C21 funded by MCIN/AEI/10.13039/501100011033 end the European Union \lq Next GenerationEU/PRTR\rq and it was  partially  supported  by  the Regional Government of Andalusia through the research group FQM-216.

\bibliographystyle{plain}
\bibliography{./biblio}

\section{Appendices} \label{sec:appendix}
\subsection{Fourth order version -- CAT4}
\label{sec:CAT4}
\begin{figure}[!ht]
    	\centering\includegraphics[scale=0.7]{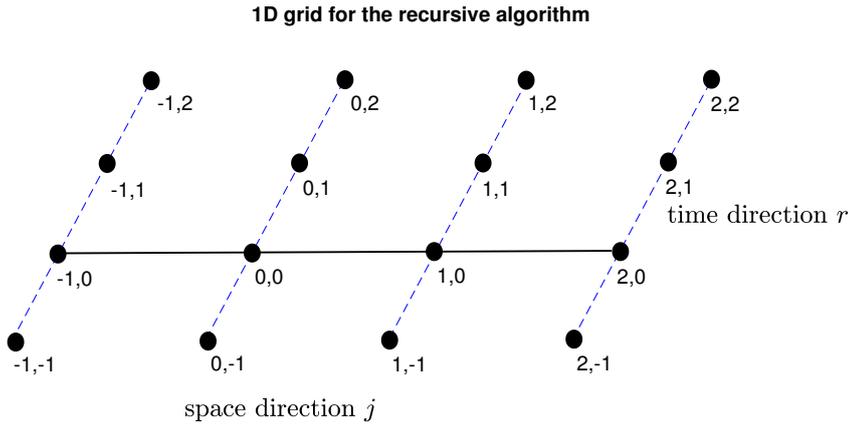}
    	\vspace{-1.6 cm}
    	\caption{Local space-time  grid where approximations of $U$ are computed to calculate $F^P_{i+1/2}$ with $P = 2$. For simplicity a pair $j,r$ represents the point $(x_{i+j}, t_{n+r})$. Taylor expansions in time are used to obtain these approximations following the blue lines. These Taylor expansions are centered in the points lying on the black line.}
    	\label{1D_grid}
\end{figure}

For a complete presentation of the general procedure to compute the $k-$th time derivative of flux, we detail hereafter the CAT4 scheme, that is with $P=2$.
Unlike CAT2 description, CAT4 scheme's description presents the high-order CAT automatism by introducing in detail the development of the iterative procedure used to compute the time derivatives of the flux.
With this in mind, see Figure \ref{1D_grid} to have a graphic idea of the necessary stencil.

The iterative CAT4 algorithm proceeds with the following steps to compute the flux $F_{i+\ha}^P$:
\begin{itemize}
    \item[Step 1:] Compute $f^{(0)}_{i+\ha}$ adopting the interpolation formula on stencil $\mathcal{S}^2_{i+\ha}=\left\{ u_{i-1}^n, u_i^n, u_{i+1}^n, u_{i+2}^n \right\}$, so that $-P+1=-1\leq p \leq 2=P$, as:
    $$ f^{(0)}_{i+\ha} = \sum_{p=-1}^{P=2}\gamma_{2,p}^{0,\ha} \, f_{i+p}^n, $$
    where $\gamma_{2,p}^{0,\ha}$ is the interpolation coefficient for $P=2$ and label $p$.
    \item[Step 2:] Compute the first time derivative of $u$ at time $t^n$ at position $x_{i+p}$ with $p = -1,\ldots,2$ through the numerical compact Cauchy-Kovaleskaia identity
    \begin{equation}
        \label{CK}
        \partial^k_tu = -\partial^k_xf^{(k-1)}(u), \quad k = 1,2,\ldots  
    \end{equation} 
    For $k=1$ we have:
    $$ u_{i,j}^{(1)} = -\frac{1}{\Delta x}\sum_{j = -1}^2\gamma_{2,j}^{1,0}f_{i+j}^n . $$
    \item[Step 3:] Compute the Taylor expansion of $u$ in time truncated at first term at time $t^{n+r}$ with $r=-1,\ldots,2$ for each position $x_{i+j}$ with $j=-1,\ldots,2$ as:
    $$ u_{i,j}^{1,n+r} = u_{i+j}^n + r\Delta t \,u_{i,j}^{(1)} . $$
    \item[Step 4:] Compute the first time derivative of $f$ at time $t^n$ at position $x_{i+j}$ with $j = -1,\ldots,2$ using the four fluxes $f_{i,j}^{1,n+r}=f\left( u_{i,j}^{1,n+r} \right)$ as:
    $$ f_{i,j}^{(1)} = \frac{1}{\Delta t}\sum_{r=-1}^2 \gamma_{2,r}^{1,0} \, f_{i,j}^{1,n+r} .$$
    \item[Step 5:] Compute $f^{(1)}_{i+\ha}$ adopting the interpolation formula on stencil $\mathcal{S}^2_{i+\ha}$ as:
    $$ f^{(1)}_{i+\ha} = \sum_{j=-1}^2\gamma_{2,j}^{0,\ha}f_{i,j}^{(1)}; $$
    \item[Step 6:] Compute the second time derivative of $u$ at time $t^n$ and position $x_{i+j}$ with $j=\nobreak -1,\ldots,2$ from the first time derivatives of $f$ using the \eqref{CK} as:
    $$ u_{i,j}^{(2)} = -\frac{1}{\Delta x}\sum_{s = -1}^2\gamma_{2,s}^{1,j}f_{i,s}^{(1)}; $$
    \item[Step 7:] Compute the Taylor expansion in time truncated at second term at time $t^{n+r}$ with $r=-1,\ldots,2$ for each position $x_{i+j}$ with $j=-1,\ldots,2$ as:
    $$ u_{i,j}^{2,n+r} = u_{i+j}^n + r\Delta t u_{i,j}^{(1)} + \frac{(r\Delta t)^2}{2}u^{(2)}_{i,j}; $$
    \item[Step 8:] Compute the second time derivative of $f$ at time $t^n$ at position $x_{i+j}$ with $j = -1,\ldots,2$ using the four fluxes $f_{i,j}^{2,n+r}$ as:
    $$ f_{i,j}^{(2)} = \frac{1}{\Delta t^2}\sum_{r=-1}^2\gamma_{2,r}^{2,0}f_{i,j}^{2,n+r} . $$
    \item[Step 9:] Compute $f^{(2)}_{i+\ha}$ adopting the interpolation formula on stencil $\mathcal{S}^2_{i+\ha}$ as:
    $$ f^{(2)}_{i+\ha} = \sum_{j=-1}^2\gamma_{2,j}^{0,\ha}f_{i,j}^{(2)} .$$
    \item[Step 10:] Compute the third time derivative of $U$ at time $t_n$ and position $x_{i+j}$ with $j=-1,\ldots,2$ from the second time derivatives of $f$ using the \eqref{CK} as:
    $$ u_{i,j}^{(3)} = -\frac{1}{\Delta x}\sum_{s = -1}^2\gamma_{2,s}^{1,j}f_{i,s}^{(2)}; $$
    \item[Step 11:] Compute the Taylor expansion in time truncated at third term at time $t^{n+r}$ with $r=-1,\ldots,2$ for each position $x_{i+j}$ with $j=-1,\ldots,2$ as:
    $$ u_{i,j}^{3,n+r} = u_{i+j}^n + r\Delta t u_{i,j}^{(1)} + \frac{(r\Delta t)^2}{2}u^{(2)}_{i,j} + \frac{(r\Delta t)^3}{6}u^{(3)}_{i,j} . $$
    \item[Step 12:] Compute the third and last time derivative of $f$ at time $t^n$ at position $x_{i+j}$ with $j = -1,\ldots,2$ using the four fluxes $f_{i,j}^{3,n+r}$ as:
    $$ f_{i,j}^{(3)} = \frac{1}{\Delta t^3}\sum_{r=-1}^2\gamma_{2,r}^{3,0}f_{i,j}^{3,n+r} . $$
    \item[Step 13:] Compute $f^{(3)}_{i+\ha}$ adopting the interpolation formula on stencil $\mathcal{S}^2_{i+\ha}$ as:
    $$ f^{(3)}_{i+\ha} = \sum_{j=-1}^2\gamma_{2,j}^{0,\ha}f_{i,j}^{(3)}; $$
    \item[Step 14:] Deduce $F^2_{i+\ha}$ from \eqref{FP} as:
    $$ F^2_{i+\ha} = f_{i+\ha}^{(0)} + \Delta t \, f_{i+\ha}^{(1)} + \frac{\Delta t^2}{2} \, f_{i+\ha}^{(2)} + \frac{\Delta t^3}{6} \, f_{i+\ha}^{(3)}   $$
\end{itemize}

\subsection{Computational complexity}
\label{sec:comp_compl}
In this section we focus on the local computational complexity for the CAT$2P$ scheme applied to scalar case. With this in mind, observing that the size of the stencil $S_{i+\ha}^P$ is $2P,$ we can divide the CAT algorithm into four parts and notice that:
\begin{enumerate}
    \item $f_{i,j}^n = f(u_{i+j}^n)$ for all $j=-P+1,\ldots,P.$
    
    \item For $k=1,\ldots,2P-1,$ the procedure to compute $u_{i,j}^{(k)},$ $u_{i,j}^{k,n+r},$ $f(u_{i,j}^{k,n+r})$ and $f_{i,j}^{(k)}$ does not depend explicitly on $k$ (see expressions (\ref{c2.a}-\ref{c2.d}) below). 
    
    \item For $k=0,\ldots,2P-1,$ the computational formula for the approximation of the $(k-1)$-th time derivative of the flux at position $x = x_{i+1/2},$   $f_{i+1/2}^{(k-1)} =  \mathcal{A}^{0,1/2}_{P}\left(f_{i,*}^{(k-1)},\Delta x\right) = \sum_{j=-P+1}^P \gamma^{0,\ha}_{P,j}f_{i,j}^{(k-1)},$ does not depend explicitly on $k,$ since the interpolatory formula is invariant for each $k.$ 
    
    \item Compute  $F^P_{i+\ha},$ $u_i^{n+1}$ and precomputed constants such as $\Delta t^k$ or $\Delta t^k/k!$, etc.
\end{enumerate} 
For this reason, the computational complexity is so structured:
\begin{enumerate}
    \item \label{aa} $f_{i,j}^n = f(u_{i+j}^n)$ for all $j=-P+1,\ldots,P$. $1$ evaluation multiplied by the size of $u$\footnote{$M\times N$ for a discrete system with $M-$variables and $N$ cells.}.
    
    \item \label{c2} For $k=1,\ldots,2P-1$
    \begin{enumerate}
        \item \label{c2.a} $u_{i,j}^{(k)} = -\frac{1}{\Delta x}\sum_{s = -P+1}^{P}\gamma_{P,s}^{1,j}f_{i,s}^{(k-1)}.$ \;
        $2P(2P+1)$ flop for each $k$, $2P$ related to $j$ and $2P+1$ related to the summation.
        
        \item \label{c2.b} $ u_{i,j}^{k,n+r} = u^n_{i+j} +\sum_{m=1}^k\frac{(r\Delta t)^m}{m!} u^{(m)}_{i,j} $ under the assumption $ u_{i,j}^{k,n+r} = u^{k-1,n+r}_{i,j} +\frac{c_r^k}{k!} u^{(k)}_{i,j}$ where  $c_r=r\Delta t$ and $u^{0,n+r}_{i,j} = u^{n}_{i+j}.$ \; $(2P)(2P-1)$ flop for each $k$ under the assumption that $c_r^k/k!$ are precomputed and, when $r=0,$ $u_{i,j}^{k,n} \equiv u_{i+j}^n$.
        
        \item \label{c2.c} $f_{i,j}^{k,n+r}=f(u_{i,j}^{k,n+r})$. \; $(2P)(2P - 1)$ evaluations for each $k$\footnote{The evaluations of $f(u_{i,j}^{k,n+r})$ are much expensive than $f(u_{i+j}^n).$ Indeed, for each $k$ the size of $f(u_{i,j}^{k,n+r})$ is $M\times N,$ the size of $f(u_{i,j}^{k,n+r})$ is $M\times N \times 2P \times 2P-1,$ respectively for $j$ and $r$. For the scalar case $M=1.$}, $(2P)$ related to $j$ and $(2P-1)$ related to $r.$ When $r = 0,$ $f(u_{i,j}^{k,n+r}) = f(u_{i+j}^{n})$ and the computational complexity relapses to expression \eqref{aa}. 
        
        \item \label{c2.d} $f^{(k)}_{i,j} = \mathcal{A}^{k,j}_{P}\left(f_{i,j}^{k,*},\Delta t\right) = \frac{1}{\Delta t^k}\sum_{r=-P+1}^P\gamma_{P,r}^{k,0}f_{i,j}^{k,n+r} .$ \;
        $2P(2P+1)$ flop for each $k$ under the assumption that $\Delta t^k$ are precomputed.
    \end{enumerate}
    
    \item For $k=1,\ldots,2P,$ $f_{i+\ha}^{(k-1)}.$ \; $2P$ flop.
    
    \item The computational complexity of $F^P_{i+\ha},$  $u_i^{n+1}$ and precomputed terms 
    \begin{enumerate}
        \item $ F^P_{i+1/2} = \sum_{k=1}^{2P}\frac{\Delta t^{k-1}}{k!}f_{i+1/2}^{(k-1)} $ under the assumption that $\Delta t^{k-1}/k!$ are precomputed. \; $(2P-1)$ flop since for $k=1$ there are no products.
        
        \item $ u_i^{n+1} = u_i^n - \frac{\Delta t}{\Delta x}\Bigl(F^P_{i+1/2} - F^P_{i-1/2}\Bigr).$ \; $2$ flop.
        
        \item for $k=1,\ldots,2P-1$ $\Delta t^k.$ \; $k-1$ flop.
        
        \item for $k=1,\ldots,2P-1$ and for $r = -P+1,\ldots,P,$ with $r\neq0,$ $c_r^k/k!$ where $c_r =r\Delta t.$ \; $(k+1)(2P-1)$ flop.
        
        \item for $k=1,\ldots,2P-1$ $\Delta t^{k-1}/k!.$ \; $1$ flop.        
    \end{enumerate} 
\end{enumerate}
\begin{enumerate}
    \item $n_{\rm func}\leftarrow1$ evaluation.
    \item\begin{enumerate}
        \item $n_{\rm flop}\leftarrow (2P)^3 - (2P)$ flop;
        \item $n_{\rm flop}\leftarrow n_{\rm flop} + (2P)^3 - 2(2P)^2 + (2P)$ flop;
        \item $n_{\rm func} \leftarrow n_{\rm func} +(2P)^3 - 2(2P)^2 + (2P)$ evaluations;
        \item $n_{\rm flop}\leftarrow n_{\rm flop} + (2P)^3 - (2P)$ flop.
    \end{enumerate}
    \item $n_{\rm flop}\leftarrow n_{\rm flop} + (2P)^2$ flop.
    \item\begin{enumerate}
        \item $n_{\rm flop}\leftarrow n_{\rm flop} + (2P)-1$ flop;
        \item $n_{\rm flop}\leftarrow n_{\rm flop} + 2$ flop;
        \item $n_{\rm flop}\leftarrow n_{\rm flop} + 0.5(2P)^2 - 1.5(2P)$ flop;
        \item $n_{\rm flop}\leftarrow n_{\rm flop} + 0.5(2P^3) - 0.5(2P)$ flop;
        \item $n_{\rm flop}\leftarrow n_{\rm flop} + (2P) - 1$ flop.
    \end{enumerate}
\end{enumerate}
Finally, the operation count per cell per time step for CAT$2P,$ $P>1,$ applied to the scalar case is $3.5(2P)^3 - 1.5(2P)^2 + (2P)$ flop plus $(2P)^3 - 2(2P)^2 + (2P) + 1$ function evaluations. For CAT2 this gives $24$ flop and $3$ function evaluations, while for CAT4 we get $204$ flop and $37$ function evaluations.

\begin{table}[!ht]
    \centering
    \numerikNine
    \begin{tabular}{|cc||c|c|c|}
    \hline
    Scheme & $P$ & Cost per flux & Ratio & Ratio \\
    & & {\footnotesize $C(P)=3.5(2P)^3 - 1.5(2P)^2 + (2P)$} & {\footnotesize$C(P+1)/C(P)$} & {\footnotesize CAT2P/CAT2} \\
    \hline
    CAT2 & $1$ & 29   & 1 & 1 \\
    CAT4 & $2$ & 299  & 10.31 & 10.31\\
    CAT6 & $3$ & 1097 &  3.67 & 37.83 \\
    CAT8 & $4$ & 2711 &  2.47 & 93.48 \\
    CAT10 &$5$ & 5429 &  2.00 &187.2 \\
    $\vdots$ & $\vdots$ & $\vdots$ & $\vdots$ & $\vdots$ \\
    CAT2P & $P$ & $C(P)$ & 1.00 & $\infty$ \\ 
    \hline
    \end{tabular}
    \caption{Cost of the CAT schemes in term of number of operations to update one interface flux as a function of $P$. }
    \label{tab:Cost}
\end{table}

\end{document}